\documentclass[12pt]{article}
\usepackage{amsmath}
\usepackage{amssymb}
\usepackage[dvips]{graphicx}
\usepackage{epsfig}
\usepackage{subfigure}
\usepackage{psfrag}
\usepackage{vector}
\def\eps{\varepsilon}
\font\tencmmib=cmmib10 \skewchar\tencmmib '60
\newfam\cmmibfam
\textfont\cmmibfam=\tencmmib

\def\bbox{\quad\hbox{\vrule \vbox{\hrule \vskip2pt \hbox{\hskip2pt
\vbox{\hsize=1pt}\hskip2pt} \vskip2pt\hrule}\vrule}}
\def\lessim{\ \lower4pt\hbox{$
\buildrel{\displaystyle <}\over\sim$}\ }
\def\gessim{\ \lower4pt\hbox{$\buildrel{\displaystyle >}
\over\sim$}\ }
\def\P{{\cal P}}

\def\R{{\cal R}}

\def\V{{\cal V}}

\def\eps{{\varepsilon}}
\def\ch{{\mbox{ch}}}
\def\sh{{\mbox{sh}}}

\def\zz{{\mathfrak{z}}}
\def\nn{{\mathfrak{n}}}
\def\m{{\mathfrak{m}}}

\def\RS{{\mbox{RS}}}

\def\PM{{\mbox{PM}}}
\def\card{{\mbox{\rm{card}}}}

\def\la{{\Bigl\langle}}
\def\ra{{\Bigr\rangle}}

\def\qed{\hfill\break\rightline{$\bbox$}}
\newcommand{\e}{\mathbb{E}}
\newcommand{\p}{\mathbb{P}}
\newcommand{\Reals}{\mathbb{R}}

\newcommand{\vsi}{{\vec{\sigma}}}
\newcommand{\vrho}{{\vec{\rho}}}

\newtheorem{proposition}{Proposition}
\newtheorem{lemma}{Lemma}
\newtheorem{theorem}{Theorem}
\newtheorem{corollary}{Corollary}
\makeatletter
\@addtoreset{equation}{section}

\makeatother

\font\tencmmib=cmmib10 \skewchar\tencmmib '60
\newfam\cmmibfam
\textfont\cmmibfam=\tencmmib

\def\bbox{\quad\hbox{\vrule \vbox{\hrule \vskip2pt \hbox{\hskip2pt
\vbox{\hsize=1pt}\hskip2pt} \vskip2pt\hrule}\vrule}}
\def\lessim{\ \lower4pt\hbox{$
\buildrel{\displaystyle <}\over\sim$}\ }
\def\gessim{\ \lower4pt\hbox{$\buildrel{\displaystyle >}
\over\sim$}\ }

\def\eps{\varepsilon}

\def\go0{\to 0}

\def\la{\langle}

\def\leftitem#1{\item{\hbox to\parindent{\enspace#1\hfill}}}

\def\qed{\hfill\break\rightline{$\bbox$}}

\def\ra{\rangle}

\def\sg{\sigma}

\def\sg2{\sigma^2}

\def\__{_{\infty}}
\begin{document}

\title{
Free energy in the generalized Sherrington-Kirkpatrick 
mean field model.}

\author{ 
Dmitry Panchenko\thanks{Department of Mathematics, Massachusetts Institute
of Technology, 77 Massachusetts Ave, Cambridge, MA 02139
email: panchenk@math.mit.edu}\\
{\it Department of Mathematics}\\
{\it Massachusetts Institute of Technology}\\
}
\date{}

\maketitle

\begin{abstract}
In \cite{T-P} Michel Talagrand gave a rigorous proof of the Parisi 
formula in the classical Sherrington-Kirkpatrick (SK) model. 
In this paper we build upon the methodology developed in \cite{T-P}
and extend Talagrand's result to the class of SK type models in which 
the spins have arbitrary prior distribution on a bounded subset 
of the real line. 

\end{abstract}
\vspace{0.5cm}


\section{Introduction and main results.}

In \cite{T-P} Michel Talagrand invented a rigorous proof of
the Parisi formula for the free energy in the Sherrington-Kirkpatrick model
\cite{SherK}.
The methodology developed by Talagrand was based upon a deep 
extension of Guerra's interpolation method in \cite{Guerra} 
to coupled systems of spins which provided necessary control 
of the remainder terms in Guerra's interpolation. 
The same methodology was successfully used in \cite{T-Sph}
to compute the free energy in the spherical model and
in the present paper we will utilize it in the setting of 
a generalized Sherrington-Kirkpatrick model in which the prior
distribution of the spins is given by an arbitrary probability 
measure with bounded support on the real line. 

Let us start by introducing all necessary notations and definitions.
Consider a {\it bounded} set $\Sigma\subseteq\Reals$
and a probability measure $\nu$ on the Borel $\sigma$-algebra
on $\Sigma.$
Given $N\geq 1,$ consider a product space $(\Sigma^N,\nu^N)$ 
which will be called the space of configurations. 
A configuration $\vsi\in\Sigma^N$ is a vector 
$(\sigma_1,\ldots,\sigma_N)$ of spins $\sigma_i$ 
that take values in $\Sigma.$ 
For simplicity of notations we will omit index $N$ in 
$\nu^N$ since it will always be clear from the context 
whether we consider measure $\nu$ on $\Sigma$ or the product
measure on $\Sigma^N.$
For each $N$ we consider a Hamiltonian $H_N(\vsi)$ on 
$\Sigma^N$ that is a Gaussian process indexed by $\vsi\in \Sigma^N.$ 
We will assume that $H_N(\vsi)$ is jointly measurable in $(\vsi,\vec{g}),$
where $\vec{g}$ is the generic point of the underlying probability 
space on which the Gaussian process $H_N(\vsi)$ is defined.
We assume that for a certain sequence $c(N)\to 0$
and a certain function $\xi:\Reals\to\Reals$ we have
\begin{equation}
\forall \vsi^1, \vsi^2 \in \Sigma^N ,\,\,\,\,\,
\Bigl|\frac{1}{N}\e H_N(\vsi^1) H_N(\vsi^2) - \xi(R_{1,2})\Bigr|
\leq c(N),
\label{correlation}
\end{equation}
where
\begin{equation}
R_{1,2}=\frac{1}{N}\sum_{i\leq N}\sigma_i^1\sigma_i^2
\label{overlap}
\end{equation}
is called {\it the overlap} of the configurations $\vsi^1, \vsi^2.$
We will assume that $\xi(x)$ is three times continuously differentiable
and satisfies the following conditions
\begin{equation}
\xi(0)=0,\, \xi(x)=\xi(-x),\,\xi''(x)>0 \mbox{ if } x>0.
\label{xi}
\end{equation}
We will denote {\it the self-overlap} of $\vsi$ by
$$
R_{1,1}=\frac{1}{N}\sum_{i\leq N} (\sigma_i)^2.
$$
One defines the Gibbs measure $G_N$ on $\Sigma^N$ by
\begin{equation}
d G_N(\vsi)=\frac{1}{Z_N}
\exp H_N(\vsi) d\nu(\vsi)
\label{Gibbs}
\end{equation}
where the normalizing factor
$$
Z_N=\int_{\Sigma^N}
\exp H_N(\vsi) d\nu(\vsi)
$$
is called the {\it partition function.} 
This definition of the Gibbs measure also includes the
case of models with external field because given a
measurable function $h(\sigma)$ on $\Sigma$ and the
Gibbs measure defined by
$$
d G_N(\vsi)=\frac{1}{Z_N}
\exp\Bigl(H_N(\vsi) +\sum_{i\leq N}h(\sigma_i)\Bigr) d\nu(\vsi)
$$
we can simply make the change of measure 
$d\nu'(\sigma) \sim \exp h(\sigma) d\nu(\sigma)$
to represent this Gibbs measure as (\ref{Gibbs}).
The only assumption that we need to make on $h(\sigma)$ is that
$$
\int\exp h(\sigma)d\nu(\sigma)<\infty
$$
which holds, for example, when $h(\sigma)$ is uniformly bounded.
This general model includes the original 
Sherrington-Kirkpatrick (SK) model in \cite{SherK} and the
Ghatak-Sherrington (GS) model in \cite{GS} that will be 
considered in more detail in Section \ref{GSsec}. 
In both cases the Hamiltonian $H_N(\vsi)$ is given by
\begin{equation}
H_N(\vsi)=\frac{\beta}{\sqrt{N}}\sum_{i<j}g_{ij}\sigma_i\sigma_j 
\label{Hamorig}
\end{equation}
for some $\beta>0$ and i.i.d. Gaussian r.v. $g_{ij},$
the set $\Sigma$ is equal to $\{-1,+1\}$ in the SK model and 
$\{0,\pm 1,\ldots,\pm S\}$ for some integer $S$ in the GS model 
and in both cases the measure $\nu$ is uniform on $\Sigma.$ 
In the case of the SK model the function $h(\sigma)$ is given by
$h(\sigma)=h\sigma$ with the {\it external field} parameter $h\in \Reals,$
and in the case of the GS model it is given by
$h(\sigma)=h\sigma^2$ with the {\it crystal field} parameter $h\in\Reals.$
Let us define
\begin{equation}
F_N=\frac{1}{N}\e \log Z_N
=
\frac{1}{N}\e \log
\int_{\Sigma^N}
\exp H_N(\vsi)
d\nu(\vsi),
\label{FE1}
\end{equation}
which (usually, with the factor $-\beta^{-1}$ which we omit for simplicity
of notations) is called the {\it free energy} of the system $(\Sigma^N, G_N).$ 
The main goal of this paper is to find the limit
$\lim_{N\to\infty} F_N.$
It will soon become clear that the main difference of the above model
from the classical SK model lies in the fact that 
in the classical model the length of any configuration 
$\vsi\in\{-1,+1\}^N$ was constant, $|\vsi| = \sqrt{N},$ 
which is not always true here.
In general, if $\Sigma$ is not of the type $\{-a,+a\}$ 
for some $a\in\Reals$ then the length of the configuration 
or self-overlap $R_{1,1}=|\vsi|^2/N$ will become variable. 
As a result, in order to make the methodology of Guerra and Talagrand work,
we will first have to compute the {\it local} free energy 
of the set of configurations with constrained self-overlap. 

Let us now describe the analogue of the Parisi formula that
gives the limit of (\ref{FE1}).
Let $[d,D]$ be the smallest interval such that 
\begin{equation}
\nu(\{\sigma : \sigma^2\in[d,D]\})=1.
\label{diam1}
\end{equation} 
In other words, $d\leq \sigma^2\leq D$ with probability one and
$\sigma^2$ can take values arbitrarily close to $d$ and $D$
with positive probability. 
From now on we will simply say that $d\leq \sigma^2\leq D$ for all
$\sigma\in\Sigma.$
Let us consider $u\in [d,D]$ and a sequence $(\eps_N)$ such that
$\eps_N>0$ and $\lim_{N\to\infty}\eps_N = 0,$ 
and consider a sequence of sets
\begin{equation}
U_N=
\Bigl\{\vsi\in\Sigma^N : R_{1,1}\in
[u-\eps_N,u+\eps_N]\Bigr\}.
\label{const}
\end{equation}
We define
\begin{equation}
F_N(u,\eps_N)=
\frac{1}{N}\e \log Z_N(u,\eps_N), 
\label{FE2}
\end{equation}
where
$$
Z_N(u,\eps_N) = \int_{U_N}\exp H_N(\vsi) d\nu(\vsi).
$$
$F_N(u,\eps_N)$ is the free energy of the subset of configurations in (\ref{const}).  
We will first compute $\lim_{N\to\infty}F_N(u,\eps_N)$ for some sequence $(\eps_N)$
for each $u\in[d,D].$
Consider an integer $k\geq 1,$ numbers
\begin{equation}
0=m_0\leq m_1\leq \ldots \leq m_{k-1}\leq m_k=1
\label{m}
\end{equation}
and, given $u\in[d,D]$,
\begin{equation}
0=q_0\leq q_1\leq \ldots \leq q_{k}\leq q_{k+1} = u.
\label{q}
\end{equation}
We will write $\vec{m}=(m_0,\ldots,m_k)$ and $\vec{q}=(q_0,\ldots,q_{k+1}).$
Consider independent centered Gaussian r.v. 
$z_p$ for $0\leq p\leq k$ with 
\begin{equation}
\e z_p^2 = \xi'(q_{p+1}) - \xi'(q_p).
\label{z}
\end{equation}
Given $\lambda\in\Reals,$ we define the r.v.
\begin{equation}
X_{k+1}=\log 
\int_{\Sigma}
\exp\Bigl(\sigma\sum_{0\leq p\leq k} z_p 
+ \lambda \sigma^2 \Bigr)d\nu(\sigma)
\label{Xlast}
\end{equation}
and, recursively for $l\geq 0,$ define
\begin{equation}
X_l=\frac{1}{m_l}\log \e_l\exp m_l X_{l+1},
\label{Xl}
\end{equation}
where $\e_l$ denotes the expectation in the r.v. $(z_p)_{p\geq l}.$
When $m_l=0$ this means $X_l=\e_l X_{l+1}.$ 
Clearly, $X_0=X_0(\vec{m},\vec{q},\lambda)$ 
is a non-random function of the parameters
$\vec{m},\vec{q}$ and $\lambda.$ 
Whenever it does not create ambiguity
we will keep this dependence implicit.
Let us note that $X_0$ also depends on $u$ 
through $\vec{q}$ since in (\ref{q}) we have
$q_{k+1}=u.$
Let
\begin{equation}
\P_k(\vec{m},\vec{q},\lambda,u)
=
-\lambda u + X_0(\vec{m},\vec{q},\lambda)
-\frac{1}{2}
\sum_{1\leq l\leq k} m_l \bigl(\theta(q_{l+1}) - \theta(q_l)\bigr)
\label{Pk}
\end{equation}
where $\theta(q)=q\xi'(q)-\xi(q),$ and define
\begin{equation}
\P(\xi,u)=\inf \P_k(\vec{m},\vec{q},\lambda,u),
\label{Pu}
\end{equation}
where the infimum is taken over all $\lambda, k, \vec{m}$
and $\vec{q}.$ 
Finally, we define
\begin{equation}
\P(\xi) = \sup_{d\leq u\leq D} \P(\xi,u).
\label{Parisi1}
\end{equation}
We will first prove the following.
\begin{theorem}\label{Th1}
Given $u\in[d,D]$ and a sequence $(\eps_N)_{N\geq 1}$ that
goes to zero slowly enough, 
\begin{equation}
\lim_{N\to\infty} F_N(u,\eps_N) = \P(\xi,u).
\label{Parisi}
\end{equation}
\end{theorem}

Gaussian concentration of measure will imply that the limit
of the global free energy can be computed by maximizing the local free energy.
\begin{theorem}\label{Th2}
We have
\begin{equation}
\lim_{N\to\infty} F_N=\P(\xi).
\label{limit}
\end{equation}
\end{theorem}

{\bf Organization of the paper.} 
In Section \ref{RScase} we describe the replica symmetric region of the
model and discuss the example of the Ghatak-Sherrington model.
In Section \ref{ParisiFun} we introduce the construction 
(which we call the Parisi functional) that is used often throughout 
the paper and study some of its properties.
In Section \ref{GuerraInt} we prove the analogue of Guerra's interpolation
and explain why it seems to be necessary to impose the constraint on the 
self-overlap in order to utilize the methodology of Talagrand in \cite{T-P}.
Compared to the classical SK model where this problem does not occur, 
for the general model considered in this paper
a brand new argument is required to remove the constraint on the self-overlap 
at the end of Guerra's interpolation. This constitutes a certain nontrivial 
large deviation  problem that is solved in Section \ref{RR}.
In Section \ref{AprioriEst} we show how Theorem \ref{Th1} can be
reduced to certain apriori estimates on the error terms in Guerra's interpolation.  
For the most part, the proof of these apriori estimates goes along the
lines of the methodology developed by Talagrand in \cite{T-P} but, nonetheless, 
considerable effort is required to verify that the arguments and numerous 
computations in \cite{T-P} extend to this more general model. 
We carry out these computations in Appendix A.
In Section \ref{SecGlobal} we show how the global Parisi formula
of Theorem \ref{Th2} follows from the local Parisi formula of
Theorem \ref{Th1} and a certain concentration of measure result.
Finally, certain values of the parameter $u$ in Theorem \ref{Th1} require
small modifications of some arguments but, fortunately, these cases 
can be reduced to the classical model considered in \cite{T-P},
a work which is postponed until Appendix B.

\section{Replica symmetric region.}\label{RScase}

In this section we will describe a relatively simple necessary and 
sufficient condition in terms of the parameters of the model 
which guarantees that the infimum on the right hand side of 
(\ref{Pu}) is achieved when $k=1.$
If this happens then $\P(\xi,u)$ will be called a local 
replica symmetric solution.
In Section \ref{SecRed} we will explain that the cases when
$u=d$ or $u=D$ in Theorem \ref{Th1} can be reduced to the
classical SK model for which the domain of validity of 
the replica symmetric solution was described in \cite{SG}
and, hence, in the rest of the paper we will assume that 
\begin{equation}
d<u<D.
\label{ustri}
\end{equation}
If the infimum in (\ref{Pu}) is achieved when $k=1$ then 
\begin{equation}
\vec{m}=(0,1),\,\, \vec{q}=(0,q,u) \mbox{ for some } q\in[0,u]
\label{RSpar}
\end{equation}  
and $\lambda$ and $q$ are the only variables in
$\P_1(\vec{m},\vec{q},\lambda,u)$ which, hence, can be written as
\begin{equation}
\P_1(q,\lambda) =
-\lambda u -\frac{1}{2}(\theta(u) - \theta(q)) +\e\log
\int_{\Sigma} \exp H(\sigma) d\nu(\sigma),
\label{RSRSs}
\end{equation}
where
$$
H(\sigma)=\sigma z_0 
+\lambda \sigma^2 +\frac{1}{2}\sigma^2(\xi'(u)-\xi'(q)).
$$
Let us define the (local) replica symmetric solution by
\begin{equation}
\RS(u) = \inf_{\lambda, q} \P_1(q,\lambda)
\label{RSRSu}
\end{equation}
and describe the criterion which guarantees that $\P(\xi,u) = \RS(u).$
We will prove that if (\ref{ustri}) holds then the infimum on the
right hand side of (\ref{RSRSu}) is achieved on some $\lambda$ and $q$ 
which, therefore, must satisfy the critical point conditions
\begin{equation}
\frac{\partial \P_1}{\partial \lambda} = 
\frac{\partial \P_1}{\partial q} = 0.
\label{RSstable}
\end{equation}
From now on let $(q,\lambda)$ be such a pair, i.e. $\RS(u)=\P_1(q,\lambda).$
Suppose that $\P(\xi,u)=\RS(u).$ Then taking $k=2$ in (\ref{Pu})
should not decrease the infimum on the right hand side.
Let us take 
\begin{equation}
k=2,\,\, \vec{m}=(0,m,1) \,\,\mbox{ and }\,\, \vec{q}=(0,q,a,u) 
\,\,\mbox{ for }\,\, a\in [q,u].
\label{RSBpar}
\end{equation}
With this choice of parameters $\P_k(\vec{m},\vec{q},\lambda,u)$ becomes
\begin{equation}
\Phi(m,a) =
-\lambda u 
-\frac{1}{2}m(\theta(a)-\theta(q))
-\frac{1}{2}(\theta(u) - \theta(a)) + 
\frac{1}{m}\e\log
\e_1 X^m,
\label{RSBP}
\end{equation}
where 
\begin{equation}
X= \int_{\Sigma}\exp H'(\sigma) d\nu(\sigma)
\mbox{ and }
H'(\sigma) =\sigma (z_0 + z_1) 
+ \lambda \sigma^2 +\frac{1}{2}\sigma^2(\xi'(u)-\xi'(a))
\label{RSH}
\end{equation}
and where $\e z_0^2 = \xi'(q)$ and $\e z_1^2 = \xi'(a) - \xi'(q).$
It should be obvious that $\Phi(1,a)=\P_1(q,\lambda)$ for any 
$q\leq a\leq u.$ The derivative of $\Phi(m,a)$ with respect to $m$ at 
$m=1$ can not be positive because, otherwise, by decreasing $m$ slightly
we could decrease $\Phi(m,a)$ that would imply
$$
\P(\xi,u)\leq \Phi(m,a)<\P_1(q,\lambda) = \RS(u).
$$ 
A simple computation gives
\begin{equation}
f(a)=
\frac{\partial \Phi}{\partial m}(m,a)\bigr|_{m=1} = 
-\frac{1}{2}(\theta(a) -\theta(q)) +
\e \frac{X}{\e_1 X}\log \frac{X}{\e_1 X}
\label{RSV}
\end{equation}
and, hence, the following condition is necessary if $\P(\xi,u)=\RS(u),$
\begin{equation} 
f(a)\leq 0 \,\,\mbox{ for all }\,\, q\leq a\leq u.
\label{RSRScond}
\end{equation}
This derivative $f(a)$ represents what is usually called the
replica symmetry breaking fluctuations.
Let us note that since $\Phi(m,q)=\P_1(q,\lambda)$ does not depend
on $m,$ we have $f(q)=0.$  Also, it is easy to check that
(\ref{RSstable}) implies that $f'(q)=0.$ 
Therefore, if (\ref{RSRScond}) holds then we must have
\begin{equation}
f''(q)\leq 0,
\label{RSAT}
\end{equation}
which in the SK model is called the {\it Almeida-Thouless} condition.
It is believed (and numerical computations show) that in the 
classical SK model (\ref{RSAT}) implies (\ref{RSRScond}). 
However, we will give an example below where this is not the case
and, therefore, condition (\ref{RSRScond}) can not be weakened to
(\ref{RSAT}) in general. We will prove that (\ref{RSRScond}) 
is (necessary and) sufficient for $\P(\xi,u)=\RS(u)$.

\begin{theorem}\label{RST1}
If a pair $(\lambda, q)$ satisfies (\ref{RSstable})
and (\ref{RSRScond}) then $\P(\xi,u) = \P_1(q,\lambda)$
and such pair $(\lambda,q)$ is unique.
\end{theorem}
The proof of this theorem goes in parallel with the proof of
Theorem \ref{Th1} as will be explained in Section \ref{AprioriEst}.
However, its proof would be immediate if we knew that the
functional $\P_k$ defined in (\ref{Pk}) was convex in $\vec{m}.$
The conjecture that $\P_k$ is indeed convex in $\vec{m}$ 
was made in \cite{T-PM} and in \cite{P-PM} where
a partial result was proved. We do not give the details here but,
shortly speaking, the convexity of $\P_k$ would imply the uniqueness
of the minimum in the optimization problem (\ref{Pu}) and since
(\ref{RSRScond}) means that the replica symmetric choice of parameters 
is a local minimum in (\ref{Pu}), hence, it would be a global minimum.

\subsection{Ghatak-Sherrington model.}\label{GSsec} 

\begin{figure}[t]
\centerline{\psfig{figure=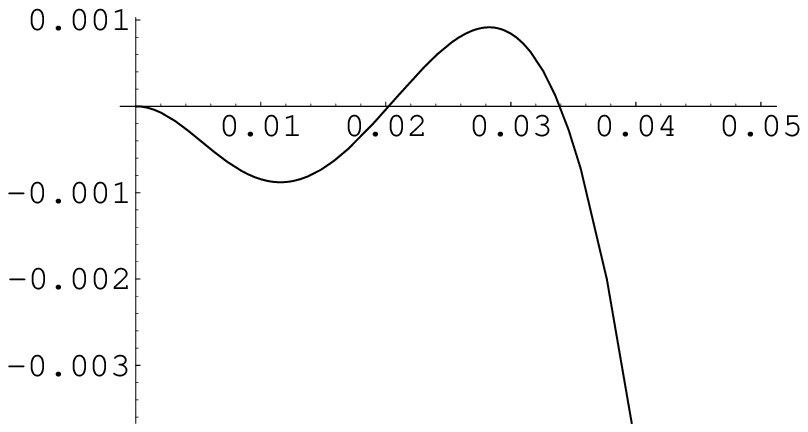,width=7cm,height=4cm}}
\caption{
\small
A function $f(a)$ for $u=0.05, \beta=17.5$ and $a\leq u.$}
\label{RSATcounter}
\end{figure}

Let us consider the Ghatak-Sherrington model introduced in \cite{GS} with
$\Sigma=\{-1,0,+1\}$, the Hamiltonian $H_N(\vsi)$ defined in (\ref{Hamorig}),
the measure $\nu$ is the counting measure on $\Sigma$ and the external field 
$h(\sigma)= h \sigma^2$ for some $h \in\Reals.$ 
This choice of parameters gives
$$
\P_1(q,\lambda) = -\lambda u -\frac{\beta^2}{4}(u^2 - q^2)
+\e \log\Bigl(
1+2\ch(z\beta\sqrt{q})\exp\bigl(\lambda + h+ \frac{1}{2}\beta^2(u-q)\bigr)
\Bigr), 
$$
where $z$ is a standard normal r.v. and we keep the dependence of $\P_1$
on $u$ implicit.
Because of the symmetry of the model, rather than the replica symmetric 
solution, one is usually interested in the case when 
$\P(\xi,u) = \P_1(q,\lambda)$ for $q=0$ and some $\lambda.$
It is easy to check that
$$
\frac{\partial \P_1}{\partial q} = \frac{\beta^2 q}{2}
- \frac{\beta^2}{2}\e\Bigl(
\frac{2\sh(z\beta\sqrt{q})\exp(\lambda +h+\beta^2(u-q)/2)}
{1+2\ch(z\beta\sqrt{q})\exp(\lambda +h+\beta^2(u-q)/2)}
\Bigr)^2
$$
and, therefore, $q=0$ always satisfies the critical point condition
(\ref{RSstable}).
If $\P(\xi,u) = \P_1(0,\lambda)$ for some $\lambda$ then we will
call $\P(\xi,u)$ a {\it paramagnetic solution} and denote it by $\PM(u).$ 
\begin{eqnarray*}
\PM(u) 
=
\inf_{\lambda} \P_1(0,\lambda)
&=& 
\inf_{\lambda}\Bigl(
-\lambda u -\frac{1}{4} \beta^2 u^2 +
\log (1 + 2e^{\lambda + h +\frac{\beta^2}{2} u})
\Bigr)
\\
&=&
hu + \frac{1}{4}\beta^2 u^2  
+\inf_{\lambda'} \Bigl(
-\lambda' u + \log(1+2e^{\lambda'}) 
\Bigr) 
\\
&=&
hu + \frac{1}{4}\beta^2 u^2 +u\log\frac{2}{u} 
+ (1-u)\log\frac{1}{1-u},
\end{eqnarray*}
where we made the change of variable 
$\lambda' = \lambda + h + \beta^2 u/2.$
The infimum is achieved on 
\begin{equation}
\lambda = -h - \frac{1}{2}\beta^2 u + \log \frac{u}{2(1-u)}.
\label{lamchange}
\end{equation}
Theorem \ref{RST1} can be applied to this model to describe
when the local free energy $\P(\xi,u)$ is given
by the paramagnetic solution $\PM(u).$
When $q=0,$ the definition (\ref{RSH}) implies that
$$
X=1 + 2e^{\lambda + h +\frac{\beta^2}{2}
(u-a)}
\ch(z_1\beta\sqrt{a})
= 1 + \frac{u}{1-u} e^{-\frac{\beta^2}{2} a} 
\ch(z_1\beta\sqrt{a}).
$$
Using the fact that
\begin{equation}
\e e^{-\frac{\beta^2}{2} a} \ch(z\beta\sqrt{a}) = 1
\label{RSik}
\end{equation}
we get $\e_1 X = 1 + u/(1-u) = 1/(1-u)$ and for $0\leq a\leq u,$ 
\begin{equation}
f(a) = -\frac{1}{4}\beta^2 a^2 + 
\e \bigl(1-u+u e^{-\frac{\beta^2}{2} a} \ch(z\beta\sqrt{a})\bigr)
\log  \bigl(1-u+u e^{-\frac{\beta^2}{2} a} \ch(z\beta\sqrt{a})\bigr).
\label{RSVGS}
\end{equation}
Since $\lambda$ in (\ref{lamchange}) and $q=0$ satisfy (\ref{RSstable}),
Theorem \ref{RST1} implies that the subset of configurations
with constrained self-overlap $R_{1,1}\approx u$ will be in 
the paramagnetic phase, $\P(\xi, u) = \PM(u)$, 
if and only if 
\begin{equation}
f(a) \leq 0 \,\,\mbox{ for }\,\, 0\leq a \leq u.
\label{RSRSGS}
\end{equation}
It is easy to check that 
$f''(0) = \beta^2(-1+\beta^2 u^2)/2$
and, therefore, (\ref{RSAT}) implies
\begin{equation}
-1 + \beta^2 u^2 \leq 0 \,\,\,\mbox{ or }\,\,\,
\beta u \leq 1.
\label{RSATGS}
\end{equation}
It is tempting to conjecture that (\ref{RSATGS}) implies
(\ref{RSRSGS}) but, unfortunately, even though it is expected to be
true in the classical SK model it is not always true here.
For example, one can check that for $u=0.05$ and
$\beta = 17.5,$ (\ref{RSATGS}) holds but (\ref{RSRSGS}) fails
(see figure \ref{RSATcounter}).

\begin{figure}[t]
\centerline{
%\psfrag{A}{\it Region $\theta_1$}
%\includegraphics{phasesmain.eps}
\psfig{figure=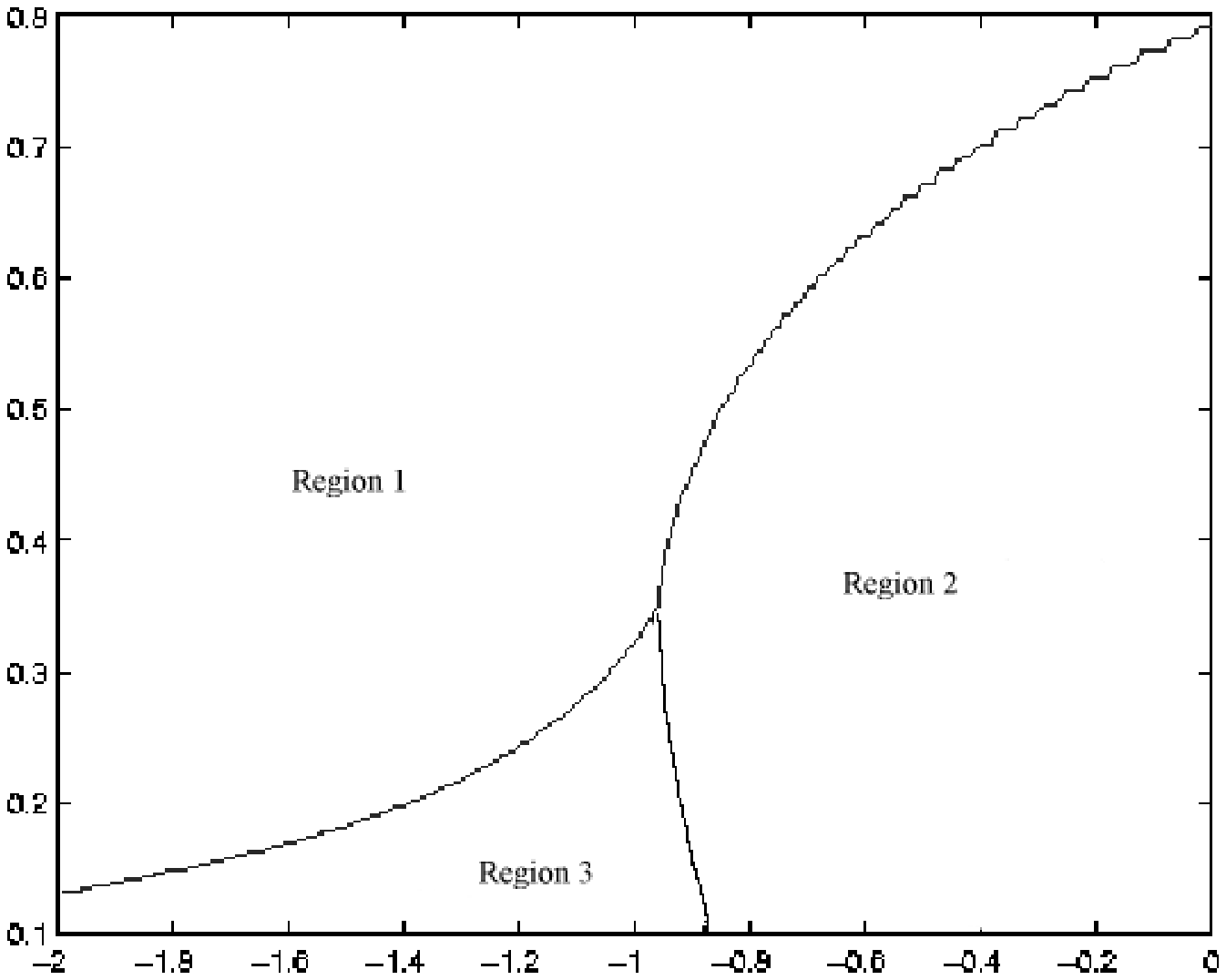,width=14cm,height=8cm}}
\caption{
\small
A phase diagram in $(h\beta^{-1}, \beta^{-1})$ coordinates.
Regions 1 and 3 are a paramagnetic phase and,
according to \cite{GS}, Region 2 is a spin glass phase.   
}
\label{RSphases}
\end{figure}

Next, we would like to describe the region of parameters 
$\beta$ and $h$ such that the system as a whole is in 
the paramagnetic phase in the sense that
\begin{equation}
\P(\xi)=\sup_{u\in[0,1]} \P(\xi,u) = \P(\xi,u') = \PM(u'),
\label{Para}
\end{equation}
i.e. the local free energy $\P(\xi,u)$ is maximized
at some point $u'$ where $\P(\xi,u')=\PM(u')$ 
and, thus, the global free energy is given by the paramagnetic 
solution, $\P(\xi)=\PM(u').$
In figure \ref{RSphases} we show a phase diagram in the coordinates
$(h\beta^{-1}, \beta^{-1})$ to compare it with \cite{GS} where 
the phase diagram was given in these coordinates. 
According to \cite{GS}, Regions 1 and 3 constitute the paramagnetic phase
where (\ref{Para}) holds, and Region 2 is the spin glass phase where
(\ref{Para}) fails, i.e. $\P(\xi,u')<\PM(u')$ for $u'$ such that
$\sup_{u} \P(\xi,u) = \P(\xi,u').$
We will explain how these regions were defined in \cite{GS} 
and argue that in Regions 1 and 3 (\ref{Para}) holds. 
We consider Regions 1 and 3 separately 
because Region 1 can be treated rigorously.

The way figure \ref{RSphases} was obtained in \cite{GS} 
is apparently as follows.
The authors used the replica symmetric approximation 
which means that instead of looking at the free energy 
$\P(\xi)$ which by Theorem \ref{Th2} is given by $\sup_{u}\P(\xi, u)$ 
they considered a replica symmetric approximation of $\P(\xi)$ given by 
\begin{equation}
\sup_{u}\RS(u) = \sup_{u}\inf_{q,\lambda}\P_1(q,\lambda).
\label{RSapprox}
\end{equation}
The value provided by this approximation, in general, is not
equal to the actual free energy $\P(\xi)$ and is only an upper
bound. However, this optimization problem is much easier than
the case of the general Parisi formula in Theorem \ref{Th2}
since (\ref{RSapprox}) depends only on three parameters $(u,q,\lambda)$.
The saddle point conditions for the solution $(u,q,\lambda)$ 
of (\ref{RSapprox}) are given by 
\begin{equation}
\frac{\partial \P_1}{\partial \lambda} = 0,\,\,
\frac{\partial \P_1}{\partial q} = 0 \,\,
\mbox{ and }\,\, \lambda=0,
\label{stable2}
\end{equation}
since it is easy to check that maximizing over $u$ and using 
$\partial \P_1/\partial \lambda = 0$ gives $\lambda=0$. 
Hence, (\ref{stable2}) reduces to solving the system of two equations.
It was predicted in \cite{GS} that the paramagnetic phase
coincides with the set of parameters $(\beta,h)$ 
for which the infimum in (\ref{RSapprox}) is attained 
at the saddle point $(u,q,\lambda)$ such that $q=0.$
This set is given by the union of Regions 1 and 3. 
On the complement, Region 2, the replica symmetric approximation
of the free energy $\P(\xi)$ is given by $\P_1(q,\lambda)$
with $q\not=0$ and the authors in \cite{GS} concluded that
it is, therefore, a spin glass phase in the sense that (\ref{Para})
fails.  However, this conclusion
in general requires further justification because (\ref{RSapprox})
is only an approximation of the general Parisi formula. 

{\bf Region 1.}
We will define Region 1 below after we explain 
a simple but important property of this model.
The fact that Region 1 is a paramagnetic phase will 
follow from this property.
First, let us observe that for any fixed $a$ 
the function $f_u(a) = f(a)$ is increasing in $u,$ where
for a moment we made the dependence on $u$ explicit. 
We have 
$$
\frac{\partial f_u(a)}{\partial u} =
\e \bigl(-1+e^{-\frac{\beta^2}{2} a} \ch(z\beta\sqrt{a})\bigr)
\log  \bigl(1-u+u e^{-\frac{\beta^2}{2} a} \ch(z\beta\sqrt{a})\bigr).
$$
It is easy to check that for any $u\in [0,1],$ the function
$x\to (-1+x)\log(1-u + ux)$ is convex for $x\geq 0$ and,
therefore, (\ref{RSik}) and Jensen's inequality imply that
$\partial f_u(a)/\partial u \geq 0.$
This means that if we take $u_1\leq u_2$ then
\begin{equation}
f_{u_2}(a) \leq 0 \mbox{ for } a\leq u_2 
\Longrightarrow
f_{u_1}(a) \leq 0 \mbox{ for } a\leq u_1. 
\label{abu}
\end{equation}
If for some $u_2\in[0,1]$ we have $\P(\xi,u_2)=\PM(u_2)$ 
then (\ref{RSRSGS}) holds for $u=u_2.$ By (\ref{abu}), 
(\ref{RSRSGS}) holds for any $u=u_1\leq u_2$ 
and, therefore, $\P(\xi,u_1)=\PM(u_1).$
This proves the following important property:
for any $\beta$ and $h$
there exists $u_0 = u_0(\beta) \in [0,1]$ such that 
\begin{equation}
\P(\xi, u) = \PM(u) \mbox{ if and only if }
u\leq u_0.
\label{RSuo}
\end{equation}
Suppose that the maximum of $\PM(u)$ is achieved on some $u_1\leq u_0.$
Then the system as a whole will be in the paramagnetic phase because
\begin{equation}
\sup_{u\in[0,1]} \P(\xi,u) \geq
\sup_{u\leq u_0} \P(\xi,u) = \sup_{u\leq u_0} \PM(u) = 
\sup_{u\in [0,1]} \PM(u) \geq \sup_{u\in[0,1]} \P(\xi,u)
\label{RSphase1cond}
\end{equation}
and, therefore, 
\begin{equation}
\P(\xi) = \sup_{u\in[0,1]} \P(\xi,u) 
= \sup_{u\leq u_0} \PM(u) = \PM(u_1).
\label{RSphase1}
\end{equation}
Region 1 is precisely where the maximum of $\PM(u)$ is achieved 
on some $u_1\leq u_0$ and, thus, it is a subset of the 
paramagnetic phase.

{\bf Region 3.}
As we mentioned above, in Region 3 the optimization problem
(\ref{RSapprox}) is solved at the saddle point $(u',q',\lambda')$ 
such that $q'=0.$ This means that for this $u',$
\begin{equation}
\sup_{u}\RS(u)=\RS(u')=\inf_{q,\lambda}\P_1(q,\lambda)=
\inf_{\lambda}\P_1(0,\lambda) = \PM(u').
\label{R3u}
\end{equation}
By itself this fact does allow us to conclude that we are in 
the paramagnetic phase, but apparently
in this particular model when this happens we also 
have $u'\leq u_0$ where $u_0$ was defined in (\ref{RSuo})
and we do not see how to prove this using calculus. 
Should this numerical observation reflect the true situation,
as seems likely, then (\ref{R3u}) would imply that
$$
\PM(u')=\P(\xi,u')\leq
\sup_{u}\P(\xi,u)\leq \sup_{u}\RS(u)=\PM(u')
$$
and we would again be in the paramagnetic phase.

{\bf Region 2.}
The fact that the infimum in  (\ref{RSapprox}) is attained
at the saddle point $(u',q',\lambda')$ with $q'\not = 0$ 
implies that $u'> u_0$ and
\begin{equation}
\sup_{u\leq u_0} \PM(u) < \RS(u').
\label{R2u}
\end{equation}
However, since $\RS(u)$ is only an upper bound
on $\P(\xi,u),$ in general, (\ref{R2u}) does not
exclude the possibility that
$$
\P(\xi, u) < \sup_{u\leq u_0} \PM(u)\,\,
\mbox{ for all }\,\, u>u_0
$$
which would imply (\ref{Para}).
At this moment we do not see how to prove that Region 2 is a
spin glass phase except by checking directly that (\ref{Para})
fails.

\section{Parisi functional.}\label{ParisiFun}

We will often consider iterative constructions 
similar to (\ref{Xlast}) and (\ref{Xl}), so it will be convenient
to define an operator that implements this recursion. In each case
we will only need to specify the parameters of the operator.
We will call this operator the {Parisi functional}.
Given $k\geq 1,$ consider a vector
\begin{equation}
\vec{m} = (m_0, \ldots, m_{k})
\label{Pmono}
\end{equation}
such that all coordinates $m_i\geq 0.$ 
Let $\vec{z} = (z_{l})_{0\leq l\leq k}$ be a collection
of independent random vectors.

Suppose that we are given a random variable $F$
that is a function of $\vec{z},$ $F = F(\vec{z}).$
Then we let $F_{k+1}=F$ and for $0\leq l\leq k$ define
iteratively
\begin{equation}
F_l = \frac{1}{m_{l}}\log\e_l \exp m_{l} F_{l+1},
\label{Piter}
\end{equation}
where $\e_l$ denotes the expectation in $(z_{i})_{i\geq l}.$
When $m_{l} = 0,$ this means $F_l = \e_l F_{l+1}.$ 

{\bf Definition.}
We define the Parisi functional by
\begin{equation}
\P(\vec{m})F = F_0.
\label{PO}
\end{equation}

With these notations the definition of $X_0$ given by
(\ref{m}) - (\ref{Xl}) can be written as
\begin{equation}
X_0 = \P(\vec{m})X_{k+1}.
\label{PXlast}
\end{equation}
Let us describe several immediate properties of the Parisi
functional that will be often used throughout the paper.
It is obvious by induction in (\ref{Piter}) 
that for any constant $c\in\Reals,$
\begin{equation}
\P(\vec{m})(c+F)=c+\P(\vec{m})F.
\label{Pprop1}
\end{equation}
Similarly, if $\vec{z}^j = (z_{p}^j)_{0\leq p\leq k}$ 
are independent for $j=1,2$ and
$F^j=F^j(\vec{z}^j)$ then
\begin{equation}
\P(\vec{m})(F^1+F^2)=\P(\vec{m})F^1 + \P(\vec{m})F^2.
\label{Pprop2}
\end{equation}
If $F\leq F'$ then
\begin{equation}
\P(\vec{m})F\leq \P(\vec{m})F'.
\label{Pprop3}
\end{equation}
The next property plays an important role in Talagrand's interpolation
for two copies of the system. Suppose that we have two random variables
\begin{equation}
F^j = F\Bigl(\sum_{p\leq k} z_{p}^j\Bigr) \mbox{ for } j=1,2
\label{Pprop40}
\end{equation}
such that 
\begin{equation}
z_{l}^1 = z_{l}^2 \mbox{ for } l< r
\mbox{ and } 
z_{l}^1 \mbox{ and } z_{l}^2 \mbox{ are independent copies for }
l\geq r.
\label{Pprop41}
\end{equation}
Let us define $\vec{n}=(n_{p})_{p\leq k}$
such that
\begin{equation}
n_{l} = \frac{m_{l}}{2} \mbox{ for } l< r
\mbox{ and }
n_{l} = m_{l} \mbox{ for } l\geq r.
\label{Pprop42}
\end{equation}

\begin{lemma}\label{Plem1}
Given (\ref{Pprop40}), (\ref{Pprop41}) and (\ref{Pprop42})
let $F= F^1 + F^2$ and define $F_l^j$ by (\ref{Piter}) 
and $F_l$ by (\ref{Piter}) with $\vec{m}$ replaced by $\vec{n}.$
Then,
\begin{equation}
F_l = F_l^1+F_l^2 \mbox{ for } l> r \mbox{ and }
F_l = 2 F_l^1 = 2 F_l^2 \mbox{ for } l \leq r.
\label{Pprop43}
\end{equation}
In particular,
\begin{equation}
\P(\vec{n})(F^1+F^2) = 2\P(\vec{m})F^1.
\label{Pprop4}
\end{equation}
\end{lemma}
{\bf Proof.}
The proof follows by induction in (\ref{Piter}).
For $l> r,$ using independence of
$z_{l}^1$ and $z_{l}^2$ we get
$$
F_l=
\frac{1}{n_{l}}\log\e_l \exp n_{l} F_{l+1}
=
\frac{1}{m_{l}}\log\e_l \exp m_{l} 
(F_{l+1}^{1}+F_{l+1}^2)
= F_l^1 + F_l^2.
$$
For $l=r$ all independent copies already have been averaged and since
$z_{l}^1 = z_{l}^2$ for $l< r,$ 
$
F_{r}^1 = F_{r}^2 \mbox{ and } 
F_{r} = 2F_{r}^1.
$
Finally, by induction for $l < r$ we get
$$
F_l =
\frac{1}{n_{l}}\log\e_l \exp n_{l} F_{l+1}
=
\frac{2}{m_{l}}\log\e_l \exp \frac{m_{l}}{2} 
2F_{l+1}^{1}
= 2F_l^1
$$
which for $l=0$ proves (\ref{Pprop4}).
\qed

The next property concerns the computation of the derivatives
of $\P(\vec{m})F.$ With the notations of (\ref{Pmono}) and 
(\ref{Piter}), let us define 
\begin{equation}
W_l = \exp m_{l}(F_{l+1} - F_{l}).
\label{PW}
\end{equation}
Note that by definition of $F_l$ in (\ref{Piter}) we have
$\e_l W_l = 1.$ Also, since $F_l,$ $F_{l+1}$ and $W_l$ do not depend
on $z_{i}$ for $i\geq l+1,$ we can write
$$
\e_l W_l W_{l+1} = \e_l \e_{l+1} W_l W_{l+1} =\e_l W_l \e_{l+1} W_{l+1}
= 1.
$$
Repeating the same argument,
\begin{equation}
\e_l W_l\ldots W_{k} = 1.
\label{Pden1}
\end{equation}
This fact will be used often below.
\begin{lemma}
For a generic variable $x,$
\begin{equation}
\frac{\partial \P(\vec{m})}{\partial x}
=
\e_0 W_0\ldots W_k \frac{\partial F}{\partial x}.
\label{Pder}
\end{equation}
\end{lemma}
{\bf Proof.}
By (\ref{Piter}), 
$\exp m_{l} F_{l} = \e_l \exp m_{l} F_{l+1}$
and, therefore,
$$
m_{l} \exp m_{l} F_{l} 
\frac{\partial F_{l}}{\partial x}
=
m_{l} \e_l \exp m_{l} F_{l+1} 
\frac{\partial F_{l+1}}{\partial x}.
$$ 
Since $F_{l}$ does not depend on $z_{i}$ for $i\geq l,$
$$
\frac{\partial F_{l}}{\partial x}
=
\e_l W_l \frac{\partial F_{l+1}}{\partial x},
$$
where $W_l$ was defined in (\ref{PW}). Applying the
same equation to $F_{l+1}$ gives
$$
\frac{\partial F_{l}}{\partial x}
=
\e_l W_l \e_{l+1} W_{l+1}\frac{\partial F_{l+2}}{\partial x}
=
\e_l W_l W_{l+1}\frac{\partial F_{l+2}}{\partial x}
$$
since $F_{l},F_{l+1}$ and, therefore, $W_l$ do not
depend on $(z_{i})$ for $i\geq l+1.$
Repeating the same argument inductively we get
\begin{equation}
\frac{\partial F_{l}}{\partial x}
=
\e_l W_l\ldots W_k \frac{\partial F}{\partial x},
\label{Pder2}
\end{equation}
which for $l=0$ implies (\ref{Pder}).
\qed

We will often assume that the coordinates
of $\vec{m}$ are arranged in a nondecreasing order,
\begin{equation}
m_0\leq\ldots \leq m_k.
\label{arrange}
\end{equation}
The following lemma provides a useful control of the expressions of 
the type (\ref{Pder}) that appear as the derivatives of the Parisi functional.

\begin{lemma}\label{Plemlast}
Suppose (\ref{arrange}) holds.
Let $f=f(\vec{z})$ and $F=F(\vec{z})$ and let $W_l$ be
defined by (\ref{PW}). Suppose that $f\leq F$ and let 
$m=m_{r}$ be the first nonzero element in (\ref{arrange}).
Then
\begin{equation}
\e_0 \log \e_r W_r\ldots W_k \exp m_{k}(f-F) \leq 
 m (\P(\vec{m})f - \P(\vec{m})F).
\label{Pprop5}
\end{equation}
\end{lemma}
{\bf Proof.}
Let us define $U=\exp m_{k}(f-F).$ 
Using (\ref{PW}), we can write
$$
W_k U = \exp m_{k} (f - F_k) 
$$
and since $F_k$ does not depend on $z_{k},$ this implies that
$$
\e_k W_k U = \exp(-m_{{k}} F_k)\e_k   \exp m_{k} f
=  \exp m_{k} (f_k - F_k).
$$
We will proceed by induction to show that for $r\leq l\leq k,$
\begin{equation}
\e_l W_l\ldots W_k U  \leq  \exp m_{{l}} (f_l - F_l).
\label{Plhere}
\end{equation}
As in (\ref{Pprop3}), $f\leq F$ implies that $f_l \leq F_l$ 
and since $m_{{l-1}} \leq m_{l},$ (\ref{Plhere}) implies
that
$$
\e_l W_l\ldots W_k U  \leq  \exp m_{{l-1}} (f_l - F_l).
$$ 
Multiplying both sides by $W_{l-1}$ gives
$$
\e_l W_{l-1} W_l\ldots W_k U  \leq  \exp m_{{l-1}} (f_l - F_{l-1})
$$
since $W_l$ does not depend on $z_{i}$ for $i\geq l.$
Taking the expectation $\e_{l-1}$ and using that 
$F_{l-1}$ does not depend on $z_{i}$ for $i\geq l-1$
we can write
$$
\e_{l-1} W_{l-1} W_l\ldots W_k U  \leq  
\exp(- m_{{l-1}} F_{l-1}) \e_{l-1} \exp m_{{l-1}} f_{l} 
=  \exp m_{{l-1}} (f_{l-1} - F_{l-1}).
$$
This finishes the proof of the induction step. 
For $l=r,$ (\ref{Plhere}) implies that
\begin{equation}
\log \e_r W_r\ldots W_k U  \leq  m (f_r - F_r).
\label{ppp}
\end{equation}
Since for $l<r,$ $m_{l} = 0,$ (\ref{Piter}) implies that
$$
\P(\vec{m}) f = f_0 = \e_0 f_r \mbox{ and } \P(\vec{m})F = F_0 = \e_0 F_r
$$  
and, therefore, taking the expectation of both sides of (\ref{ppp})
proves (\ref{Pprop5}).
\qed

\section{Guerra's interpolation.}\label{GuerraInt}

The first step of the proof of Theorem \ref{Th1} 
is the analogue of Guerra's interpolation method in \cite{Guerra}.
For $1\leq i\leq N,$ we consider independent copies $(z_{i,p})_{0\leq p\leq k}$ 
of the sequence $(z_p)_{0\leq p\leq k}$ defined in (\ref{z})
that are also independent of the randomness of the Hamiltonian $H_N(\vsi).$ 
Consider $\vec{m}=(m_p)_{0\leq p\leq k}$ as in (\ref{m}). 
We denote by $\e_l$ the expectation 
in the r.v. $(z_{i,p})_{i\leq N, p\geq l}.$
Consider the Hamiltonian
\begin{equation}
H_t(\vsi)=\sqrt{t}H_N(\vsi) + \sqrt{1-t}\sum_{i\leq N} \sigma_i
\Bigl(\sum_{0\leq p\leq k} z_{i,p}\Bigr).
\label{GHt}
\end{equation}
Given $U_N$ defined in (\ref{const}), let
\begin{equation}
F = \log \int_{U_N} \exp H_t(\vsi) d\nu(\vsi)
\label{GFtlast}
\end{equation}
and let
\begin{equation} 
\varphi_N(t)=\frac{1}{N}\e \P(\vec{m}) F,
\label{Gphi}
\end{equation}
where $\e$ denotes the expectation in all random variables 
including the randomness of the Hamiltonian $H_N(\vsi)$.

For a function $h:\Sigma^N\to\Reals,$
let $\la h \ra_t$ denote its average 
with respect to the Gibbs measure with Hamiltonian $H_t(\vsi)$ in (\ref{GHt})
on the set $U_N,$  i.e.
\begin{equation}
\la h \ra_t \exp F = \int_{U_N} 
h(\vsi) \exp H_t(\vsi) d\nu(\vsi).
\label{Gaverage}
\end{equation}
(\ref{Pden1}) implies that the functional
\begin{equation}
h\to \e_l(W_l\ldots W_k\la h\ra_t) 
\label{Ggammal}
\end{equation}
is a probability $\gamma_l$ on $U_N.$
We denote by $\gamma_l^{\otimes 2}$ its product on 
$U_N\times U_N,$ and for a function $h:U_N\times U_N\to\Reals$
we set
\begin{equation}
\mu_l(h)=\e (W_1\ldots W_{l-1}\gamma_l^{\otimes 2}(h)).
\label{Gmul}
\end{equation}
The following Gaussian integration by parts will be commonly used below. 
If $g$ is a Gaussian random variable then for a function 
$F:\Reals\to\Reals$ of moderate growth we have (A.40 in \cite{SG}),
\begin{equation}
\e g F(g) = \e g^2 \e F'(g).
\label{GGI0}
\end{equation}
This can be generalized as follows.
If $\vec{g}=(g_1,\ldots,g_n)$ is a jointly Gaussian family 
of random variables then for a function 
$F:\Reals^n\to\Reals$ of moderate growth
we have (see for example, A.41 in \cite{SG}),
\begin{equation}
\e g_i F(\vec{g}) = \sum_{j\leq n} \e (g_i g_j)
\e \frac{\partial F}{\partial g_j}(\vec{g}).
\label{GGIold}
\end{equation}
We will need a similar statement for functionals of
not necessarily finite Gaussian families, for example,
for a random process $H_N(\vsi)$ indexed by $\vsi$ in a possibly infinite
set $\Sigma^N.$ The following is a simple consequence of (\ref{GGI0}).
\begin{lemma}\label{GLem4}
Let $\vec{g}=(g(\vrho))_{\vrho\in U}$ be a 
Gaussian process indexed by $U\subseteq \Reals^n$ 
and let $F(\vec{g})$ be a differentiable functional on $\Reals^U.$ 
Given $\vsi\in U,$ we have
\begin{equation}
\e g(\vsi) F(\vec{g}) = 
\e \frac{\delta F}{\delta \vec{g}}[\e g(\vsi) g(\vrho)]
\label{GGI}
\end{equation}
- the expectation of the variational derivative of $F$ in the direction 
$h(\vrho) = \e g(\vsi) g(\vrho).$
\end{lemma}
{\bf Proof.}
Consider a process $\vec{g}'=(g'(\vrho))_{\vrho\in U}$ defined by
$$
g'(\vrho)=g(\vrho)-g(\vsi)\frac{\e g(\vrho)g(\vsi)}{\e g(\vsi)^2},
$$
which is, obviously, independent of the r.v. $g(\vsi).$
If we fix $\vec{g}'$ and denote by $\e'$ the expectation
with respect to $g(\vsi)$ then, using (\ref{GGI0}) with $g=g(\vsi)$
gives
$$
\e' g(\vsi) F(\vec{g}) =\e' g(\vsi) F\Bigl(
g'(\vrho) + g(\vsi)\frac{\e g(\vrho)g(\vsi)}{\e g(\vsi)^2}
\Bigr)
=
\e' \frac{\delta F}{\delta \vec{g}}[\e g(\vsi) g(\vrho)].
$$
Taking the expectation in $\vec{g}'$ proves (\ref{GGI}).
\qed

We are ready to prove the main result of this section.
The proof will clarify why we first compute the free energy of 
the set of configurations with constrained self-overlap 
$R_{1,1}\in[u-\eps_N,u+\eps_N].$ 

\begin{theorem}\label{GTh3} (Guerra's interpolation).
For $t\in [0,1]$ we have
\begin{eqnarray}
\varphi_N'(t)
&=&
-\frac{1}{2}\sum_{1\leq l\leq k}m_l(\theta(q_{l+1})-\theta(q_l))
\nonumber
\\
&&
-\frac{1}{2}\sum_{1\leq l\leq k} (m_{l} - m_{l-1})
\mu_l\bigl(
\xi(R_{1,2})-R_{1,2}\xi'(q_l) +\theta(q_l) 
\bigr) +{\cal R},
\label{Gphider}
\end{eqnarray}
where $|{\cal R}|\leq c(N)+ L\eps_N.$
\end{theorem}
{\bf Proof.}
The proof of this theorem repeats the proof of the main result in \cite{Guerra}
(see also Theorem 2.1 in \cite{T-P}) with some necessary modifications. 
We will give the detailed proof in order to demonstrate how 
Lemma \ref{GLem4} replaces (\ref{GGIold}) and to show that the
constraint $q_{k+1}=u$ in (\ref{q}), in some sense, matches
the constraint on the self overlap $R_{1,1}\in [u-\eps_N,u+\eps_N]$ 
in the definition of local free energy (\ref{const}), (\ref{FE2}).
In particular, we will see in (\ref{GThfinish}) below that
such choice of parameters allows us to get rid of a certain
term in the derivative $\varphi_N'(t)$ that, otherwise, 
would be problematic to control.
First of all, (\ref{Pder}) implies that
$$
\varphi_N'(t)=\frac{1}{N}\e W_1\ldots W_k 
\frac{\partial F}{\partial t}.
$$
Using the fact that $m_{k}=1,$ one can write
\begin{eqnarray}
W_1\ldots W_k 
&=& 
\exp \sum_{1\leq l\leq k} m_l(F_{l+1} - F_{l})
\nonumber
\\
&=&
\exp \Bigl(F +
\sum_{1\leq l\leq k} (m_{l-1} - m_{l}) F_{l} \Bigr)
= T \exp F,
\label{GWT}
\end{eqnarray}
where 
$
T=T_1\ldots T_k \mbox{ and }
T_l = \exp (m_{l-1} - m_{l}) F_{l}.
$
Using (\ref{GFtlast}) we can write
$$
\varphi_N'(t)=
\frac{1}{N}\e T \exp F
\frac{\partial F}{\partial t}
=\mbox{I}-\mbox{II},
$$
where
\begin{equation}
\mbox{I}=
\frac{1}{2N\sqrt{t}} 
\int_{U_N} \e T H_N(\vsi) \exp H_t(\vsi) d\nu(\vsi).
\label{GI}
\end{equation}
and
\begin{equation}
\mbox{II}=
\frac{1}{2N\sqrt{1-t}} 
\int_{U_N}
\e T \sum_{i\leq N}\sigma_i \sum_{p\leq k} z_{i,p} \exp H_t(\vsi)
d\nu(\vsi).
\label{GII}
\end{equation}
To compute I, we will use (\ref{GGI}) for the family 
$\vec{g} = (H_N(\vrho))_{\vrho\in U_N}$ and
we will think of each factor $T_l$ in $T$ as the functional of
$\vec{g}.$ Let us denote
\begin{equation}
\zeta(\vsi,\vrho)=\frac{1}{N}\e H_N(\vsi)H_N(\vrho).
\label{Gzeta}
\end{equation}
Then (\ref{GGI}) and (\ref{Gzeta}) imply
\begin{eqnarray} 
\mbox{I}
&=&
\frac{1}{2\sqrt{t}}
\int_{U_N} \e T
\frac{\partial \exp H_t(\vsi)}{\partial H_N(\vsi)} 
\zeta(\vsi,\vsi) d\nu(\vsi)
\label{GI1}
\\
&+&
\frac{1}{2\sqrt{t}}\sum_{l\leq k}
\int_{U_N}
\e T\exp H_t(\vsi)\frac{1}{T_l}
\frac{\delta T_l}{\delta \vec{g}}[\zeta(\vsi,\vrho)] d\nu(\vsi) .
\nonumber
\end{eqnarray}
First of all,
$$
\frac{\partial \exp H_t(\vsi)}{\partial H_N(\vsi)}
=\sqrt{t} \exp H_t(\vsi),
$$
and, therefore, the first line in (\ref{GI1})
can be written as
\begin{equation}
\frac{1}{2}
\int_{U_N} \e T \exp H_t(\vsi) \zeta(\vsi,\vsi)d\nu(\vsi)
=\frac{1}{2}\e W_1\ldots W_k \bigl\la
\zeta(\vsi,\vsi)
\bigr\ra_t.
\label{GI2}
\end{equation}
Using the definition of $T_l$ we can write
$$
\frac{1}{T_l}\frac{\delta T_l}{\delta \vec{g}}[\zeta(\vsi,\vrho)]
=
(m_{l-1}-m_l) \frac{\delta F_{l}}{\delta \vec{g}}[\zeta(\vsi,\vrho)]
=
(m_{l-1}-m_l)\e_l W_l\ldots W_k 
\frac{\delta F}{\delta \vec{g}}[\zeta(\vsi,\vrho)],
$$
where we used (\ref{Pder2}).
We have
$$
\frac{\delta F}{\delta \vec{g}}[\zeta(\vsi,\vrho)]
=
\sqrt{t}\exp(-F)\int_{U_N}\zeta(\vsi,\vrho)
\exp H_t(\vrho) d\nu(\vrho)
=
\sqrt{t}\bigl\la \zeta(\vsi,\vrho)\bigr\ra_t',
$$
where $\la\cdot\ra_t'$ denotes the Gibbs average with respect to
$\vrho$ for a fixed $\vsi.$ 
Therefore, for a fixed $\vsi$ we get
$$
\frac{1}{T_l}\frac{\delta T_l}{\delta \vec{g}}[\zeta(\vsi,\vrho)]
=
\sqrt{t}(m_{l-1}-m_l)\e_l W_l\ldots W_k 
\bigl\la \zeta(\vsi,\vrho)\bigr\ra_t'
=
\sqrt{t}(m_{l-1}-m_l)\gamma_l(\zeta(\vsi,\vrho)),
$$
where $\gamma_l$ was defined in (\ref{Ggammal}). 
Hence, the second line in (\ref{GI1}) is equal to
\begin{eqnarray}
&&
\frac{1}{2}\sum_{l\leq k}(m_{l-1}-m_l)
\e W_1\ldots W_k \exp(-F)
\int_{U_N} \gamma_l(\zeta(\vsi,\vrho))
\exp H_t(\vsi) d\nu(\vsi)
\label{GI3}
\\
&&
=
\frac{1}{2}\sum_{l\leq k}(m_{l-1}-m_l)
\e W_1\ldots W_{l-1}\gamma_l^{\otimes 2}
(\zeta(\vsi,\vrho))
=\frac{1}{2}\sum_{l\leq k}(m_{l-1}-m_l)
\mu_l(\zeta(\vsi,\vrho)),
\nonumber
\end{eqnarray}
where $\mu_l$ was defined in (\ref{Gmul}).
Combining (\ref{GI2}) and (\ref{GI3}) we get
\begin{equation}
\mbox{I} = 
\frac{1}{2}\e W_1\ldots W_k \bigl\la
\zeta(\vsi,\vsi)
\bigr\ra_t 
+\frac{1}{2}\sum_{l\leq k}(m_{l-1}-m_l)
\mu_l(\zeta(\vsi,\vrho)).
\label{GI4}
\end{equation}
By (\ref{correlation}), for any $\vsi^1,\vsi^2\in\Sigma^N$ we have
$$
|\zeta(\vsi^1,\vsi^2) - \xi(R_{1,2})|\leq c(N)
$$
and, therefore, (\ref{GI4}) implies that
\begin{equation}
\mbox{I} = 
\frac{1}{2}\e W_1\ldots W_k \bigl\la
\xi(R_{1,1})
\bigr\ra_t 
+\frac{1}{2}\sum_{l\leq k}(m_{l-1}-m_l)
\mu_l(\xi(R_{1,2})) + {\cal R},
\label{GI5}
\end{equation}
where $|{\cal R}|\leq c(N).$
The computation of $\mbox{II}$ is very similar, one only
needs to note that $F_{l}$ does not depend on
$z_{i,p}$ for $l\leq p.$ We have
\begin{equation}
\mbox{II} = 
\frac{1}{2} \xi'(q_{k+1})\e W_1\ldots W_k \bigl\la
R_{1,1}
\bigr\ra_t 
+\frac{1}{2}\sum_{l\leq k}(m_{l-1}-m_l)\xi'(q_l)
\mu_l(R_{1,2}).
\label{GI6}
\end{equation}
Combining (\ref{GI5}) and (\ref{GI6}) and rearranging terms,
it is easy to see that
\begin{eqnarray}
\varphi_N'(t) 
&=& 
\frac{1}{2}\e W_1\ldots W_k \bigl\la
\xi(R_{1,1})-R_{1,1}\xi'(q_{k+1}) +\theta(q_{k+1})
\bigr\ra_t 
-\frac{1}{2}\sum_{l\leq k}m_l(\theta(q_{l+1})-\theta(q_l))
\nonumber
\\
&-&
\frac{1}{2}\sum_{1\leq l\leq k} (m_{l} - m_{l-1})
\mu_l\bigl(
\xi(R_{1,2})-R_{1,2}\xi'(q_l) +\theta(q_l) 
\bigr) +{\cal R}.
\label{GThfinish}
\end{eqnarray}
It now suffices to notice that
$\la\cdot\ra_t$ is restricted to the set $U_N$ and
$q_{k+1}=u$ by (\ref{q}), which implies that
$$
0\leq \xi(R_{1,1})-R_{1,1}\xi'(u) +\theta(u) \leq
L|R_{1,1} - u| \leq L\eps_N.
$$
Since $\e W_1\ldots W_k =1,$ this finishes the proof 
of Theorem \ref{GTh3}.
It is to control the first term on the right hand side of (\ref{GThfinish}) 
that we impose the constraint on self-overlap.
In the classical Sherrington-Kirkpatrick model this problem did not occur
because $R_{1,1}$ was always $1.$ 
\qed

\section{Removing the constraint on the self-overlap.}\label{RR}

The main goal of this section is to compute $\lim_{N\to\infty}\varphi_N(0).$
The convexity of $\xi$ implies that $\xi(a) - a\xi'(b) + \theta(b) \geq 0$
for any $a,b\in\Reals$ and, therefore, (\ref{Gphider}) implies
$$
\varphi_N'(t)
\leq
-\frac{1}{2}\sum_{1\leq l\leq k}m_l(\theta(q_{l+1})-\theta(q_l))
+\R
$$
and, hence,
\begin{equation}
F_N(u,\eps_N)=\varphi_N(1)\leq \varphi_N(0) 
-\frac{1}{2}\sum_{1\leq l\leq k}m_l(\theta(q_{l+1})-\theta(q_l))
+\R.
\label{GuerraRSB}
\end{equation}
In the SK model $\varphi_N(0)$ was very easy to compute and, in fact,
it was independent of $N$ because at the end of Guerra's interpolation the
spins became decoupled. The situation is different here
because of the constraint (\ref{const}). Let us recall that for $t=0,$ 
$$
\varphi_N(0)
%=\frac{1}{N} \e\P(\vec{m}) F
= \frac{1}{N} \P(\vec{m}) F
\,\,\mbox{ where }\,\,
F = \log \int_{U_N} \exp \sum_{i\leq N} \sigma_i
\Bigl(\sum_{0\leq p\leq k} z_{i,p}\Bigr)d\nu(\vsi).
$$
Removing the constraint $U_N$ constitutes a large deviation problem 
that will be addressed in Theorem \ref{Th4} below.
First of all, let us give an easy upper bound on $\varphi_N(0).$
Since by (\ref{diam1}) the self-overlap $R_{1,1}\in[d, D],$ 
for $\vsi\in U_N$ we have $R_{1,1}\in [d,D]\cap [u-\eps_N,u+\eps_N]$
and, therefore, 
\begin{equation}
-\lambda R_{1,1}\leq -\lambda u_N(\lambda),
\end{equation}
where 
\begin{equation}
u_N(\lambda)=\max(d,u-\eps_N) \mbox{ for } \lambda\geq 0
\mbox{ and }
u_N(\lambda)=\min(D,u+\eps_N) \mbox{ for } \lambda<0.
\label{uenn}
\end{equation}
Using this, we can bound $F$ as follows,
\begin{eqnarray*}
F
&=&
\log \int_{U_N} \exp 
\Bigl(
\sum_{i\leq N}\sigma_i
\sum_{0\leq p\leq k} z_{p,i}
\Bigr)
d\nu(\vsi)
\\
&\leq&
-N\lambda u_N(\lambda) + \log\int_{U_N} \exp\Bigl(
\sum_{i\leq N}\sigma_i\sum_{0\leq p\leq k} z_{p,i} 
+\lambda\sum_{i\leq N}\sigma_i^2 
\Bigr)d\nu(\vsi)
\\
&\leq&
-N\lambda u_N(\lambda) + \log \int_{\Sigma^N} \exp\Bigl(
\sum_{i\leq N}\sigma_i\sum_{0\leq p\leq k} z_{p,i} 
+\lambda\sum_{i\leq N}\sigma_i^2
\Bigr)d\nu(\vsi)
\\
&=&
-N\lambda u_N(\lambda)  +\sum_{i\leq N} X_{k+1,i}
\end{eqnarray*}
where 
\begin{equation}
X_{k+1,i}=\log \int_{\Sigma}
\Bigl(\sigma\sum_{0\leq p\leq k} z_{p,i} 
+\lambda \sigma^2\Bigr)d\nu(\sigma)
\label{Xcopies}
\end{equation}
are independent copies of $X_{k+1}$ defined in (\ref{Xlast}).
Using (\ref{Pprop1}) and (\ref{Pprop3}),
$$
\varphi_N(0) = \frac{1}{N} \P(\vec{m})F \leq 
-\lambda u_N(\lambda)  + \P(\vec{m})X_{k+1}
=-\lambda u_N(\lambda) +  X_0.
$$
The arbitrary choice of $\lambda$ here implies that
\begin{equation}
\varphi_N(0)\leq \inf_{\lambda} (-\lambda u_N(\lambda) + X_0).
\label{var0}
\end{equation}
Combining (\ref{GuerraRSB}) and (\ref{var0}) we get
\begin{equation}
F_N(u,\eps_N) \leq
\inf\Bigl(
-\lambda u_N(\lambda) + X_0(\vec{m},\vec{q},\lambda)
-\frac{1}{2}\sum_{1\leq l\leq k} m_l(\theta(q_{l+1})-\theta(q_l))
\Bigr) +\R,
\label{oneside}
\end{equation}
where the infimum is over all choices of parameters
$k,\vec{m},\vec{q}$ and $\lambda.$ 
The bound (\ref{oneside}) is the analogue of Guerra's 
replica symmetry breaking bound in \cite{Guerra}. 
If instead of $u_N(\lambda)$ we had $u$ in (\ref{oneside})
then the infimum would be equal to $\P(\xi,u)$ which would
prove the upper bound in Theorem \ref{Th1}. We will now show 
that for 
\begin{equation}
d<u<D
\label{strict}
\end{equation}
this infimum is not changed much by replacing
$u_N(\lambda)$ with $u$. We will need the following.

\begin{lemma}\label{Lextra1}
There exists a function $a(\lambda)$ such that
for any $k,\vec{m},\vec{q},$
\begin{equation}
a(\lambda)\lambda\leq X_0(\vec{m},\vec{q},\lambda)
\label{ab}
\end{equation}
and such that 
\begin{equation}
\lim_{\lambda\to +\infty}a(\lambda)=D 
\mbox{ and }
\lim_{\lambda\to -\infty}a(\lambda)=d.
\label{alim}
\end{equation}
\end{lemma}
{\bf Proof.} Indeed, if in the recursive construction (\ref{Xl}) 
one takes all $m_l=0$ then H\"older's inequality 
yields that for any sequence $\vec{m}$
\begin{equation}
\e \log\int_{\Sigma}
\exp \Bigl(\sigma\sum_{0\leq p\leq k} z_p 
+ \lambda\sigma^2 \Bigr)d\nu(\sigma)
\leq X_0(\vec{m},\vec{q},\lambda) 
\label{twochoices}
\end{equation}
Since the function
$$
x\to \log\int_{\Sigma}
\exp \Bigl(\sigma x 
+ \lambda\sigma^2 \Bigr)d\nu(\sigma)
$$
is convex by H\"older's inequality, (\ref{twochoices})
and Jensen's inequality imply that
\begin{equation}
\log\int_{\Sigma}
\exp \lambda\sigma^2 d\nu(\sigma)
\leq X_0(\vec{m},\vec{q},\lambda). 
\label{above}
\end{equation}
It is clear that (\ref{diam1}) implies that the left hand side of (\ref{above})
is asymptotically equivalent to $D\lambda$ for $\lambda\to + \infty$
and to $d\lambda$ for $\lambda\to -\infty$ and this proves Lemma \ref{Lextra1}.
\qed

Next we will show that
the estimate (\ref{ab}) and equation (\ref{oneside}) imply the upper
bound in Theorem \ref{Th1}.
\begin{lemma}\label{Lextra2}
For $d<u<D$ we have
\begin{equation}
\limsup_{N\to\infty }F_N(u,\eps_N)\leq \P(\xi,u).
\label{halft1}
\end{equation}
\end{lemma}
{\bf Proof.} Lemma \ref{Lextra1} implies that
\begin{equation}
-\lambda u_N(\lambda) + X_0(\vec{m},\vec{q},\lambda) \geq
\lambda (a(\lambda) - u_N(\lambda)).
\label{aminu}
\end{equation}
For $d<u<D,$ the definition (\ref{uenn}) implies that
for any $\lambda,$ $\lim_{N\to\infty} u_N(\lambda) = u.$
Combining this with (\ref{alim}) yields that 
for $N$ large enough the right hand side of (\ref{aminu})
goes to infinity as $\lambda\to\pm\infty$ and,
thus, the infimum in (\ref{oneside}) is achieved for
$|\lambda|\leq \Lambda$ where $\Lambda$ is a large enough constant
independent of $k,\vec{m},\vec{q}.$ This means that in (\ref{oneside}) 
restricting minimization over $\lambda$ to the set 
$\{|\lambda|\leq \Lambda\}$ does not change the infimum and, therefore,
\begin{eqnarray*}
F_N(u,\eps_N) 
&\leq&
\inf\Bigl(
-\lambda u + X_0(\vec{m},\vec{q},\lambda)
-\frac{1}{2}\sum_{1\leq l\leq k} m_l(\theta(q_{l+1})-\theta(q_l))
\Bigr) + \Lambda |u_N(\lambda) - u| + \R,
\\
&=&
\P(\xi,u) + \Lambda |u_N(\lambda) - u| + \R.
\end{eqnarray*}
and this finishes the proof.
\qed

In the rest of the section we will show that the bound 
in (\ref{var0}) is exact in the limit. 

\begin{theorem}\label{Th4}
For any $d<u<D,$ if the sequence $\eps_N$ goes to zero
slowly enough then 
\begin{equation}
\lim_{N\to \infty}\varphi_N(0) = 
\varphi_0(\vec{m},\vec{q}):= \inf_{\lambda}(-\lambda u 
+X_0(\vec{m},\vec{q},\lambda)).
\label{philimit}
\end{equation}
\end{theorem}
{\bf Proof.}
Given a measurable set $A\subseteq [d,D]$
and $\lambda\in\Reals,$ we define
\begin{equation}
F(A,\lambda)=\log\int_{\{R_{1,1}\in A\}}\exp\Bigl(
\sum_{i=1}^N \sigma_i \sum_{0\leq p\leq k} z_{i,p} 
+ \lambda \sum_{i=1}^N \sigma_i^2 
\Bigr)d\nu(\vsi)
\label{Ylast}
\end{equation}
and 
\begin{equation} 
\Phi(A,\lambda)=\frac{1}{N} \P(\vec{m}) F(A,\lambda).
\label{fN}
\end{equation}
When $A=[d,D],$ the set $\{R_{1,1}\in A\} = \Sigma^N$ and, therefore,
$$
F([d,D],\lambda) = \sum_{i\leq N} X_{k+1,i}
$$
where $X_{k+1,i}$ were defined in (\ref{Xcopies}). 
Using (\ref{Pprop2}),
\begin{equation}
\Phi([d,D],\lambda) = \P(\vec{m})X_{k+1} = X_0(\lambda).
\label{simple}
\end{equation}
Here we made the dependence of $X_0$ on $\lambda$ explicit
while keeping the dependence on other parameters
implicit. 
For simplicity of notations we will write
$$
F(\lambda): = F([d,D],\lambda) \mbox{ and }
\Phi(\lambda): = \Phi([d,D],\lambda)=X_0(\lambda).
$$
In these notations Theorem \ref{Th4} states that
there exists a sequence $\eps_N \to 0$ such that
\begin{equation}
\lim_{N\to\infty} \Phi(U_N,0) = \inf_{\lambda}(-\lambda u + \Phi(\lambda))
= \inf_{\lambda}(-\lambda u +X_0(\lambda)).
\label{reform}
\end{equation}
Let us define $\lambda(u)$ to be the point where the infimum in
(\ref{reform}) is achieved, i.e.
\begin{equation}
X_{0}(\lambda(u)) - \lambda(u) u =
\inf_{\lambda}\bigl(X_{0}(\lambda) - \lambda u\bigr).
\label{lambdau}
\end{equation}
The infimum is, indeed, achieved because of the following argument.
A function $-\lambda u + X_{0}(\lambda)$ is convex in $\lambda$
by H\"older's inequality. Lemma \ref{Lextra1} implies 
that $X_{0}(\lambda)\geq a(\lambda)\lambda$ for some function $a(\lambda)$ such that
$\lim_{\lambda\to+\infty}=D$ and $\lim_{\lambda\to-\infty}=d.$
Hence, for $d<u<D,$ the convex function
$-\lambda u + X_{0}(\lambda)\to+\infty$ as $\lambda\to\pm\infty$
and, therefore, it has a unique minimum.
The critical point condition for $\lambda(u)$ is 
\begin{equation}
\frac{\partial X_0}{\partial \lambda}(\lambda(u)) = u.
\label{lamt}
\end{equation}
Consider a fixed small enough $\eps>0$ such that
$d<u-\eps$ and $u+\eps<D$ and let 
$U_{\eps}= [u-\eps,u+\eps].$
Let us analyze $\Phi(U_{\eps},\lambda(u)).$
Consider a set $V$ equal to either
$[d,u-\eps]$ or $[u+\eps,D]$ and note that
$$
[d,D]=U_{\eps}\cup [d,u-\eps]\cup [u+\eps,D].
$$ 
We will start by proving an upper bound on $\Phi(V,\lambda(u)).$ 
We will only consider the case of $V=[u+\eps,D]$ 
since the case $V=[d,u-\eps]$ can be treated similarly.
Since
$$
-\gamma R_{1,1} \leq -\gamma (u+\eps) 
\,\,\,\mbox{ for }\,\,\,
R_{1,1}\in V \mbox{ and } \gamma\geq 0,
$$
we get that for $\gamma\geq 0,$
$$
F(V,\lambda(u)) \leq -N\gamma (u+\eps) + F(\lambda(u) + \gamma).
$$
Using (\ref{Pprop1}) and (\ref{Pprop3}) we get
\begin{equation}
\Phi(V,\lambda(u)) \leq U(\gamma) := -\gamma (u+\eps) + \Phi(\lambda(u)+\gamma)
= -\gamma (u+\eps) + X_0(\lambda(u)+\gamma).
\label{Ugamma}
\end{equation}
Setting $\gamma=0$ gives
$$
U(0)= \Phi(\lambda(u)) = X_0(\lambda(u)).
$$
Next, the right derivative of $U(\gamma)$ at zero is
$$
\frac{\partial U}{\partial\gamma}\Bigr|_{\gamma=0^{+}} = 
-(u+\eps)+\frac{\partial X_0}{\partial\lambda}(\lambda(u)) =
-(u+\eps)+u = -\eps
$$
using (\ref{lamt}). Finally, it follows from a 
tedious but straightforward computation which we will omit here that 
$$
\Bigl| \frac{\partial^2 U}{\partial\gamma^2}\Bigr|\leq L
$$
for some constant $L$ that depends only on the parameters of the model
$\xi$ and $\nu.$ 
Therefore, minimizing over 
$\gamma\geq 0$ in the right hand side of (\ref{Ugamma}) gives
\begin{equation}
\Phi(V,\lambda(u)) \leq \Phi(\lambda(u))-\frac{\eps^2}{L}.
\label{boundvu}
\end{equation}
The same bound holds for $V=[d,u-\eps].$
For $l\geq 1,$ let
$$
W_l = \exp m_l \bigl(F_{l+1}(\lambda(u)) - F_l(\lambda(u))\bigr)
$$
be defined by (\ref{PW}) with $F=F(\lambda(u)).$
Given a set $A\subseteq [d,D]$ let
$\la I(R_{1,1}\in A)\ra$ denote the Gibbs 
average defined by
$$
\la I(R_{1,1}\in A)\ra = 
\exp\bigl(F(A,\lambda(u)) - F(\lambda(u))\bigr).
$$
The following Proposition is the crucial step in the proof
of Theorem \ref{Th4}. This type of computation was invented by Talagrand
in \cite{T-P} in order to control the remainder terms in Guerra's 
interpolation and we will use this argument with the same purpose 
later in the paper as well. 

\begin{proposition}\label{Rprop1}
Assume that for $A\subseteq [d,D]$ and for some $\eps'>0$ we have
\begin{equation}
\Phi(A,\lambda(u))\leq \Phi(\lambda(u))-\eps'.
\label{fNeps}
\end{equation}
Then,
\begin{equation}
\e W_1\ldots W_k \la I(R_{1,1}\in A) \ra
\leq L\exp\Bigl(-\frac{N}{L}\Bigr),
\label{Vcontrol}
\end{equation}
where $L$ does not depend on $N.$ 
\end{proposition}
{\bf Proof.} 
The proof is based on the property of the Parisi functional described in Lemma
\ref{Plemlast}. For simplicity of notations let us assume that $m_1>0.$
The case when several elements of the sequence $\vec{m}$ are zeroes 
can be handled in exactly the same way. Let
$$
f=F(A,\lambda(u)) 
\mbox{ and } 
F=F(\lambda(u))
$$ 
so that the condition 
$f\leq F$ of Lemma \ref{Plemlast} is satisfied.
Lemma \ref{Plemlast} and (\ref{fNeps}) imply that
\begin{eqnarray*}
\e \log\e_1 W_1\ldots W_k \exp (f-F) 
&\leq& 
m_1 (\P(\vec{m})f - \P(\vec{m})F)
\\
&=& 
N m_1 (\Phi(A,\lambda(u)) - \Phi(\lambda(u)))
\leq - N m_1 \eps'.
\end{eqnarray*}
Let us consider a function
\begin{equation}
\phi(Z)= \log\e_1 W_1\ldots W_k \exp (f-F) \mbox{ where }
Z=(z_{i,0})_{i\leq N}.
\label{Rphit}
\end{equation}
We will show that $\phi(Z)$ is a Lipschitz function of $Z:$
\begin{equation}
|\phi(Z)-\phi(Z')| \leq L\sqrt{N}|Z-Z'|.
\label{Rlip}
\end{equation}
First of all, 
$$
\sum_{i\leq N} \sigma_i (z_{i,0} - z_{i,0}') \leq
|Z-Z'|\Bigl(\sum_{i\leq N} \sigma_i^2\Bigr)^{1/2} \leq \sqrt{ND}|Z-Z'|,
$$
using (\ref{diam1}). 
Definition (\ref{Ylast}) implies that for any set $A$ 
and any $\lambda,$
$$
|F(A,\lambda)(Z) - F(A,\lambda)(Z')|\leq \sqrt{ND}|Z-Z'|,
$$
where we made the dependence of $F(A,\lambda)$ on $Z$ explicit.
In particular, this holds for $f$ and $F.$
It is also clear from the properties (\ref{Pprop1}) and (\ref{Pprop3})
that iteration (\ref{Piter}) in the definition of the Parisi functional
preserves Lipschitz condition and, therefore,
$$
|F_l(Z) - F_l(Z')|\leq \sqrt{ND}|Z-Z'| \mbox{ for }
l\leq k.
$$
Using (\ref{GWT}) we can rewrite $\phi(Z)$ as
$$
\phi(Z)=\log \e_1
\exp \Bigl(f(Z) +
\sum_{1\leq l\leq k} (m_{l-1} - m_{l}) F_{l}(Z) \Bigr)
$$
and (\ref{Rlip}) is now obvious. Gaussian concentration of measure
(see, for example, Theorem 2.2.4 in \cite{SG}) implies that 
for any $t\geq 0,$
\begin{equation}
\p(|\phi(Z) - \e\phi| \geq Nt) \leq 2\exp(-N t^2/L).
\label{Rconc}
\end{equation}
In particular, with probability at least $1-2\exp(-N/L),$
$$
\phi\leq \e\phi + \frac{1}{2}Nm_1\eps' \leq -\frac{1}{2}Nm_1\eps'
$$
and, therefore, with probability at least $1-2\exp(-N/L),$
\begin{equation}
\e_1 W_1\ldots W_k \exp (f-F) \leq \exp(-N/L).
\label{probpart}
\end{equation}
Since $f\leq F,$
$$
\e_1 W_1\ldots W_k \exp (f-F)\leq \e_1 W_1\ldots W_k \leq 1
$$
using (\ref{Pden1}) which together with (\ref{probpart}) implies that
$$
\e W_1\ldots W_k \exp (f-F) \leq L\exp(-N/L).
$$
\qed

\begin{corollary}
For any $\eps>0$ we have
\begin{equation}
\Phi(U_{\eps},\lambda(u))\geq \Phi(\lambda(u)) - \delta_N
= X_0(\lambda(u)) - \delta_N,
\label{Udelta}
\end{equation}
where $\lim_{N\to\infty}\delta_N = 0.$
\end{corollary}
{\bf Proof.}
(\ref{boundvu}) and (\ref{Vcontrol}) imply that 
for $V$ equal to either $[d,u-\eps]$ or $[u+\eps,D]$
we have
$$
\e W_1\ldots W_k \la I(R_{1,1}\in V) \ra
\leq L\exp\Bigl(-\frac{N}{L}\Bigr)
$$
and, therefore, 
\begin{equation}
\e W_1\ldots W_k \la I(R_{1,1}\not \in U_{\eps}) \ra
\leq L\exp\Bigl(-\frac{N}{L}\Bigr).
\label{Udelta1}
\end{equation}
Suppose that (\ref{Udelta}) is not true which means that
for some positive $\eps'>0$ we have
$$
\Phi(U_{\eps},\lambda(u))\leq \Phi(\lambda(u)) - \eps',
$$
for some arbitrarily large  $N.$
Then again (\ref{boundvu}) and (\ref{Vcontrol}) would imply that 
$$
\e W_1\ldots W_k \la I(R_{1,1}\in U_{\eps}) \ra
\leq L\exp\Bigl(-\frac{N}{L}\Bigr).
$$
Combining with (\ref{Udelta1}) we would get
$$
1 = \e W_1\ldots W_k 
\leq L\exp\Bigl(-\frac{N}{L}\Bigr)
$$
and we arrive at contradiction.
\qed

In order to bound $\Phi(U_{\eps},\lambda(u))$ in terms of
$\Phi(U_{\eps},0),$ we can write
\begin{eqnarray*}
F(U_{\eps},\lambda(u))
&=&
\log\int_{\{R_{1,1}\in U_{\eps}\}}
\exp\Bigl(\sum_{i=1}^N \sigma_i \sum_{0\leq p\leq k} z_{i,p} 
+ N \lambda(u) R_{1,1} 
\Bigr) d\nu(\vsi) 
\\
&\leq &
N\lambda(u) u + N|\lambda(u)|\eps + F(U_{\eps},0).
\end{eqnarray*}
Using (\ref{Pprop1}) and (\ref{Pprop3}),
$$
\Phi(U_{\eps},\lambda(u))\leq \Phi(U_{\eps},0)+\lambda(u) u 
+ |\lambda(u)|\eps
$$
and, therefore,
\begin{equation}
\Phi(U_{\eps},0) \geq - \lambda(u) u + X_0(\lambda(u)) 
- |\lambda(u)|\eps - \delta_N.
\label{Ude}
\end{equation}
This implies that for any $\eps>0,$
$$
\liminf_{N\to\infty} \Phi(U_{\eps},0)\geq -\lambda(u)u + X_0(\lambda(u))
-L\eps.
$$
Clearly, this means that one can choose a sequence $\eps_N\to 0$
such that for  $U_N = U_{\eps_N},$
$$
\liminf_{N\to\infty} \Phi(U_{N},0)\geq -\lambda(u)u + X_0(\lambda(u)).
$$
Since a similar upper bound is obvious this finishes the proof of Theorem
\ref{Th4}.
\qed

\section{Reduction of the main results to apriori estimates.}\label{AprioriEst}

Now that we understood what happens at the end of Guerra's interpolation
we will turn to analyzing the remainder terms in the second line of
(\ref{Gphider}) and, in particular, the functional $\mu_r$ defined
in (\ref{Gmul}). First of all, for a function 
$h$ on $U_N\times U_N$ the definition of $\mu_r(h)$
can be written equivalently as follows.
Let $(z_{p}^j)$ for $j=1,2$ be two copies of the random vector
$(z_{p})$ defined in (\ref{z}) such that
\begin{equation}
z_{p}^1 = z_{p}^2 \mbox{ for }
p<r \mbox{ and } z_{p}^1,  z_{p}^2 \mbox{ are independent for }
p\geq r.
\label{Rzz}
\end{equation}
Let $\vec{z}^j=(z_{i,p}^j)_{i\leq N,p\leq k}$ where $(z_{i,p}^j)_{p\leq k}$
are independent copies of the vector $(z_p^j)$ for $i\leq N.$
Let $F^j$ be defined by (\ref{GFtlast}) in terms of $\vec{z}^j$
and let $W_l^j$ be defined by (\ref{PW}) in terms of $F^j.$
Let us consider the Hamiltonian
\begin{equation}
H_t(\vsi^1,\vsi^2)=\sum_{j=1,2}\Bigl(
\sqrt{t}H_N(\vsi^j) + \sqrt{1-t}\sum_{i\leq N} \sigma_i^j
\Bigl(\sum_{0\leq p\leq k} z_{i,p}^j\Bigr)
\Bigr)
\label{RHt}
\end{equation}
and define the Gibbs average $\la h \ra$ of $h$ by
$$
\la h \ra \exp(F^1 + F^2) = \int_{U_N\times U_N}
h(\vsi^1,\vsi^2) \exp H_t(\vsi^1,\vsi^2) d\nu(\vsi^1) d\nu(\vsi^2).
$$
Then, the definition (\ref{Gmul}) is equivalent to 
\begin{equation}
\mu_r(h)=\e W_1^1\ldots W_{r-1}^1 W_{r}^1 W_{r}^2 \ldots W_k^1 W_k^2 \la h\ra.
\label{eq1}
\end{equation}
We simply decoupled the measure $\gamma_r^{\otimes 2}$ by using independent
copies of $z_{i,p}$ for $p\geq r.$
Using Lemma \ref{Plem1} we can rewrite this in a more compact way. 
Let $F = F^1 + F^2$ and define $\vec{n}$ as in (\ref{Pprop42}), i.e.
\begin{equation}
n_{p} = \frac{m_{p}}{2} \mbox{ for } p< r
\mbox{ and }
n_{p} = m_{p} \mbox{ for } p\geq r.
\label{Rn}
\end{equation}
Lemma \ref{Plem1} then implies that for $l\geq r$
$$
W_l^1 W_l^2 = \exp m_l(F_{l+1}^1 - F_{l}^1) 
\exp m_l( F_{l+1}^2 - F_l^2)  = \exp n_l (F_{l+1} - F_{l})
$$
and for  $l<r$
$$
W_l^1 = \exp m_l (F_{l+1}^1 - F_l^1) = \exp n_l (F_{l+1} - F_l).
$$
Therefore, if we define
$$
W_l = \exp n_l(F_{l+1}-F_l),
$$
(\ref{eq1}) becomes
\begin{equation}
\mu_r(h)=\e W_1\ldots W_{k} \la h\ra.
\label{eq2}
\end{equation}
In particular, if $h = I(A)$ is an indicator of  
a measurable subset $A\subseteq U_N\times U_N$ and
\begin{equation}
F(A)= \log \int_{A} \exp H_t(\vsi^1,\vsi^2) d\nu(\vsi^1) d\nu(\vsi^2)
\label{RRfa}
\end{equation}
then, since
\begin{equation}
F= F^1 + F^2 =  
\log \int_{U_N\times U_N} \exp H_t(\vsi^1,\vsi^2) d\nu(\vsi^1) d\nu(\vsi^2),
\label{RRf}
\end{equation}
we get
$$
\mu_r(I(A))  = \e W_1\ldots W_{k} \exp(f -F).
$$
Lemma \ref{Plemlast} provided the methodology to control this expression. 
\begin{lemma}\label{Alem1}
Let $F(A)$ and $F$ be defined by (\ref{RRfa}) and (\ref{RRf}). 
If for some $\eps>0$ we have
\begin{equation}
\frac{1}{N}\e\P(\vec{n})F(A) \leq 2\varphi_N(t) - \eps
\label{RRcont}
\end{equation}
then for some constant $K$ independent of $N$ and the set $A,$
\begin{equation}
\mu_r(I(A))\leq K\exp(-N/K).
\label{RRcont2}
\end{equation}
\end{lemma}
{\bf Proof.} 
Using (\ref{Pprop4}) and (\ref{Gphi}), we can write
$$
\frac{1}{N}\e \P(\vec{n}) F = \frac{2}{N}\e \P(\vec{m}) F^1
= 2 \varphi_N(t),
$$
so that (\ref{RRcont}) can be written as
\begin{equation}
\e\P(\vec{n})F(A) - \e\P(\vec{n})F \leq - N\eps,
\label{RRcont1}
\end{equation}
which is the type of condition used in Lemma \ref{Plemlast}.
The main idea was already explained in detail in the proof of
Proposition \ref{Pprop1}. However, the function $\phi$ that was
defined in (\ref{Rphit}) was a function of the finite Gaussian
vector $Z$ in $\Reals^N$ with independent coordinates
which allowed us to use the classical
Gaussian concentration of measure inequality in (\ref{Rconc}).
Now, however, both $F(A)$ and $F$ depend on the entire Gaussian
process $H_N(\vsi)$ indexed by $\vsi\in\Sigma^N$ and the only
information that we specified about this process was the covariance
operator in (\ref{correlation}). Still, using the specific definition
of $F(A)$ and $F$ and (\ref{correlation}) one can prove the same 
concentration inequality as (\ref{Rconc}) but it would require a tedious
computation repeating the proof of (\ref{Rconc}) in \cite{SG}.
We will actually carry out this computation in a relatively easier
situation, below Lemma \ref{Lextra4}, so the idea will be clear and we will
omit this computation here. 
The proof becomes more transparent when the Hamiltonian is expressed
explicitly in terms of an i.i.d. Gaussian sequence. For example, one often
considers a Hamiltonian of the type
\begin{equation}
H_{N}(\vsi)=\sum_{p\geq 1}\frac{a_p}{N^{(p-1)/2}}
\sum_{i_1,\ldots,i_p} g_{i_1,\ldots,i_p} \sigma_{i_1}\ldots\sigma_{i_p},
\label{Ham}
\end{equation}
where $\vec{g} = (g_{i_1,\ldots,i_p})$ 
is a sequence of standard Gaussian random variables 
independent for all $p\geq 1$ and all $(i_1,\ldots,i_p).$
In this case,
$$
\frac{1}{N}\e H_N(\vsi^1) H_N(\vsi^2)=\xi(R_{1,2})
\mbox{ where }
\xi(x)=\sum_{p\geq 1} a_p^2 x^p.
$$ 
The sequence $(a_p)_{p\geq 1}$ should be such that $\xi(R_{1,2})$
is well defined for all $\vsi^1,\vsi^2$ and, comparing with 
(\ref{diam1}) this means that $\xi(D)<\infty.$ 
It is easy to check that
for two sequences $\vec{g}$ and $\vec{g}'$ we have
$$
|H_N(\vsi)(\vec{g}) - H_N(\vsi)(\vec{g}')| \leq \sqrt{N \xi(D)}
|\vec{g} - \vec{g}'|
$$
and since inequality (\ref{Rconc}) is dimension independent, 
it applies to a sequence $\vec{g}$ and
the rest of the proof repeats the proof of Proposition \ref{Pprop1}.
\qed

In Lemma \ref{Lextra2} we explained why $\P(\xi,u)$ is an upper bound
on local free energy and the main reason was that the remainder terms
in Guerra's interpolation were nonnegative. In order to show that this bound
is exact in the limit we must show that these remainder terms are small
along the interpolation for some choices of the parameters $k,\vec{m}$
and $\vec{q}$ and that $\varphi_N(t)$ can be approximated by 
\begin{equation}
\psi(t) = \varphi_0(\vec{m},\vec{q}) - \frac{t}{2}\sum_{1\leq l\leq k}
m_l(\theta(q_{l+1}) - \theta(q_l)),
\label{RRpsi}
\end{equation}
where $\varphi_0(\vec{m},\vec{q})$ was defined in (\ref{philimit}).
It is also clear that these parameters should approximate the
infimum in the definition (\ref{Pu}) of $\P(\xi,u).$

{\bf Definition.} We will call a vector $(k,\vec{m},\vec{q},\lambda)$
an $\eps$-minimizer if
\begin{equation}
\P_k(\vec{m},\vec{q},\lambda,u)\leq \P(\xi,u) + \eps
\,\,\,\mbox{ and }\,\,\, (\vec{m},\vec{q},\lambda) 
\mbox{ is the minimizer of (\ref{Pk}). }
\label{mini}
\end{equation}
For any $v\in[-D,D],$ let us define a set
\begin{equation}
A(v) = \Bigl\{ |R_{1,2} - v| \leq \frac{1}{N}\Bigr\}\bigcap (U_N\times U_N).
\label{UN}
\end{equation}
The following apriori estimate will allow us to control the remainder
terms in Guerra's interpolation.

\begin{theorem}\label{Apriori}
For any $t_0<1$ there exists $\eps>0$ that depends on $t_0,\xi,\nu, u$ only
such that if (\ref{mini}) holds then for $t\leq t_0$ and
for large enough $N,$
\begin{equation}
\frac{1}{N}\e \P(\vec{n})F(A(v))\leq 2\psi(t) - \frac{(v-q_r)^2}{K}+\R,
\label{apri}
\end{equation}
where $K$ is a constant independent of $N, t$ and $v$
and $|\R|\leq a_N$ for a sequence $(a_N)$ independent of $t$ and $v$ and 
$\lim_{N\to\infty} a_N = 0.$
\end{theorem}

In the case of the replica symmetric region of Theorem \ref{RST1}
the condition on $\eps$-minimizer is replaced by the condition of
stability to replica symmetry breaking fluctuations defined in 
(\ref{RSRScond}).

\begin{theorem}\label{RSApriori}
Suppose that all functions are defined in terms of parameters in (\ref{RSpar})
and that (\ref{RSstable}) and (\ref{RSRScond}) hold. Then
for any $t_0<1$  and for any $t\leq t_0,$ for large enough $N,$
\begin{equation}
\frac{1}{N}\e \P(\vec{n})F(A(v))\leq 2\psi(t) - \frac{(v-q)^2}{K}+\R,
\label{RSapri}
\end{equation}
where $K$ is a constant independent of $N, t$ and $v$
and $|\R|\leq a_N$ for a sequence $(a_N)$ independent of $t$ and $v$ and 
$\lim_{N\to\infty} a_N = 0.$
\end{theorem}

The proof of these apriori estimates will be postponed until Appendix A.
First, let us show how they imply Theorems \ref{Th1} and \ref{RST1}.

{\bf Proof of Theorem \ref{Th1}.}
Given $t_0<1$ let us take $\eps>0$ as in Theorem \ref{Apriori}
and let $(k,\vec{m},\vec{q},\lambda)$ be an $\eps$-minimizer
defined by (\ref{mini}). Clearly, in this case,
\begin{equation}
\psi(1) = \P_k(\vec{m},\vec{q},\lambda,u) \,\,\mbox{ and }\,\,
|\psi(1)- \P(\xi,u)| \leq \eps.
\label{psi1}
\end{equation}
Let us take $K$ as in (\ref{apri}) and for $\eps_1>0$ define a set
$$
\V = \{v\in[-D,D] : (v-q_r)^2\geq 2K(\psi(t) - \varphi_N(t)) + 2K \eps_1\}.
$$
For any $v\in\V,$ (\ref{apri}) implies that
$$
\frac{1}{N}\e \P(\vec{n})F(A(v)) \leq
2\psi(t) - \frac{(v-q_r)^2}{K}+\R \leq 2\varphi_N(t) - \eps_1
$$
for large enough $N.$ Everywhere below let $L$ denote a constant
that might depend on $\eps_1$ and $K$ denote a constant independent
of $\eps_1.$ Applying Lemma \ref{Alem1}, we get 
\begin{equation}
\mu_r(I(A(v)))\leq L\exp(-N/L)
\label{muav}
\end{equation}
and the constant $L$ here does not depend on $v$. Let us consider
a set
$$
A = \Bigl\{(R_{1,2}-q_r)^2\geq 2K(\psi(t) - \varphi_N(t)) +2K \eps_1\Bigr\}
\bigcap (U_N\times U_N).
$$
We can choose the points $v_1,\ldots, v_M \in \V$
with $M\leq K N$ such that $A\subseteq \bigcup_{i\leq M} A(v_i)$
and (\ref{muav}) implies that
\begin{equation}
\mu_r(I(A))\leq LN\exp(-N/L)\leq L\exp(-N/L).
\label{muA}
\end{equation}
Using the definition of $\psi$ in (\ref{RRpsi}) and (\ref{Gphider}),
\begin{equation}
(\psi(t)-\varphi_N(t))' = 
\frac{1}{2}\sum_{1\leq l\leq k} (m_{l} - m_{l-1})
\mu_l\bigl(
\xi(R_{1,2})-R_{1,2}\xi'(q_l) +\theta(q_l) 
\bigr) +{\cal R}.
\label{derdif}
\end{equation}
Since the second derivative of $\xi$ is bounded on $[-D,D],$
$$
\xi(R_{1,2})-R_{1,2}\xi'(q_l) +\theta(q_l) \leq K(R_{1,2} - q_l)^2
$$
and because on the complement $A^c$ of $A$ we have 
$$
(R_{1,2}-q_r)^2\leq 2K(\psi(t) - \varphi_N(t)) +2K \eps_1,
$$
(\ref{muA}) for $r=l$ implies that
\begin{eqnarray*}
&&
\mu_l(\xi(R_{1,2})-R_{1,2}\xi'(q_l) +\theta(q_l)) \leq 
K\mu_l((R_{1,2} - q_l)^2) 
\\
&&
\leq 
K\Bigl(
(\psi(t) -\varphi_N(t)) + \eps_1 +  \mu_l(I(A))
\Bigr)
\leq 
K(\psi(t) -\varphi_N(t)) + K\eps_1 + L\exp(-N/L).
\end{eqnarray*}
(\ref{derdif}) now implies that
$$
(\psi(t)-\varphi_N(t))' \leq K(\psi(t)-\varphi_N(t)) + K\eps_1 + L\exp(-N/L).
$$
Since by Theorem \ref{Th4}, $\lim_{N\to\infty} \varphi_N(0) = \psi(0),$
solving this differential inequality and then letting $N\to\infty$
and $\eps_1\to 0$ implies that 
\begin{equation}
\lim_{N\to\infty}\varphi_N(t) = \psi(t) \mbox{ for } t\leq t_0.
\label{limee}
\end{equation}
The derivatives $\psi'(t)$ and $\varphi_N'(t)$ are both bounded,
which is apparent from (\ref{RRpsi}) and (\ref{Gphider}),
and we get
$$
\limsup_{N\to\infty} |\varphi_N(1) -\psi(1)|\leq K(1-t_0).
$$
Using (\ref{psi1}),
$$
\limsup_{N\to\infty} |\varphi_N(1) -\P(\xi,u)|\leq K(1-t_0) + \eps
$$
and letting $\eps\to 0$ and $t_0\to 1$ finishes the proof of 
Theorem \ref{Th1}.
\qed

{\bf Proof of Theorem \ref{RST1}.} 
The proof of the the first part follows from Theorem \ref{RSApriori}
in exactly the same way as Theorem \ref{Th1} follows from Theorem
\ref{Apriori}. The uniqueness of $(q,\lambda)$ follows from the
following simple argument. (\ref{limee}) implies that for $t\leq t_0$ 
$$
\lim_{N\to\infty}\mu_1((R_{1,2} - q)^2)=0
$$ 
and since the definition of $\mu_1$ does not depend on $\lambda$
this $q$ must be unique. Since $\P_1(q,\lambda)$ is convex in $\lambda,$
this implies the uniqueness of $\lambda.$
\qed

\section{ Computing global free energy.}\label{SecGlobal}

In this section we will prove Theorem \ref{Th2} which will follow 
from Theorem \ref{Th1} and Gaussian concentration of measure.
Let us start by proving the following concentration inequality.

\begin{lemma}\label{Lextra4}
For any measurable subset $\Omega\subseteq \Sigma^N$ 
let us consider a r.v.
\begin{equation}
X=\log\int_{\Omega}\exp H_N(\vsi) d\nu(\vsi).
\label{Xconc}
\end{equation}
Then, for any $t\geq 0,$
\begin{equation}
\p\Bigl(|X - \e X|\geq 2\sqrt{LNt}\Bigr)\leq 2\exp(-t),
\label{concineq}
\end{equation}
where $L=\max\{\xi(x): x\in [d,D]\} + c(N).$
\end{lemma}
{\bf Proof.}
The proof is a simple modification of Theorem 2.2.4 in \cite{SG}.
Unfortunately, Lemma \ref{Lextra4} does not fall into the framework
of Theorem 2.2.4 in \cite{SG} directly, but the same argument still
works if we utilize the particular definition of $X$ and
the covariance structure (\ref{correlation}) of the Hamiltonian 
$H_N(\vsi).$

Let $H_N^1$ and $H_N^2$ be two independent copies of the Hamiltonian
$H_N.$ For $t\in [0,1],$ we define,
\begin{equation}
F_j = \log \int_{\Omega}\exp(\sqrt{t} H_N^j(\vsi)
+\sqrt{1-t}H_N(\vsi)) d\nu(\vsi) 
\label{Pfi}
\end{equation}
for $j=1,2$ and let $F=F_1+F_2.$ For $s\geq 0,$ let
$$
\varphi(t) = \e \exp s(F_2 -F_1).
$$
If we define
$$
H_t(\vsi^1,\vsi^2) = 
\sqrt{t} (H_N^1(\vsi^1)+ H_N^2(\vsi^1))
+\sqrt{1-t}(H_N(\vsi^1)+H_N(\vsi^2))
$$
then  straightforward computation as in Theorem \ref{GTh3} gives
\begin{equation}
\varphi'(t) = -N s^2 
\e \exp s(F_1 - F_2) \exp(-F)
\int_{\Omega^2} \zeta(\vsi^1,\vsi^2)\exp H_t(\vsi^1,\vsi^2)
d\nu(\vsi^1)d\nu(\vsi^2),
\label{pfider}
\end{equation}
where $\zeta$ was defined in (\ref{Gzeta}).
Since $|\zeta(\vsi^1,\vsi^2)|\leq L,$ we get
\begin{equation}
\varphi'(t) \leq LN s^2 
\e \exp s(F_1 - F_2) \exp(-F)
\int_{\Omega^2} \exp H_td\nu d\nu
= LNs^2 \varphi(t).
\label{diff}
\end{equation}
By construction $\varphi(0)=1,$ so that (\ref{diff}) implies that
$\varphi(1)\leq \exp NLs^2.$ 
On the other hand, by construction, $\varphi(1)=\e \exp s(X-X'),$
where $X$ is defined in (\ref{Xconc}) and $X'$ is an independent
copy of $X$ and, thus,
$$
\e \exp s(X-X')\leq \exp NLs^2.
$$
By Jensen's inequality, this implies that,
$\e \exp s(X-\e X)\leq \exp NLs^2$ and
using Markov's inequality we get that for $t>0,$
$$
\p\Bigl( X-\e X\geq t\Bigr)\leq \inf_{s\geq 0} \exp(NLs^2 - st)
=\exp\Bigl(-\frac{t^2}{4NL}\Bigr).
$$
Obviously, a similar inequality can be written for $\e X -X$ and, 
therefore,
$$
\p\Bigl( |X-\e X|\geq t\Bigr)\leq 2\exp\Bigl(-\frac{t^2}{4NL}\Bigr).
$$
This is equivalent to (\ref{concineq}).
\qed

The proof of Theorem \ref{Th2} will be based on the concentration
inequality of Lemma \ref{Lextra4} and the following result.

\begin{theorem}\label{ThE}
For $\eps_N=N^{-1}$ there exists a set 
$\{u_1,\ldots,u_M\}\subseteq [d,D]$ of cardinality
$M\leq LN$ such that
\begin{equation}
[d,D]\subseteq \bigcup_{i\leq M}[u_i-\eps_N,u_i+\eps_N]
\label{unio}
\end{equation}
and such that for all $i\leq M,$
\begin{equation}
F_N(u_i,\eps_N) \leq \sup_{d\leq u\leq D}\P(\xi,u) + \R,
\label{Eps}
\end{equation}
where $|\R|\leq c_N + L\eps_N.$
\end{theorem}

We will first show that Theorem \ref{Th2} follows from
Theorem \ref{Th1} and Theorem \ref{ThE}.

{\bf Proof of Theorem 2.}
It follows from the definitions that 
for any $d\leq u\leq D$ and any sequence $(\eps_N),$
$F_N(u,\eps_N)\leq F_N$ and, therefore, Theorem \ref{Th1} implies
$$
\liminf_{N\to\infty} F_N\geq \sup_{d\leq u\leq D}\P(\xi,u).
$$
To prove Theorem \ref{Th2} we need to show that
$$
\limsup_{N\to\infty} F_N\leq \sup_{d\leq u\leq D}\P(\xi,u).
$$ 
Suppose not. Then for some $\eps>0,$
$$
\limsup_{N\to\infty} F_N > \sup_{d\leq u\leq D}\P(\xi,u) +\eps.
$$
To simplify the notations, instead of considering a subsequence
of $N$ we simply assume that for $N$ large enough we have
\begin{equation}
F_N \geq \sup_{d\leq u\leq D}\P(\xi,u) +\eps.
\label{1}
\end{equation}
On the other hand, let $\eps_N=N^{-1}$ and consider a set of points
$\{u_1,\ldots,u_M\}$ as in Theorem \ref{ThE}, so that for
$N$ large enough for all $i\leq M$,
\begin{equation} 
F_{N}(u_i,\eps_N) \leq \sup_{d\leq u\leq D}\P(\xi,u) + \frac{\eps}{2}.
\label{2}
\end{equation}
Lemma \ref{Lextra4} implies that for any $t>0$ with probability at least 
$1-LN\exp(-Nt^2/4L)$ we have
\begin{equation}
F_N\leq \frac{1}{N}\log Z_N + t\,\,\,
\mbox{ and }\,\,\,
\frac{1}{N}\log Z_N(u_i,\eps_N) \leq F_N(u_i,\eps_N)+t
\label{3}
\end{equation}
for all $i\leq M.$
Therefore, (\ref{1}), (\ref{2}) and (\ref{3}) imply that
for $N$ large enough with probability at least $1-LN\exp(-Nt^2/4L)$ we have
\begin{eqnarray*}
&&
\frac{1}{N}\log Z_N(u_i,\eps_N) 
\leq 
F_N(u,\eps_N)+t
\leq F_N +t -\frac{\eps}{2} 
\leq \frac{1}{N}\log Z_N +2t -\frac{\eps}{2}.
\end{eqnarray*}
Taking $t=\eps/K$ we get that for $N$ large enough 
with probability at least $1-LN\exp(-N\eps^2/K)$ 
we have
\begin{eqnarray*}
&&
\forall i\leq M, \,\,\,\,
\log Z_N(u_i,\eps_N)
\leq
\log Z_N -\frac{N\eps}{K}.
\end{eqnarray*}
This yields that for $i\leq M$ the Gibbs measure
$$
G_N(\{R_{1,1}\in[u_i-\eps_N,u_i+\eps_N]\})
\leq \exp\Bigl(- \frac{N\eps}{K}\Bigr)
$$
and, therefore, by (\ref{unio}),
$$
1=G_N(\Sigma^N)=G_N(\{R_{1,1}\in[d,D]\})
\leq LN \exp\Bigl(-\frac{N\eps}{K}\Bigr),
$$
which is impossible for large $N.$
\qed

The statement of Theorem \ref{ThE} is intuitively obvious considering
(\ref{oneside}). We, basically, need to show that the term $u_N(\lambda)$
in (\ref{oneside}) can be substituted by $u$ in a controlled manner.
This will be based on three technical lemmas.
To formulate the first lemma, it will be convenient to think of the pair
$(\vec{m}, \vec{q})$ in terms of the function $m=m(q)$ defined 
in (\ref{remark}).

The following Lemma is the analogue of a well known continuity
property of $X_0(\vec{m},\vec{q},\lambda)$ with respect 
to the functional order parameter $m(q)$
in the SK model (the statement can be found in \cite{Guerra} 
and the proof is given in \cite{T-PM})
and since its proof is exactly the same we will
not reproduce it here.

\begin{lemma}\label{LemmaL1}
For any $\lambda$ and for any functions $m(q)$ and $m'(q)$
defined by (\ref{remark}) and corresponding to pairs $(\vec{m},\vec{q})$
and $(\vec{m}',\vec{q}')$ we have
\begin{equation}
|X_0(\vec{m},\vec{q},\lambda)-X_0(\vec{m}',\vec{q}',\lambda)|
\leq
L\int_0^D |m(q)-m'(q)|dq,
\label{L1}
\end{equation}
where the constant $L$ depends on $\xi$ and $D$ only.
\end{lemma}

Next, we will describe several properties of the function
$\P_k(\vec{m},\vec{q},\lambda,u)$ defined in (\ref{Pk}).
For any $k,$ any vectors $\vec{m}, \vec{q}$ and any
$d\leq u\leq D$ let $\lambda(u)$ be defined by (\ref{lambdau}). 

\begin{lemma}\label{Lem1}
For any $\delta>0$ there exists a constant $\Lambda(\delta)$
such that for any vectors $\vec{m},\vec{q}$ 
and any $u\in [d+\delta, D-\delta]$ we have
\begin{equation}
|\lambda(u)|\leq \Lambda(\delta).
\label{uopt}
\end{equation}
\end{lemma}
{\bf Proof.}
For a fixed $\vec{m}$ and $\vec{q}$ minimizing 
$\P_k(\vec{m},\vec{q},\lambda,u)$ over $\lambda$
is equivalent to minimizing $-\lambda u + X_0(\vec{m},\vec{q},\lambda)$
over $\lambda$. By Lemma \ref{Lextra1},
$$
-\lambda u + X_0(\vec{m},\vec{q},\lambda) \geq 
\lambda (a(\lambda) - u),
$$ 
where $a(\lambda)$ satisfies (\ref{alim}). 
Therefore, for any $\delta>0$ the right hand side goes to infinity
uniformly over $d+\delta\leq u\leq D-\delta$ and, thus,
the left hand side goes to infinity uniformly over
$\vec{m},\vec{q}$ and $d+\delta\leq u\leq D-\delta.$
This, obviously, implies that there exists $\Lambda(\delta)$
such that (\ref{uopt}) holds.
\qed

\begin{lemma}\label{Lem2}
There exists $\delta>0$ such that for all
$\vec{m}$ and $\vec{q}$ we have
\begin{equation}
\lambda(u)<0 \,\,\mbox{ for }\,\, u\in [d,d+\delta),
\,\,\,
\lambda(u)>0 \,\,\mbox{ for}\,\,  u\in(D-\delta,D].
\label{del}
\end{equation}
\end{lemma}
{\bf Proof.} 
Using (\ref{Pder}) and (\ref{lamt}),
\begin{equation}
\frac{\partial X_0}{\partial \lambda}=
\e W_1\ldots W_k \la\sigma^2\ra=u,
\label{there}
\end{equation}
where for a function $h(\sigma),$  $\la h\ra$ is defined by,
\begin{equation}
\la h \ra \exp X_{k+1} = \int_{\Sigma} h(\sigma)\exp 
\Bigl(
\sigma \sum_{p\leq k} z_p +\lambda\sigma^2 
\Bigr) d\nu(\sigma).
\label{have}
\end{equation}
Since $X_0$ is convex in $\lambda,$ $\partial X_0/\partial \lambda$ 
is increasing in $\lambda$ and, therefore, (\ref{there}) implies
that $\lambda=\lambda(u)$ is nondecreasing in $u.$
Therefore, in order to prove (\ref{del}),
it is enough to show that for some $\delta>0$
if $\lambda(u)=0$ then $u\in[d+\delta,D-\delta].$
We will only prove that $\lambda(u)=0$ implies that $u\geq d+\delta$ 
as the case $u\leq D-\delta$ is quite similar.

We set $\lambda = 0$ and denote $Z=\sum_{p\leq k}z_p.$ 
Given $\gamma>0,$ we can write
\begin{equation}
\e W_1\ldots W_k \la\sigma^2\ra \geq
\e W_1\ldots W_k \la\sigma^2\ra I(|Z|\leq \gamma).
\label{low1}
\end{equation}
First of all, let us show that if $|Z|\leq \gamma$ then
\begin{equation}
\la\sigma^2\ra \geq d + \frac{1}{L}\exp(-L\gamma),
\label{siglow}
\end{equation}
where a constant $L$ does not depend on $\gamma$
(but it depends on other parameters of the model such as measure $\nu$). 
To show this, let us first note that if $|Z|\leq \gamma$ we have
\begin{equation}
\exp X_{k+1} = \int_{\Sigma} \exp (\sigma Z) d\nu(\sigma)
\leq L\exp(L\gamma),
\label{qe1}
\end{equation}
using (\ref{diam1}).
Next,
\begin{eqnarray}
\la\sigma^2\ra \exp X_{k+1} 
&=& 
\int_{\Sigma} \sigma^2\exp (\sigma Z)d\nu(\sigma)
\nonumber
\\
&=&
d \exp X_{k+1} +
\int_{\Sigma} (\sigma^2-d)\exp (\sigma Z) d\nu(\sigma)
\nonumber
\\
&\geq&
d\exp X_{k+1}
+ \frac{1}{L}\exp(-L\gamma),
\label{qe2}
\end{eqnarray}
where the last inequality is obtained by restricting the integral
to the set $\{\sigma^2\in [(D+d)/2,D]\}$ of positive measure $\nu$
by (\ref{diam1}). Combining (\ref{qe1}) and (\ref{qe2}) yields 
(\ref{siglow}). Plugging (\ref{siglow}) into (\ref{low1})
implies
\begin{eqnarray}
\e W_1\ldots W_k \la\sigma^2\ra 
&\geq&
\Bigl(d + \frac{1}{L}\exp(-L\gamma)\Bigr)
\e W_1\ldots W_k I(|Z|\leq \gamma)
\nonumber
\\
&=&
\Bigl(d + \frac{1}{L}\exp(-L\gamma)\Bigr)
\Bigl(1- \e W_1\ldots W_k I(|Z| > \gamma)\Bigr).
\label{low4}
\end{eqnarray}
Using  H\"older's inequality, we can bound
\begin{eqnarray}
\e W_1\ldots W_k I(|Z| > \gamma) 
&\leq& 
\bigl(\p(|Z|>\gamma)\bigr)^{1/2}
\bigl(\e (W_1\ldots W_k)^2\bigr)^{1/2} 
\nonumber
\\
&\leq&
L \exp\Bigl(-\frac{\gamma^2}{L}\Bigr)
\bigl(\e (W_1\ldots W_k)^2\bigr)^{1/2},
\label{low3}
\end{eqnarray}
since $Z$ is a Gaussian r.v. with variance $\xi(u)$
uniformly bounded for all $d\leq u\leq D.$
As in (\ref{GWT}) we can write
$$
W_1\ldots W_k = \exp X_{k+1} \prod_{l\leq k} (\exp(-X_l))^{m_l-m_{l-1}} 
$$
and by H\"older's inequality,
\begin{equation}
\e (W_1\ldots W_k)^2 \leq \bigl(\e \exp (4X_{k+1})\bigr)^{1/2} 
\prod_{l\leq k} \bigl(\e \exp(-4X_l)\bigr)^{(m_l - m_{l-1})/2}.
\label{Wus}
\end{equation}
The first factor on the right hand side can be estimated as follows,
\begin{eqnarray*}
\e\exp(4X_{k+1})
&=&
\e\Bigl(
\int_{\Sigma}\exp (\sigma Z) d\nu(\sigma)
\Bigr)^4
\leq
\e\int_{\Sigma}\exp (4\sigma Z) d\nu(\sigma)
\\
&=&
\int_{\Sigma}\exp (8\sigma^2 \xi'(u)) d\nu(\sigma)
\leq L.
\end{eqnarray*}
Next, since
$$
X_l=\frac{1}{m_l}\log\e_l\exp m_l X_{l+1} \geq
\e_l X_{l+1},
$$
repeating this over $l$  yields that
$X_l\geq \e_l X_{k+1}.$
Therefore,
\begin{eqnarray*}
&&
\e \exp(-4X_l) 
\leq
\e\exp(-4\e_l X_{k+1}) \leq
\e \exp(-4X_{k+1})
=
\e\Bigl(
\int_{\Sigma}\exp (\sigma Z) d\nu(\sigma)
\Bigr)^{-4}
\\
&&
\leq
\e \int_{\Sigma}\exp (-4\sigma Z) d\nu(\sigma)
=
\int_{\Sigma}\exp (8\sigma^2 \xi'(u)) d\nu(\sigma)
\leq L.
\end{eqnarray*}
Plugging all these estimates into (\ref{Wus}) gives
$$
\e (W_1\ldots W_k)^2 \leq L
$$
and (\ref{low3}) implies
$$
\e W_1\ldots W_k I(|Z| > \gamma) 
\leq
L \exp\Bigl(-\frac{\gamma^2}{L}\Bigr).
$$
Finally, (\ref{low4}) implies that
\begin{eqnarray*}
u=\e W_1\ldots W_k \la\sigma^2\ra 
&\geq& 
\Bigl(d + \frac{1}{L}\exp(-L\gamma )\Bigr)
\Bigl(1- L \exp\Bigl(-\frac{\gamma^2}{L}\Bigr)\Bigr).
\\
&\geq& 
d + \frac{1}{L}\exp(-L\gamma )
-L\exp\Bigl(-\frac{\gamma^2}{L}\Bigr).
\end{eqnarray*}
Taking $\gamma$ large enough gives
$u\geq d+\delta$ for some $\delta>0.$
As we have already mentioned above one can
similarly show that $\lambda=0$ implies that
$u\leq D-\delta$ for some $\delta>0$ and
this finishes the proof of Lemma.
\qed

We are now ready to prove Theorem \ref{ThE}.

{\bf Proof of Theorem \ref{ThE}.}
Let us start with equation (\ref{oneside}) that states that
for any $\lambda, \vec{m}$ and $\vec{q}$ such that $q_{k+1}=u,$
\begin{equation}
F_N(u,\eps_N) \leq
-\lambda u_N(\lambda) + X_0(\vec{m},\vec{q},\lambda)
-\frac{1}{2}\sum_{1\leq l\leq k} m_l(\theta(q_{l+1})-\theta(q_l))
+\R,
\label{onesideG}
\end{equation}
where $u_N(\lambda)$ is defined in (\ref{uenn}) and
$|\R|\leq c_N + L\eps_N.$ 
Let us note that the definition (\ref{uenn}) implies that
\begin{equation}
|u_N(\lambda) - u|\leq \eps_N.
\label{dee}
\end{equation}
Let us take $\delta>0$ as in Lemma \ref{Lem2}. 
Lemma \ref{Lem1} implies that for all $d+\delta\leq u\leq D-\delta$
we have $|\lambda(u)|\leq \Lambda(\delta).$
Moreover, (\ref{onesideG}) and (\ref{dee}) imply that for
any $d+\delta\leq u\leq D-\delta$ and $|\lambda|\leq \Lambda(\delta)$
we have
\begin{eqnarray}
F_N(u,\eps_N) 
&\leq&
-\lambda u + X_0(\vec{m},\vec{q},\lambda)
-\frac{1}{2}\sum_{1\leq l\leq k} m_l(\theta(q_{l+1})-\theta(q_l))
+ \Lambda(\delta)\eps_N +\R
\nonumber
\\
&=&
\P_k(\vec{m},\vec{q},\lambda,u) + \R,
\label{aib}
\end{eqnarray}
where again $|\R|\leq c_N +L\eps_N.$
The fact that for $d+\delta\leq u\leq D-\delta$ we have
$|\lambda(u)|\leq \Lambda(\delta)$ means in this case that
for a fixed $(\vec{m},\vec{q})$ minimizing $\P_k(\vec{m},\vec{q},\lambda,u)$ 
over $\lambda$ is equivalent to minimizing it over $|\lambda|\leq \Lambda(\delta)$
and, therefore, minimizing $\P_k(\vec{m},\vec{q},\lambda,u)$ 
over $(\vec{m},\vec{q},\lambda)$ is also equivalent to minimizing 
it over $\vec{m},\vec{q}$ and $|\lambda|\leq \Lambda(\delta).$
Therefore, (\ref{aib}) implies that for $d+\delta\leq u\leq D-\delta,$
\begin{equation}
F_N(u,\eps_N) \leq \P(\xi,u) + \R
\leq \sup_{d\leq u\leq D}\P(\xi,u) + \R,
\label{almost}
\end{equation}
where $|\R|\leq c_N + L\eps_N.$

Next, let $u$ be such that $u+\eps_N < d+\delta.$
Let us consider arbitrary $k'\geq 1$ and 
arbitrary $\vec{m}'$ and $\vec{q}'$ such that
$q_{k'+1}' = u_N(\lambda).$ 
By (\ref{dee}), there exist $\vec{m}$ and $\vec{q},$
maybe, with different parameter $k,$ such that
$q_{k+1}=u$ and 
\begin{equation}
\int |m(q) - m'(q)| dq \leq |u-u_N(\lambda)| \leq \eps_N,
\label{L12}
\end{equation}
where functions $m(q)$ and $m'(q)$ are defined in (\ref{remark})
in terms of $(\vec{m},\vec{q})$ and $(\vec{m}',\vec{q}')$ 
correspondingly. This can be achieved by simply assigning
$m(q)=1$ for $q$ between $u$ and $u_N(\lambda)$ and, otherwise,
letting $m(q)=m'(q).$ Then Lemma \ref{LemmaL1} implies that
$$
|X_0(\vec{m},\vec{q},\lambda)-X_0(\vec{m}',\vec{q}',\lambda)|
\leq L\eps_N.
$$
Also, obviously, condition (\ref{L12}) implies that
$$
\Bigl|
\frac{1}{2}\sum_{1\leq l\leq k} m_l(\theta(q_{l+1})-\theta(q_l))
-\frac{1}{2}\sum_{1\leq l\leq k'} m_l'(\theta(q_{l+1}')-\theta(q_l'))
\Bigr| \leq L\eps_N.
$$
Therefore, (\ref{onesideG}) implies that
for any $k', \vec{m}'$ and $\vec{q}'$ such that
$q_{k'+1}' = u_N(\lambda),$
\begin{eqnarray}
F_N(u,\eps_N) 
&\leq&
-\lambda u_N(\lambda) + X_0(\vec{m}',\vec{q}',\lambda)
-\frac{1}{2}\sum_{1\leq l\leq k'} m_l'(\theta(q_{l+1}')-\theta(q_l'))
+\R
\nonumber
\\
&=&
\P_{k'}(\vec{m}',\vec{q}',\lambda,u_N(\lambda)) + \R,
\label{onesideG2}
\end{eqnarray}
where again $|\R|\leq c_N + L\eps_N.$
Since now $u+\eps_N<d+\delta<D,$
for $\lambda<0$ the definition (\ref{uenn}) implies that
$u_N(\lambda) = u+\eps_N$ and, therefore, for
$\lambda<0$ we have
\begin{equation}
F_N(u,\eps_N)\leq
\P_{k'}(\vec{m}',\vec{q}',\lambda,u+\eps_N) + \R.
\label{aibo}
\end{equation}
Since $u+\eps_N<d+\delta,$ Lemma \ref{Lem2} implies that
for any $\vec{m}'$ and  $\vec{q}',$
$\lambda(u)<0$ which means that 
for fixed $(\vec{m}',\vec{q}')$ minimizing $\P_{k'}(\vec{m}',\vec{q}',\lambda,u+\eps_N)$ 
over $\lambda$ is equivalent to minimizing it over $\lambda<0$
and, therefore, minimizing $\P_{k'}(\vec{m}',\vec{q}',\lambda,u+\eps_N)$ 
over $(\vec{m}',\vec{q}',\lambda)$ is also equivalent to minimizing 
it over $\vec{m}',\vec{q}'$ and $\lambda <0.$
Hence, (\ref{aibo}) yields that it $u+\eps_N<d+\delta$ then,
\begin{equation}
F_N(u,\eps_N) \leq \P(\xi,u+\eps_N) + \R \leq
\sup_{d\leq u\leq D}\P(\xi,u) +\R.
\end{equation}
Similarly, one can show that this holds for $u$ such that $u-\eps_N>D-\delta$
and, combining this with (\ref{almost}), we showed that for all
$u$ in the set
\begin{equation}
\{u+\eps_N<d+\delta\}\cup\{d+\delta\leq u\leq D-\delta\}\cup\{u-\eps_N>D-\eps_N\} 
\label{setU}
\end{equation}
we have
$$
F_N(u,\eps_N) \leq
\sup_{d\leq u\leq D}\P(\xi,u) +\R.
$$
It is obvious that for $\eps_N=N^{-1}$ the set (\ref{setU}) contains 
$\eps_N$-net of the interval $[d,D]$ of cardinality $LN.$
This finishes the proof of Theorem \ref{ThE}.
\qed

{\bf Acknowledgment.} I would like to thank David Sherrington  
for suggesting the topic of this research.

\appendix

\section{Proof of the apriori estimates. }

We will prove the apriori estimates of Section \ref{AprioriEst} in
several steps. In Section \ref{TalTwo} we obtain the analogue of Talagrand's interpolation for two copies of the system that is the main technical tool 
of the proof. 
In Section \ref{SecMin}, we summarize several properties of the parameters
in the definition of the Parisi formula in (\ref{Pk}) and (\ref{Pu}) and
in Sections \ref{SecFar} and \ref{SecClose} we prove the
apriori estimates of Section \ref{AprioriEst} by considering two
separate cases of ``far'' points and  ``close'' points.

\subsection{Talagrand's interpolation for two copies.}\label{TalTwo}

The key to proving the apriori estimate of Section \ref{AprioriEst}
is Talagrand's interpolation for two copies of the system.
The proof of this result is similar to the proof of Guerra's
interpolation in Theorem \ref{GTh3}. 
Once the Talagrand's interpolation is obtained,
the arguments in the rest of the paper will be adapted from 
\cite{T-P} with some necessary modifications. 

For $v\in [-D,D],$ let $v=\eta |v|$ for $\eta=\pm 1.$
Consider $\kappa\geq 1$ and consider a sequence $\vec{\nn}:$
\begin{equation}
0= \nn_0\leq \nn_1\leq \ldots \leq \nn_{\kappa}=1 
\label{Tens}
\end{equation}
such that $\nn_{\tau} = |v|$ for some $\tau\leq \kappa.$
Consider a sequence $\vrho:$
\begin{equation}
0=\rho_0\leq \ldots\leq \rho_{\kappa+1}=u.
\label{Trho}
\end{equation}
Consider a sequence of pairs of random variables $(Y_p^1,Y_p^2)$
independent for $0\leq p\leq \kappa$ such that
\begin{equation}
\e (Y_p^1)^2 = \e (Y_p^2)^2 = t(\xi'(\rho_{p+1}) - \xi'(\rho_{p}))
\label{yrand}
\end{equation}
and such that 
\begin{equation}
Y_p^1 = \eta Y_p^2 \mbox{ for } p<\tau
\,\,\mbox{ and } Y_p^1, Y_p^2 \mbox{ are independent for }
p\geq \tau.
\label{yrand2}
\end{equation}
Let $(Z_p^1,Z_p^2)$ be an arbitrary sequence of independent vectors 
for $0\leq p\leq \kappa.$ Let $(Z_{i,p}^1,Z_{i,p}^2)$ and 
$(Y_{i,p}^1,Y_{i,p}^2)$
be independent copies for $i\leq N$ of the sequences 
$(Z_{p}^1,Z_{p}^2)$ and $(Y_{p}^1,Y_{p}^2)$ and we assume
that they are independent of each other and the randomness in the Hamiltonian
$H_N(\vsi).$ For $s\in[0,1],$ let us define an interpolating Hamiltonian
\begin{equation}
H_s(\vsi^1,\vsi^2) = 
\sqrt{st}H_N(\vsi^1) +\sqrt{st} H_N(\vsi^2) + 
\sum_{j\leq 2}\sum_{i\leq N}\sigma_i^j \sum_{p\leq \kappa}
\Bigl(Z_{i,p}^j + \sqrt{1-s} Y_{i,p}^j\Bigr)
\label{THs}
\end{equation}
and consider 
\begin{equation}
F = \log \int_{A(v)} \exp H_s(\vsi^1,\vsi^2)d\nu(\vsi^1)d\nu(\vsi^2),
\label{TF}
\end{equation}
where $A(v)$ was defined in (\ref{UN}). 
Let
\begin{equation}
\chi(s) = \frac{1}{N}\e \P(\vec{\nn})F. 
\label{chis}
\end{equation}
\begin{theorem}\label{Tcopy}
For $s\in [0,1]$ we have
\begin{equation}
\chi'(s) \leq -2t\sum_{l<\tau}\nn_{l}(\theta(\rho_{l+1}) - \theta(\rho_l)) 
-t\sum_{l\geq \tau} \nn_l(\theta(\rho_{l+1}) - \theta(\rho_l)) + \R
\label{Tder}
\end{equation}
where $|\R| \leq K(\eps_N+c(N)).$
\end{theorem}
{\bf Proof.} 
The argument is similar to Guerra's interpolation
in Theorem \ref{GTh3}. (\ref{Pder}) implies that
$$
\chi'(s)=\frac{1}{N}\e W_1\ldots W_{\kappa} 
\frac{\partial F}{\partial s},
$$
where 
$$
W_l = \exp \nn_l(F_{l+1}-F_{l}) 
$$
and where $F_{l}$ are defined as in (\ref{Piter}).
Using the fact that $\nn_{\kappa}=1,$ one can write
\begin{eqnarray}
W_1\ldots W_{\kappa} 
&=& 
\exp \sum_{1\leq l\leq \kappa} \nn_l(F_{l+1} - F_{l})
\nonumber
\\
&=&
\exp \Bigl(F +
\sum_{1\leq l\leq \kappa} (\nn_{l-1} - \nn_{l}) F_{l} \Bigr)
= T \exp F
\label{TWT}
\end{eqnarray}
where 
$
T=T_1\ldots T_{\kappa} \mbox{ and }
T_l = \exp (\nn_{l-1} - \nn_{l}) F_{l}.
$
Using the definition of $F$ in (\ref{TF}),
$$
\chi'(s)=
\frac{1}{N}\e T \exp F
\frac{\partial F}{\partial s}
=\mbox{I}-\mbox{II}
$$
where
\begin{equation}
\mbox{I}=
\frac{\sqrt{t}}{2N\sqrt{s}} 
\int_{A(v)} \e T (H_N(\vsi^1)+H_N(\vsi^2)) \exp H_s(\vsi^1,\vsi^2) 
d\nu(\vsi^1)d\nu(\vsi^2).
\label{TI}
\end{equation}
and
\begin{equation}
\mbox{II}=
\frac{1}{2N\sqrt{1-s}} 
\int_{A(v)}
\e T \sum_{j\leq 2}\sum_{i\leq N}\sigma_i^j 
\sum_{p\leq \kappa} Y_{i,p}^j \exp H_s(\vsi^1,\vsi^2)
d\nu(\vsi^1)d\nu(\vsi^2).
\label{TII}
\end{equation}
To simplify the notations, 
let us denote $\vsi = (\vsi^1,\vsi^2)$ and $\vrho=(\vrho^1,\vrho^2),$
let $d\nu(\vsi) = d\nu(\vsi^1)d\nu(\vsi^2)$
and define
\begin{eqnarray}
\zz(\vsi,\vrho)
&=&
\frac{1}{N}\e (H_N(\vsi^1) + H_N(\vsi^2))
(H_N(\vrho^1) + H_N(\vrho^2))
\nonumber
\\
&=&
\zeta(\vsi^1,\vrho^1) + \zeta(\vsi^1,\vrho^2)
+\zeta(\vsi^2,\vrho^1) + \zeta(\vsi^2,\vrho^2)
\label{Tzeta}
\end{eqnarray}
where $\zeta$ was defined in (\ref{Gzeta}).
To compute I, we will use (\ref{GGI}) for the family 
$$
\vec{g}(\vrho) = H_N(\vrho^1)+H_N(\vrho^2) \mbox{ for }
\vrho = (\vrho^1,\vrho^2)\in A(v)
$$ 
and we will think of each factor $T_l$ in $T$ as the functional of
$\vec{g}(\vrho).$ 
Then (\ref{GGI}) and (\ref{Tzeta}) imply,
\begin{eqnarray} 
\mbox{I}
&=&
\frac{\sqrt{t}}{2\sqrt{s}}
\int_{A(v)} \e T
\frac{\partial \exp H_s(\vsi)}{\partial \vec{g}(\vsi)} 
\zz(\vsi,\vsi) d\nu(\vsi)
\label{TI1}
\\
&+&
\frac{\sqrt{t}}{2\sqrt{s}}\sum_{l\leq \kappa}
\int_{A(v)}
\e T\exp H_s(\vsi)\frac{1}{T_l}
\frac{\delta T_l}{\delta \vec{g}}[\zz(\vsi,\vrho)] d\nu(\vsi) .
\nonumber
\end{eqnarray}
First of all,
$$
\frac{\partial \exp H_s(\vsi)}{\partial \vec{g}(\vsi)}
=\sqrt{st} \exp H_s(\vsi),
$$
and, therefore, the first line in (\ref{TI1})
can be written as
\begin{equation}
\frac{t}{2}
\int_{A(v)} \e T \exp H_s(\vsi) \zz(\vsi,\vsi)d\nu(\vsi)
=\frac{t}{2}\e W_1\ldots W_{\kappa} \bigl\la
\zz(\vsi,\vsi)
\bigr\ra,
\label{TI2}
\end{equation}
where
$$
\bigl\la \zz(\vsi,\vsi) \bigr\ra = 
\exp(-F) \int_{A(v)} \zz(\vsi,\vsi) \exp H_s(\vsi) d\nu(\vsi).
$$
Using the definition of $T_l$ we can write
$$
\frac{1}{T_l}\frac{\delta T_l}{\delta \vec{g}}[\zz(\vsi,\vrho)]
=
(\nn_{l-1}-\nn_l) \frac{\delta F_{l}}{\delta \vec{g}}[\zz(\vsi,\vrho)]
=
(\nn_{l-1}-\nn_l)\e_l W_l\ldots W_{\kappa} 
\frac{\delta F}{\delta \vec{g}}[\zz(\vsi,\vrho)],
$$
where we used (\ref{Pder2}).
We have
$$
\frac{\delta F}{\delta \vec{g}}[\zz(\vsi,\vrho)]
=
\sqrt{st}\exp(-F)\int_{A(v)}\zz(\vsi,\vrho)
\exp H_s(\vrho) d\nu(\vrho)
=:
\sqrt{st}\bigl\la \zz(\vsi,\vrho)\bigr\ra',
$$
where $\la\cdot\ra'$ denotes the Gibbs average with respect to
$\vrho$ for a fixed $\vsi.$ 
Therefore, for a fixed $\vsi$ we get
$$
\frac{1}{T_l}\frac{\delta T_l}{\delta \vec{g}}[\zz(\vsi,\vrho)]
=
\sqrt{st}(\nn_{l-1}-\nn_l)\e_l W_l\ldots W_{\kappa} 
\bigl\la \zz(\vsi,\vrho)\bigr\ra'
=:
\sqrt{st}(\nn_{l-1}-\nn_l)\gamma_l(\zz(\vsi,\vrho)),
$$
where $\gamma_l$ here is defined by analogy with (\ref{Ggammal}). 
Combining this with the fact that
$$
T\exp H_s(\vsi)d\nu(\vsi)
=
W_1\ldots W_{\kappa} \exp(-F)\exp H_s(\vsi)d\nu(\vsi)
$$
we can write the second line in (\ref{TI1}) as
\begin{eqnarray}
&&
\frac{t}{2}\sum_{l\leq \kappa}(\nn_{l-1}-\nn_l)
\e W_1\ldots W_{\kappa} \exp(-F)
\int_{A(v)} \gamma_l(\zz(\vsi,\vrho))
d\nu(\vsi)
\label{TI3}
\\
&&
=
\frac{1}{2}\sum_{l\leq \kappa}(\nn_{l-1}-\nn_l)
\e W_1\ldots W_{l-1}\gamma_l^{\otimes 2}
(\zz(\vsi,\vrho))
=\frac{t}{2}\sum_{l\leq \kappa}(\nn_{l-1}-\nn_l)
\mu_l(\zz(\vsi,\vrho)),
\nonumber
\end{eqnarray}
where $\mu_l$ here is defined analogously to (\ref{Gmul}).
Combining (\ref{TI2}) and (\ref{TI3}) we get
\begin{equation}
\mbox{I} = 
\frac{t}{2}\e W_1\ldots W_{\kappa} \bigl\la
\zz(\vsi,\vsi)\bigr\ra 
+\frac{t}{2}\sum_{l\leq \kappa}(\nn_{l-1}-\nn_l)
\mu_l(\zz(\vsi,\vrho)).
\label{TI4}
\end{equation}
Let us denote by $R_{j,j'}$ the overlap of $\vsi^j$ and $\vsi^{j'}$
and by $R^{j,j'}$ the overlap of $\vsi^j$ and $\vrho^{j'}.$
From (\ref{correlation}) and (\ref{Tzeta}),
\begin{equation}
\mbox{I} = 
\frac{t}{2}\e W_1\ldots W_{\kappa} \bigl\la
\sum_{j,j'\leq 2} \xi(R_{j,j'})
\bigr\ra
+\frac{t}{2}\sum_{l\leq k}(\nn_{l-1}-\nn_l)
\mu_l(\sum_{j,j'\leq 2}\xi(R^{j,j'})) + {\cal R},
\end{equation}
where $|{\cal R}|\leq 4 c(N).$ Since the average $\la\cdot\ra$
is over the set $A(v),$ we have
\begin{equation}
|R_{j,j}-u|\leq \eps_N \mbox{ and }
|R_{1,2} - v|\leq \frac{1}{N}.
\label{Avee}
\end{equation}
Therefore,
\begin{equation}
\mbox{I} = 
t(\xi(u) + \xi(v))
+\frac{t}{2}\sum_{l\leq k}(\nn_{l-1}-\nn_l)
\mu_l(\sum_{j,j'\leq}\xi(R^{j,j'})) + {\cal R},
\label{TI5}
\end{equation}
where $|\R|\leq K(\eps_N + c(N)).$
The computation of $\mbox{II}$ is very similar.
For $p\leq \kappa,$ let us define
$$
\vec{g}_p(\vsi) = \sum_{j\leq 2}\sum_{i\leq N}\sigma_i^j Y_{i,p}^j
$$
so that $\mbox{II}=\sum_{p\leq \kappa}\mbox{II}(p)$ where
$$
\mbox{II}(p)=
\frac{1}{2N\sqrt{1-s}} 
\int_{A(v)} \e T \vec{g}_p(\vsi) \exp H_s(\vsi)d\nu(\vsi).
$$
Let us define
\begin{equation}
\zz_p(\vsi,\vrho) = \frac{1}{N}\e \vec{g}_p(\vsi) \vec{g}_p(\vrho)
=\sum_{j,j'\leq 2}R^{j,j'} \e Y_p^j Y_p^{j'}. 
\label{Tzetap}
\end{equation}
Then one can repeat the computations leading to (\ref{TI4})
with one important difference. One needs 
to note that $F_{l}$ does not depend on the r.v.
$Y_{i,p}$ for $l\leq p$ and, as a result, the summation
in the second term will be over $p<l\leq \kappa,$ i.e. 
\begin{equation}
\mbox{II}(p) = 
\frac{1}{2}\e W_1\ldots W_{\kappa} \bigl\la
\zz_p(\vsi,\vsi)\bigr\ra 
+\frac{1}{2}\sum_{p<l\leq \kappa}(\nn_{l-1}-\nn_l)
\mu_l(\zz_p(\vsi,\vrho))
\end{equation}
and
\begin{equation}
\mbox{II} = 
\frac{1}{2}\e W_1\ldots W_{\kappa} \bigl\la
\sum_{p\leq\kappa}\zz_p(\vsi,\vsi)\bigr\ra 
+\frac{1}{2}\sum_{l\leq \kappa}(\nn_{l-1}-\nn_l)
\mu_l\Bigl(\sum_{p<l}\zz_p(\vsi,\vrho)\Bigr)
\label{TI6}
\end{equation}
Using (\ref{Avee}) and (\ref{yrand}), it is easy to see
that the first term on the right hand side is 
$$
tu\xi'(u) + t|v|\xi'(|v|) + \R.
$$
Using (\ref{yrand}) and (\ref{yrand2}),
$$
\sum_{p<l} \e (Y_p^j)^2 = t\xi'(q_l)
\mbox{ and }
\sum_{p<l} \e Y_p^1 Y_p^2 = t\eta\xi'(q_{l\wedge \tau}) =
t\xi'(\eta q_{l\wedge \tau}),
$$
where in the last equality we used the fact that $\xi$ is even function.
If we define 
\begin{equation}
q_l^{j,j} = q_l \mbox{ and } q_l^{1,2} = \eta q_{l\wedge \tau}
\label{qjj}
\end{equation}
then $\sum_{p<l} \e Y_p^j Y_p^{j'} = t \xi'(q_l^{j,j'})$ and
(\ref{TI6}) can be written as
\begin{equation}
\mbox{II} = tu\xi'(u) + t|v|\xi'(|v|) + 
\frac{1}{2}\sum_{l\leq \kappa}(\nn_{l-1}-\nn_l)
\mu_l\Bigl(\sum_{j,j'\leq 2} R^{j,j'}\xi'(q_l^{j,j'})\Bigr)
+\R.
\label{TI7}
\end{equation}
Since $\theta(x) = x\xi'(x) - \xi(x)$ is also even, 
(\ref{TI5}) and (\ref{TI7})
imply that
\begin{eqnarray*}
\chi'(s) 
&=& 
\mbox{I} - \mbox{II} = -t\theta(u) - t\theta(v) 
+\frac{t}{2}\sum_{l\leq \kappa}(\nn_{l}-\nn_{l-1})\sum_{j,j'\leq 2}
\theta(q_l^{j,j'}) +\R
\\
&-&
\frac{1}{2}\sum_{l\leq \kappa}(\nn_{l}-\nn_{l-1})
\mu_l\Bigl(\sum_{j,j'\leq 2}
\Bigl(\xi(R^{j,j'}) - R^{j,j'}\xi'(q_l^{j,j'}) +\theta(q_{l}^{j,j'})
\Bigr)
\Bigr)
\\
&\leq&
-t\theta(u) - t\theta(v) 
+\frac{t}{2}\sum_{l\leq \kappa}(\nn_{l}-\nn_{l-1})\sum_{j,j'\leq 2}
\theta(q_l^{j,j'}) +\R
\\
&=&
-t\theta(u) - t\theta(v)
+\sum_{j,j'\leq 2}\theta(q_{\kappa+1}^{j,j'})
-\frac{t}{2}\sum_{j,j'\leq 2}\sum_{l\leq \kappa}\nn_{l}
(\theta(q_{l+1}^{j,j'}) - \theta(q_{l}^{j,j'}))
+\R
\\
&=&
-\frac{t}{2}\sum_{j,j'\leq 2}\sum_{l\leq \kappa}\nn_{l}
(\theta(q_{l+1}^{j,j'}) - \theta(q_{l}^{j,j'}))
+\R,
\end{eqnarray*}
since $q_{\kappa+1}^{j,j} = u$ and $q_{\kappa+1}^{1,2} = |v|.$
To finish the proof it remains to use (\ref{qjj}). We have
$$
\sum_{j,j'\leq 2}(\theta(q_{l+1}^{j,j'}) - \theta(q_{l}^{j,j'}))
=2(\theta(q_{l+1}) -\theta(q_l))
\mbox{ for } l<\tau
$$
and
$$
\sum_{j,j'\leq 2}(\theta(q_{l+1}^{j,j'}) - \theta(q_{l}^{j,j'}))
=\theta(q_{l+1}) -\theta(q_l)
\mbox{ for } l\geq\tau,
$$
since for $l\geq \tau,$
$
\theta(q_{l+1}^{1,2}) - \theta(q_{l}^{1,2})=
\theta(q_{\tau}) - \theta(q_{\tau}) = 0.
$
\qed

We can bound $\chi(0)$ as follows. 
Using (\ref{Avee}), for any $\lambda,\gamma\in\Reals,$
\begin{eqnarray*}
F\bigr|_{s=0} 
&\leq& 
-2N\lambda u - N\gamma v + 2N|\lambda| \eps_N +|\gamma|
\\
&+&
\log \int_{(\Sigma^N)^2} 
\exp 
\sum_{i\leq N}\Bigl(
\sum_{j\leq 2}\sigma_i^j \sum_{p\leq \kappa}
\bigl(Z_{i,p}^j + Y_{i,p}^j\bigr) 
+\sum_{j\leq 2}\lambda (\sigma_i^j)^2
+\gamma\sigma_i^1\sigma_i^2
\Bigr)
d\nu(\vsi^1)d\nu(\vsi^2)
\\
&=&
-2N\lambda u - N\gamma v + 2N|\lambda| \eps_N +|\gamma|
+\sum_{i\leq N}F_i(\lambda,\gamma),
\end{eqnarray*}
where $F_i(\lambda, \gamma)$ are independent copies of
\begin{equation}
F(\lambda,\gamma) = 
\log \int_{\Sigma\times\Sigma} 
\exp 
\Bigl(
\sum_{j\leq 2}\sigma_j \sum_{p\leq \kappa}
\bigl(Z_{p}^j + Y_{p}^j\bigr) 
+\sum_{j\leq 2}\lambda (\sigma_j)^2
+\gamma\sigma_1\sigma_2
\Bigr)
d\nu(\sigma_1)d\nu(\sigma_2).
\label{Flg}
\end{equation}
(\ref{Pprop1}), (\ref{Pprop2}) and (\ref{Pprop3}) now imply that
\begin{equation}
\chi(0)\leq -2\lambda u - \gamma v + 2|\lambda| \eps_N +|\gamma|N^{-1}
+\P(\vec{\nn})F(\lambda,\gamma)
\label{Tzero}
\end{equation}
and we obtained the following.

\begin{corollary}\label{corT}
If $\chi(s)$ and $F(\lambda,\gamma)$ are defined by
(\ref{chis}) and (\ref{Flg}) then for all $\lambda, \gamma\in\Reals,$
\begin{equation}
\chi(1) \leq -2\lambda u - \gamma v 
+ \P(\vec{\nn})F(\lambda,\gamma)
-2t\sum_{l<\tau}\nn_{l}(\theta(\rho_{l+1}) - \theta(\rho_l)) 
-t\sum_{l\geq \tau} \nn_l(\theta(\rho_{l+1}) - \theta(\rho_l)) + \R
\label{Tcor}
\end{equation}
where $|\R|\leq K(\eps_N + |\lambda|\eps_N +|\gamma|N^{-1} + c(N)).$
\end{corollary}
{\bf Proof.} The proof follows immediately by combining 
(\ref{Tder}) and (\ref{Tzero}).
\qed

{\bf Remark.} The remainder $\R$ in Theorem \ref{Apriori} will
be a result of application of (\ref{Tcor}). We will use it for
$\lambda = \lambda(u)$ defined in (\ref{lamt}) 
and, as in the proof of Lemma \ref{Lextra2},
$|\lambda(u)|\leq \Lambda$ for some $\Lambda$ that depends only
on $\xi,\nu$ and $u.$ Below we will use (\ref{Tcor}) 
for $|\gamma|\leq L$ for some constant $L$ independent of $N,t$ and $v.$ 
As a result, the remainder in Theorem \ref{Apriori}, $|\R|\leq a_N$
for $a_N = K(\eps_N+c(N)+N^{-1}).$

\subsection{Properties of $\eps$-minimizer.}\label{SecMin}

In this section we will describe several properties of the
sequence $(k,\vec{m},\vec{q},\lambda)$ in (\ref{mini})
that follow from its definition. Section 4 in \cite{T-P}
describes these properties in a very general setting,
with no reference to the classical SK model and
all the computations there apply to our case.
The only difference is that in \cite{T-P} certain generic computations
were applied to the function $\log\ch x$ and here
we will apply them to the function
\begin{equation}
x\to \log \int_{\Sigma} \exp(\sigma x +\lambda \sigma^2)d\nu(\sigma).
\label{lastconv}
\end{equation}
We will not reproduce some of the generic computations in \cite{T-P} 
that directly apply to our case.

\begin{figure}[t]
\centering
\psfrag{QR1}{$q_{r-1}$}
\psfrag{A}{$a$}
\psfrag{QR}{$q_r$}
\psfrag{QR11}{$q_{r+1}$}
\psfrag{MR1}{$m_{r-1}$}
\psfrag{M}{$m$}
\psfrag{MR11}{$m_r$}
\psfrag{Dot}{$\ldots$}
\psfrag{DOT}{$\ldots$}
\includegraphics[width=7cm,height=3.5cm]{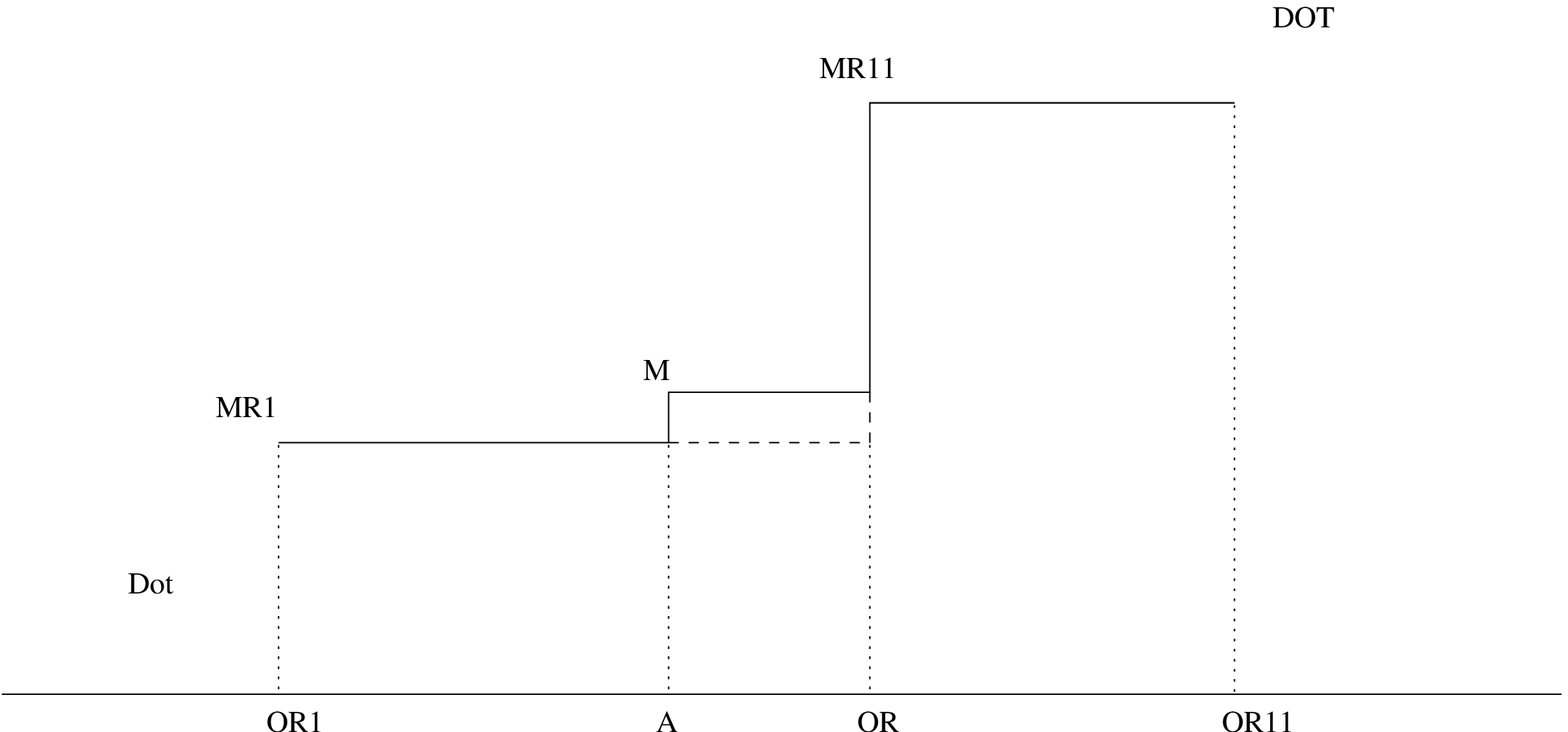}
\caption{\small Perturbing the $\eps$-minimizer.} 
\label{fig}
\end{figure}

{\bf Perturbing the sequences $\vec{m}$ and $\vec{q}$.} 
Let
$$
q_{r-1}\leq a\leq q_r, 
\,\,\, b=\xi'(a)
\mbox{ and } m_{r-1}\leq m\leq m_r
$$
and define new sequences $\vec{m}',\vec{q}'$ by inserting
$m$ and $a$ into $\vec{m}$ and $\vec{q}.$ Let us consider a
sequence $(z_p')_{p\leq k+1}$ of independent random variables
such that
$
\e (z_p')^2 = \xi'(q_{p+1}') - \xi'(q_{p}')
$
or, expressing this explicitly in terms of $b,$
\begin{eqnarray*}
&&
\e (z_p')^2 = \xi'(q_{p+1}) - \xi'(q_{p}) \mbox{ for } p<r-1,
\\
&&
\e (z_{r-1}')^2 = b - \xi'(q_{r-1}),\,\,
\e (z_r')^2 = \xi'(q_{r}) - b,
\\
&&
\e (z_{p+1}')^2 = \xi'(q_{p+1}) - \xi'(q_{p})
\mbox{ for } r\leq p\leq k.
\end{eqnarray*}
Let 
$$
F=\log\int_{\Sigma}\exp\Bigl(\sigma\sum_{p\leq k+1} z_p' 
+ \lambda\sigma^2\Bigr) d\nu(\sigma)
$$
and consider functions
\begin{equation}
T(m,b) = \P(\vec{m}')F
\label{Tmb}
\end{equation}
and
\begin{equation}
\Phi(m,a)= -\lambda u + T(m, \xi'(a))
-\frac{1}{2}\sum_{l\leq k} m_{l}(\theta(q_{l+1}) - \theta(q_l))
-\frac{1}{2} (m-m_{r-1})(\theta(q_r) - \theta(a)).
\label{Fma}
\end{equation}
Comparing with the definition (\ref{Pk}) it is clear that
\begin{equation}
T(m,\xi'(a))=X_0(\vec{m}',\vec{q}',\lambda) \mbox{ and }
\Phi(m,a) = \P_{k+1}(\vec{m}',\vec{q}',\lambda,u)
\label{TaX}
\end{equation}
and, thus, $T$ and $\Phi$ describe the behavior of $X_0,\P_k$
when we perturb the set of parameters by adding an extra point.
It will be very convenient to note that the functionals $X_0$
and $\P_k$ depend on the sequences $\vec{m}, \vec{q}$ only through
the function $m(q)$ defined by
\begin{equation}
m(q)=m_l \mbox{ for } q_l\leq q\leq q_{l+1},
\label{remark}
\end{equation}
which is called the {\it functional order parameter}.
Therefore, inserting parameters $m$ and $a$ can be visualized
as the perturbation of $m(q),$ as shown in Figure \ref{fig}.
From this point of view, it becomes obvious that for all $m$ and $a$ as above,
\begin{eqnarray}
&&
T(m_{r-1},b) = T(m,\xi'(q_r))
= X_0(\vec{m},\vec{q},\lambda),
\label{Tstable}
\\
&&
\Phi(m_{r-1},a) = \Phi(m,q_r)
= \P_{k}(\vec{m},\vec{q},\lambda,u)
\label{Fstable}
\end{eqnarray}
and, thus, it is very important to study the behavior of the 
derivatives of $T,\Phi$ in $m$ and $a$ at $m=m_{r-1}$ and $a=q_r.$

{\bf Properties of the derivatives.}
Let us define
\begin{equation}
U(b)=2\frac{\partial T}{\partial m}(m,b)\bigr|_{m=m_{r-1}}
\label{Ub}
\end{equation}
and
\begin{equation}
f(a) = \frac{\partial \Phi}{\partial m}(m,a)\bigr|_{m=m_{r-1}}
=\frac{1}{2} U(\xi'(a)) - \frac{1}{2}(\theta(q_r) - \theta(a)).
\label{Fa}
\end{equation}
The first fundamental formula is
\begin{equation}
\frac{\partial T}{\partial b}(m,b)\bigr|_{b=\xi'(q_r)}
=-\frac{1}{2}(m-m_{r-1})A 
\label{Tb}
\end{equation}
where $A$ is independent of $m.$ This can be obtained by a
straightforward computation and one can write down an explicit
formula for $A$ (see \cite{T-P}) but we will omit it here.
By definition of the $\eps$-minimizer (\ref{mini}), 
$(\vec{m},\vec{q},\lambda)$ is the minimizer of 
$\P_k(\vec{m},\vec{q},\lambda,u)$ and, therefore, 
\begin{equation}
\frac{\partial \P_k}{\partial q_r} = 0 \Longrightarrow
\frac{\partial X_0}{\partial q_r} +\frac{1}{2}(m_r -m_{r-1})
q_r\xi''(q_r) = 0,
\label{Xq}
\end{equation}
where we used (\ref{Pk}) and the fact that $\theta'(q_r) = q_r \xi''(q_r).$
When $m=m_r,$ it is apparent from Figure \ref{fig} that
the intervals $[a,q_r]$ and $[q_r,q_{r+1}]$ are glued together
in a sense that all functionals defined above become independent of
$q_r$ and, in particular, using (\ref{TaX}) 
\begin{equation}
T(m_{r},\xi'(a)) = X_0(\vec{m}',\vec{q}',\lambda)|_{m=m_r}
= X_0(\vec{m},\vec{q},\lambda)|_{q_r=a}.
\label{seen}
\end{equation}
Taking the derivative of both sides with respect to $a$ at $a=q_r$ gives 
$$
\frac{\partial X_0}{\partial q_r}(\vec{m},\vec{q},\lambda) = 
\frac{\partial T}{\partial a}(m_r,\xi'(a))\bigr|_{a=q_r}
=
\xi''(q_r)
\frac{\partial T}{\partial b}(m_r,b)\bigr|_{b=\xi'(q_r)}
=-\frac{1}{2}(m_r-m_{r-1})\xi''(q_r)A 
$$
where the last equality follows from (\ref{Tb}) for $m=m_r$.
Comparing this with (\ref{Xq}) implies the first
consequence of (\ref{mini}),
\begin{equation}
A=q_r.
\label{Aqr}
\end{equation}
Next, (\ref{Ub}) and (\ref{Tb}) imply that
$$
\frac{1}{2}U'(\xi'(q_r)) =
\frac{\partial}{\partial b}
\Bigl(
\frac{\partial T}{\partial m}(m,b)\Bigr|_{m=m_{r-1}}
\Bigr)\Bigr|_{b=\xi'(q_r)}
=
\frac{\partial}{\partial m}
\Bigl(
\frac{\partial T}{\partial b}(m,b)\Bigr|_{b=\xi'(q_r)}
\Bigr)\Bigr|_{m=m_{r-1}}
=
-\frac{1}{2}A,
$$
where, given any doubt, equality in the middle can be
checked by computing both sides. Therefore,
\begin{equation}
-U'(\xi'(q_r)) = A=q_r.
\label{UAqr}
\end{equation}
Another crucial property of $U(b)$ that can be verified
by  straightforward computation is
\begin{equation}
U''(b)\leq 0, 
\label{Ubconc}
\end{equation}
i.e. $U(b)$ is concave in $b.$
Next, let us describe several properties of $f(a)$ in (\ref{Fa}).
We have
\begin{equation}
f(q_r) = f'(q_r) = 0, \,\, f(q_{r-1}) \geq 0.
\label{fprime}
\end{equation}
The first one follows from
$$
f(q_r) = \frac{1}{2}U(\xi'(q_r)) = 
\frac{\partial T}{\partial m}(m,\xi'(q_r))\bigr|_{m=m_{r-1}}
= 0,
$$
since (\ref{Tstable}) yields that $T(m,\xi'(q_r))$ does not
depend on $m$. The second one follows from (\ref{Fa}) and
(\ref{UAqr}) since
$$
f'(q_r) = \frac{1}{2}\xi''(q_r)(U'(\xi'(q_r)) + q_r) = 0.
$$
To show that $f(q_{r-1}) \geq 0$ let us note that
$$
\Phi(m,q_{r-1}) = \P_k(\vec{m},\vec{q},\lambda,u)|_{m_{r-1}=m}
$$
since setting $a=q_{r-1}$ simply replaces $m_{r-1}$ with $m$
in the definition of $\P_k,$ which is also apparent from Figure \ref{fig}. 
If $m_{r-1}>0,$ then as in  (\ref{Xq}),
$$
\frac{\partial \P_k}{\partial m_{r-1}} = 0 \Longrightarrow
f(q_{r-1})=\frac{\partial \Phi}{\partial m}(m,q_{r-1})|_{m=m_{r-1}} =0.
$$
If $m_{r-1}=0$ then since $(\vec{m},\vec{q},\lambda)$ is
the minimizer of $\P_k,$ slightly increasing $m_{r-1}$ should not
decrease $\P_k$ and, therefore, the right derivative
$\partial \P_k/\partial m_{r-1} \geq 0$ and $f(q_{r-1})\geq 0.$

{\bf $\eps$-dependent properties of the derivatives.}
So far we have only utilized the fact that 
$(\vec{m},\vec{q},\lambda)$ is the minimizer of $\P_k$
and we have not used the condition in (\ref{mini}) that
$$
\P_k(\vec{m},\vec{q},\lambda,u)\leq \P(\xi,u) +\eps.
$$
In particular, this implies that for any $m$ and $a,$
\begin{equation}
\Phi(m_{r-1},a) = \P_{k}(\vec{m},\vec{q},\lambda,u)
\leq \P(\xi,u) + \eps \leq \Phi(m,a) +\eps,
\label{useeps}
\end{equation}
which means that we can not decrease $\P_k$ much by
varying parameters $m,a.$
This can be combined with the following fact that
plays a central role:
\begin{eqnarray}
&&
\mbox{\hspace{-1.5cm}All derivatives of $T,\Phi,U, f$ with
respect to $a,b,m$ are bounded by constants} 
\nonumber
\\
&&
\mbox{\hspace{-1.5cm} that depend
only on $\xi,\nu, u$ but not on $(k,\vec{m},\vec{q})$.}
\label{fact}
\end{eqnarray}
Let $L$ denote such constants that depend only on $\xi,\nu, u.$
The proof of (\ref{fact}) in \cite{T-P} relied on the fact
that $\e \ch(c+c' z)^L \leq L$ for a standard Gaussian $z$ 
and $|c| ,|c'|\leq L.$ 
In our case, this conditions will be replaced by an obvious
condition,
$$
\e \Bigl(\int_{\Sigma} \exp(\sigma c z + \lambda \sigma^2)d\nu(\sigma) 
\Bigr)^L\leq L \mbox{ for } |c|,|\lambda| \leq L.
$$
The following Lemma holds.
\begin{lemma}\label{epscond}
The function $f(a)$ in (\ref{Fa}) satisfies,
\begin{equation}
f(a)\geq -L\sqrt{\eps}
\label{f1}
\end{equation}
and
\begin{equation}
f''(q_r) = \frac{1}{2}\xi''(q_r)\bigl(
\xi''(q_r)U''(\xi'(q_r)) + 1
\bigr) \geq - L\eps^{1/6}.
\label{f2}
\end{equation}
\end{lemma}
{\bf Proof.}
(\ref{f1}) holds if $f(a)\geq 0$ so we can that $f(a)<0.$
Using (\ref{Fa}) and (\ref{fact}), we can write
\begin{equation}
\Phi(m,a) \leq \Phi(m_{r-1},a)+(m-m_{r-1})f(a) + L(m-m_{r-1})^2.
\label{compa}
\end{equation}
By (\ref{Fstable}) and the fact that $(\vec{m},\vec{q},\lambda)$
is a minimizer of $\P_k$,
$$
\Phi(m_{r-1},a) = \P_{k}(\vec{m},\vec{q},\lambda,u)
\leq \P_{k}(\vec{m},\vec{q},\lambda,u)|_{q_{r}=a}
= \Phi(m_r,a),
$$
and the last equality can been seen as in (\ref{seen}).
Therefore, (\ref{compa}) with $m=m_r$ implies that
$$
(m_r-m_{r-1})f(a) + L(m_r-m_{r-1})^2 \geq 0
$$
and, therefore,
\begin{equation}
m_{r}\geq m_{r-1} -\frac{f(a)}{L}\geq m_{r-1},
\label{mchoice}
\end{equation}
where the second inequality follows from our assumption that
$f(a)<0.$ (\ref{useeps}) and (\ref{compa}) imply that
for any $m$ and $a,$
$$
-\eps \leq (m-m_{r-1})f(a) + L(m-m_{r-1})^2. 
$$
Taking $m:=m_{r-1} - f(a)/2L$, which belongs to
$[m_{r-1},m_r]$ by (\ref{mchoice}), implies
$$
-\eps\leq -f(a)^2/4L
$$
and this proves (\ref{f1}).
It remains to prove (\ref{f2}). From (\ref{Fa}) and (\ref{UAqr})
we see that
\begin{eqnarray*}
f''(q_r) 
&=& 
\frac{1}{2}\xi'''(q_r)\bigl(U'(\xi'(q_r)) +q_r\bigr)
+\frac{1}{2}\xi''(q_r)\bigl(U''(\xi'(q_r))\xi''(q_r)+1\bigr)
\\
&=&
\frac{1}{2}\xi''(q_r)\bigl(U''(\xi'(q_r))\xi''(q_r)+1\bigr).
\end{eqnarray*}
(\ref{fact}) and (\ref{f1})
imply that for any $a,$
\begin{equation}
-L\sqrt{\eps}\leq f(a)\leq \frac{1}{2}f''(q_r)(a-q_r)^2 + L|a-q_r|^3. 
\label{epsq}
\end{equation}
Using (\ref{fprime}) we have
$$
0\leq f(q_{r-1})\leq \frac{1}{2}f''(q_r)(q_{r-1}-q_r)^2 + L|q_{r-1}-q_r|^3
$$
and, hence,
$$
q_r \geq q_{r-1}-\frac{f''(q_r)}{2L}.
$$
If $f''(q_r)\geq 0$ then (\ref{f2}) holds, otherwise 
$$
a= q_{r-1}-\frac{f''(q_r)}{2L} \in [q_{r-1},q_r]
$$
and using (\ref{epsq}) for this choice of $a$ again
implies (\ref{f2}).
\qed

{\bf Dual construction and the replica symmetric case.}

The construction above will be used in the proof of Theorem \ref{Apriori} 
to provide control of the points on the left hand side of $q_r,$ 
i.e. $q_{r-1}\leq v\leq q_r.$ In order to provide control of the points
on the right hand side $q_r\leq v\leq q_{r+1}$ one can consider 
a dual construction
by perturbing parameter $m_r$ on the interval $[q_r,a].$ This construction
is very similar so we will not detail it and we will only consider the
points on the left hand side in Theorem \ref{Apriori}.
In the replica symmetric case, the function $\Phi(m,a)$ defined in
(\ref{RSBP}) is the analogue of the function in (\ref{Fma}) and 
the properties  (\ref{RSRScond}) and (\ref{RSAT}) replace the
properties (\ref{f1}) and (\ref{f2}) and, in fact, are stronger because
$\eps$ is replaced by $0.$ (The change of sign in the inequalities
is simply because we consider a dual construction.) 
Therefore, the proof of the replica symmetric apriori estimate in
Theorem \ref{RSApriori} is exactly the same as the proof of Theorem 
\ref{Apriori} if we use (\ref{RSRScond}) and (\ref{RSAT}) instead of
(\ref{f1}) and (\ref{f2}), and we will not detail it.

\subsection{Control of the far points.}\label{SecFar}

In this section we will prove the easiest case of Theorem \ref{Apriori},
when the point $v$ is far from $q_r$ in the following sense. Given
$\eps>0,$ let $(k,\vec{m},\vec{q},\lambda)$ be an $\eps$-minimizer
defined by (\ref{mini}). Without loss of generality, we will assume
that all coordinates of the vector $\vec{m}$ are different and that all
coordinates of the vector $\vec{q}$ are also different. Otherwise,
we can decrease the value of $k$ by gluing equal coordinates
without changing the value of the functional $\P_k(\vec{m},\vec{q},\lambda).$
In this section we will consider the case when $v\in[-D,D]$ is such that
\begin{eqnarray}
v< q_{r-1} 
&\mbox{or}& 
v> q_{r+1}
\,\,\,\mbox{ if }\,\,\, 
r\geq 2
\nonumber
\\
v< -q_{1} 
&\mbox{or}& 
v> q_{2}
\,\,\,\mbox{ if }\,\,\, 
r=1
\label{far}
\end{eqnarray}
and we will prove the following.
\begin{proposition}\label{Propfar}
In the notations of Theorem \ref{Apriori}, if (\ref{far}) holds
then
\begin{equation}
\frac{1}{N}\e\P(\vec{n})F(A(v))\leq 2\psi(t)-K+\R,
\label{casefar}
\end{equation}
where $K>0$ is a constant independent of $N,t$ and $v.$
\end{proposition}
This implies Theorem \ref{Apriori} in the range of parameters (\ref{far})
because
$$
K\geq (v-q_r)^2\big/ 
\bigl(K^{-1}(q_r-q_{r-1})^2\wedge (q_{r+1}-q_r)^2\bigr).
$$
The proof of Proposition \ref{Propfar} is based on
the following three step construction.

{\bf 1.} Let us recall the definition of $\vec{n}$ in (\ref{Rn})
and $(z_p^1,z_p^2)_{p\leq k}$ in (\ref{Rzz}).

{\bf 2.} ({\it Inserting $|v|$.}). Let $a$ be such that 
\begin{equation}
q_{a}\leq |v| \leq q_{a+1}.
\label{Fara}
\end{equation}
Consider a vector
$$
\vec{q}' =(q_0',\ldots,q_{k+2}')
= (q_0,\ldots,q_{a-1},|v|,q_{a},\ldots,q_{k+1})
$$
which is defined by inserting $|v|$ in the vector $\vec{q}$
and define a vector 
$$
\vec{n}'=(\frac{m_0}{2},\ldots,\frac{m_{a-1}}{2},m_{a-1},m_{a}
,\ldots,m_{k}).
$$ 
Consider a sequence $(y_{p})_{p\leq k+1}$ of independent
Gaussian random variables such that 
$$
\e y_p^2 = \xi'(q_{p+1}') - \xi'(q_{p}')
$$
and let $(y_p^1,y_p^2)_{p\leq k+1}$ consist of two
copies of $(y_p)_{p\leq k+1}$ such that
\begin{equation}
y_{p}^1 = \eta y_{p}^2 \mbox{ for }
p<a \mbox{ and } y_{p}^1,  y_{p}^2 \mbox{ are independent for }
p\geq a,
\end{equation}
where $v=\eta |v|.$

{\bf 3.} ({\it Gluing two sequences together}).
Let $\kappa = 2k+1.$
Let us consider a vector  
$$
\vec{\nn} = (\nn_0,\ldots,\nn_{\kappa +1})
\mbox{ such that }
\nn_0\leq \ldots\leq\nn_{\kappa+1}
$$
and such that $\vec{\nn}$ consists of the elements of
vectors $\vec{n}$ and $\vec{n}'.$ More precisely, there
exists a partition $I,J$ of the set $\{0,\ldots, \kappa+1\}$ 
with $\card I = k+1$ and $\card J =k+2$ such that
the elements $\nn_p$ are the elements of $\vec{n}$ 
for $p\in I $ and the elements of $\vec{n}'$ for $p\in J.$
For $0\leq p\leq \kappa,$ define
\begin{equation}
(Z_p^1,Z_p^2)=0 \mbox{ for } p\in J \mbox{ and }
(Z_p^1,Z_p^2) = \sqrt{1-t}(z_l^1,z_l^2) \mbox{ for }
p\in I
\label{FarZ}
\end{equation}
and where $l$ is such that $\nn_p = n_l.$ Similarly,
define
\begin{equation}
(Y_p^1,Y_p^2)=0 \mbox{ for } p\in I \mbox{ and }
(Y_p^1,Y_p^2) = \sqrt{t}(y_l^1,y_l^2) \mbox{ for }
p\in J
\label{FarY}
\end{equation}
and where $l$ is such that $\nn_p = n_l'.$
For $0\leq p\leq\kappa,$ let us define
\begin{equation}
(g_p^1,g_p^2)=(Z_p^1,Z_p^2)+(Y_p^1,Y_p^2)
\label{gees}
\end{equation}
or, in other words,
$$
(g_p^1,g_p^2)=(Z_p^1,Z_p^2) \mbox{ for } p\in I \mbox{ and }
(g_p^1,g_p^2) = (Y_p^1,Y_p^2) \mbox{ for } p\in J.
$$
In order to match this definition of $(Y_p^j)$
with (\ref{Trho}) and (\ref{yrand}), let us define
a sequence
$$
\rho_0\leq \ldots\leq \rho_{\kappa+1}
$$  
such that $\rho_p = q_l'$ for $l$ such that $\nn_p=n_l'.$
Define $\tau$ by $\nn_{\tau} = n_{a}' = m_{a-1}.$

We will now apply Corollary \ref{corT} 
to these choices of $\vec{\nn},$ $(Z_p^j)$ and $(Y_p^j).$
First of all, 
\begin{equation}
\chi(1) = \frac{1}{N}\e\P(\vec{n})F(A(v))
\label{chione}
\end{equation}
since for $s=1,$ the random variables $(Y_p^j)$ will disappear
in the definition of $H_s(\vsi^1,\vsi^2)$ in (\ref{THs}),
the Hamiltonian $H_s(\vsi^1,\vsi^2)$ will coincide with
Hamiltonian $H_t(\vsi^1,\vsi^2)$ in (\ref{RHt}) and,
as a result, the definition of $F$ in (\ref{TF}) will
coincide with $F(A(v))$ in Theorem \ref{Apriori}.
Next, it is clear from the construction that
\begin{eqnarray*}
&&
-2t\sum_{l<\tau}\nn_{l}(\theta(\rho_{l+1}) - \theta(\rho_l)) 
-t\sum_{l\geq \tau} \nn_l(\theta(\rho_{l+1}) - \theta(\rho_l)) 
\\
&=&
-2t\sum_{l<a} n_{l}'(\theta(q_{l+1}') - \theta(q_l')) 
-t\sum_{l\geq a} n_l'(\theta(q_{l+1}') - \theta(q_l')) 
\\
&=&
-t\sum_{l\leq k} m_{l}(\theta(q_{l+1}) - \theta(q_l)). 
\end{eqnarray*}
Corollary \ref{corT} now implies that
$$
\frac{1}{N}\e\P(\vec{n})F(A(v)) \leq
-2\lambda u - \gamma v + \P(\vec{\nn})F(\lambda,\gamma)
-t\sum_{l\leq k} m_{l}(\theta(q_{l+1}) - \theta(q_l))
+\R.
$$
We will use this bound for $\lambda$ as in the $\eps$-minimizer 
$(k,\vec{m},\vec{q},\lambda)$ and $\gamma=0,$ i.e.
$$
\frac{1}{N}\e\P(\vec{n})F(A(v)) \leq
-2\lambda u + \P(\vec{\nn})F(\lambda,0)
-t\sum_{l\leq k} m_{l}(\theta(q_{l+1}) - \theta(q_l))
+\R.
$$ 
The argument in Lemma \ref{Lextra2} shows that 
$|\lambda|\leq \Lambda$ for a constant $\Lambda$ 
that depends only on $\xi, \nu$ and $u$ 
and, hence, $|\R|\leq K(\eps_N + c(N)).$ 
Recalling the definition of $\psi(t)$ in (\ref{RRpsi}), 
in order to prove Proposition \ref{Propfar},
it remains to show that
$$
\P(\vec{\nn})F(\lambda,0) \leq 2X_0(\vec{m},\vec{q},\lambda) -K
=2\P(\vec{m})X_{k+1} -K,
$$
where we used (\ref{PXlast}) and
where $K>0$ is a constant independent of $t$ and $v.$
In fact, it is enough to show that
\begin{equation}
\P(\vec{\nn})F(\lambda,0) < 2X_0(\vec{m},\vec{q},\lambda)
=2\P(\vec{m})X_{k+1},
\label{Farstrict}
\end{equation}
for all parameters $t\in[0,1-t_0]$ and $v$ as in (\ref{far}) because the
functionals on both sides are continuous in these parameters
and, even though the set defined in (\ref{far}) is not a compact,
the case of the end points will be proved in the following sections
and (\ref{Farstrict}) holds on the closure of (\ref{far}).
Therefore, by continuity and compactness, strict inequality
for each $(t,v)$ will imply strict inequality
uniformly over the entire set of parameters. 
The proof of (\ref{Farstrict}) repeats the proof of Proposition 5.7
in \cite{T-P} with only one modification that instead of $\log \ch x$
we consider (\ref{lastconv})
and note that this function is also strictly convex in $x$
because we eliminated the case when $\nu$ is concentrated
on one point in Section \ref{SecRed}.
Instead of reproducing
the proof in its entirety we will explain a very clear idea behind
it by looking at a few cases.
Let us first show that a nonstrict version of (\ref{Farstrict}), i.e.
\begin{equation}
\P(\vec{\nn})F(\lambda,0) \leq 2X_0(\vec{m},\vec{q},\lambda)
\label{nonstrict}
\end{equation}
always holds, even without the assumption (\ref{far}).
(\ref{Flg}) gives that
$F(\lambda,0)=F^1 + F^2$ where
$$
F^j = 
\log \int_{\Sigma} 
\exp 
\Bigl(
\sigma \sum_{p\leq \kappa} g_{p}^j 
+\lambda \sigma^2
\Bigr)d\nu(\sigma).
$$
It is clear from the construction that for $p\in I$
\begin{equation}
\nn_p = n_l=\frac{m_l}{2} \mbox{ for } l<r \Longleftrightarrow
g_p^1 = g_p^2
\label{FarI}
\end{equation}
and for $p\in J$
\begin{equation}
\nn_p = n_l'=\frac{m_l}{2} \mbox{ for } l<a \Longleftrightarrow
g_p^1 = \eta g_p^2.
\label{FarJ}
\end{equation}
In other words, $\nn_p$ is of the type $m_l/2$ whenever the
corresponding random pair is fully correlated.
Let us define a vector $\vec{\m}=(\m_0,\ldots,\m_{\kappa})$ by
\begin{equation}
\m_p = 2\nn_p \mbox{ for $p$ as in (\ref{FarI}) or (\ref{FarJ}) and } 
\m_p = \nn_p \mbox{ otherwise}.
\label{Farm}
\end{equation}
A fact that plays a very important role below is that
coordinates of $\vec{\m}$ are not necessarily arranged
in an increasing order. Let us first prove the following.
\begin{lemma}\label{FarNT}
We have
\begin{equation}
\P(\vec{\nn})F(\lambda,0) = 
\P(\vec{\nn})(F^1+F^2)\leq \P(\vec{\m}) F^1 + \P(\vec{\m}) F^2.
\label{Farntom}
\end{equation}
\end{lemma}
{\bf Proof.}
This follows by induction in (\ref{Piter}).
For $p$ such that $\m_p = \nn_p$ and $g_p^1,g_p^2$ are
independent we have
\begin{eqnarray}
(F^1+F^2)_{p} 
&=& 
\frac{1}{\nn_p} \log \e_p \exp \nn_p (F^1+F^2)_{p+1}
\leq \frac{1}{\nn_p} \log \e_p \exp \nn_p (F_{p+1}^1+F_{p+1}^2)
\nonumber
\\
&=& 
\frac{1}{\m_p} \log \e_p \exp \m_p F_{p+1}^1 + 
\frac{1}{\m_p} \log \e_p \exp \m_p F_{p+1}^2
= F_p^1 + F_p^2.
\end{eqnarray}
For $p$ such that $\m_p = 2\nn_p$ and $g_p^1=\pm g_p^2$ 
we have
\begin{eqnarray}
(F^1+F^2)_{p} 
&=& 
\frac{1}{\nn_p} \log \e_p \exp \nn_p (F^1+F^2)_{p+1}
\leq \frac{2}{\m_p} \log \e_p \exp \frac{\m_p}{2} (F_{p+1}^1+F_{p+1}^2)
\nonumber
\\
&\leq& 
\frac{1}{\m_p} \log \e_p \exp \m_p F_{p+1}^1 + 
\frac{1}{\m_p} \log \e_p \exp \m_p F_{p+1}^2
= F_p^1 + F_p^2,
\label{FarHolder}
\end{eqnarray}
where in the second line we used H\"older's inequality.
For $p=0$ this gives (\ref{Farntom}). 
\qed

\begin{lemma}\label{FarPcomp}
Let $\vec{\m}'$ be a nondecreasing permutation of the vector $\vec{\m}.$
Then,
\begin{equation}
\P(\vec{\m})F^j \leq \P(\vec{\m}') F^j = \P(\vec{m}) X_{k+1} = 
X_0(\vec{m},\vec{q},\lambda).
\label{Farperm}
\end{equation}
\end{lemma}

The first inequality means that the Parisi functional will
decrease if $\vec{m}$ is not arranged in an increasing order.
Lemma \ref{FarPcomp} 
follows from Lemma 5.12 in \cite{T-P} which states the following.
Given a function $Q$ and numbers $a\geq 0$ and $m>0,$ let
$$
T_{m,a}(Q)(x)=\frac{1}{m}\log \e \exp m Q(x+g \sqrt{a})
$$
where $g$ is standard Gaussian. 
\begin{lemma}\label{FarTlem}(\cite{T-P})
If $a,a'\geq 0$ and $m\geq m'$ then for each $x$ we have
\begin{equation}
T_{m,a}\circ T_{m',a'}(Q)(x) \leq 
T_{m',a'}\circ T_{m,a}(Q)(x).
\label{FarT} 
\end{equation}
If $a,a'>0$ and $m>m'$ then we can have equality only if $Q$ is constant.
\end{lemma}

{\bf Proof of Lemma \ref{FarPcomp}.}
The first inequality is obvious by Lemma \ref{FarTlem}.
Equality in (\ref{Farperm}) follows by construction.
The elements of $\vec{\m}'$ are precisely the elements
of $\vec{m}.$ The random variables $g_p^j$ for $p$
such that $\m_p' = m_l$ are exactly
$$
\sqrt{1-t} z_l^j \mbox{ and } \sqrt{t} y_l^j
\mbox{ for } l\not= a-1
$$
and
$$
\sqrt{1-t} z_l^j \mbox{ and } \sqrt{t} y_l^j,\,
\sqrt{t} y_{l+1}^j \mbox{ for } l = a-1.
$$
Obviously, 
\begin{equation}
T_{m,a}\circ T_{m,a'} = T_{m,a+a'}, 
\label{Ttwice}
\end{equation}
which means that
we can combine the random variables corresponding to the same value
$m_l$ and since
it is easy to check that in both cases the sum of these
random variables is equal in distribution to $z_l$
defined in (\ref{z}), (\ref{Farperm}) follows.
\qed

Combining Lemma \ref{FarNT} and \ref{FarPcomp}, we proved
(\ref{nonstrict}). From the proof it is clear
that there are only two places, (\ref{FarHolder}) and (\ref{FarT}),
where the inequality could become strict. 
It turns out that condition (\ref{far}) ensures that gluing two
sequences together occurs in such a way that at least
in one of these two steps the inequality will become strict.
We will not present the detailed proof here and refer
a reader to Proposition 5.7 in \cite{T-P}. We will explain the main
idea by looking at several typical cases.

Let us consider the case $r\geq 2$ in (\ref{far}). The case
$r=1$ is quite similar with the exception that the interval
$-q_1\leq v\leq q_0=0$  was excluded in (\ref{far}) because it requires
a different approach and it will be postponed until the following sections.

For $r\geq 2$ and $v$ as in (\ref{far}) we have
$$
\mbox{ (a) } |v|\not\in [q_{r-1},q_{r+1}]
\,\,\mbox{ or\,\, (b) } |v|\in [q_{r-1},q_{r+1}] \mbox{ and } 
v=-|v|, \mbox{ i.e. } \eta=-1.
$$ 
{\bf Case (a).} Let us assume for simplicity 
that $q_{r-2}\leq |v|<q_{r-1}$ since other cases are similar.
This corresponds to the case $a=r-2$ in (\ref{Fara}).
Then we will split case (a) into two subcases:
\begin{equation}
q_{r-2}\leq |v|<q_{r-1} \mbox{ and } m_{r-2}<m_{r-1}/2;
\label{a'}
\end{equation}
\begin{equation}
q_{r-2}\leq |v|<q_{r-1} \mbox{ and } m_{r-1}/2 \leq m_{r-2}.
\label{a''}
\end{equation}
Let us now see what happens when we combine the sequences at step 3
above. First of all, at step 1 the sequences $\vec{n}$ and
$(z_p^1,z_p^2)$ will have subsequences
\begin{equation}
\begin{array}{cccccc}
\ldots & m_{r-2}/2 & m_{r-1}/2 & m_r & m_{r+1} & \ldots
\\
\ldots & (z_{r-2},z_{r-2}) & (z_{r-1},z_{r-1}) & 
(z_{r}^1,z_r^2) & (z_{r+1}^1,z_{r+1}^2) & \ldots
%\\
%\ldots & z_{r-2} & z_{r-1} & z_{r}^2 & z_{r+1}^2 & \ldots
\end{array}
\label{step1seq}
\end{equation}
At step 2 the sequences $\vec{n}'$ and
$(y_p^1,y_p^2)$ will have subsequences
\begin{equation}
\begin{array}{cccccc}
\ldots & m_{r-2}/2 & m_{r-2} & m_{r-1} & m_{r} & \ldots
\\
\ldots & (y_{r-2},\eta y_{r-2}) & (y_{r-1}^1,y_{r-1}^2) & 
(y_{r}^1,y_r^2) & (y_{r+1}^1,y_{r+1}^2) & \ldots
%\\
%\ldots & \eta y_{r-2} & y_{r-1}^2 & y_{r}^2 & y_{r+1}^2 & \ldots
\end{array}
\label{step2seq}
\end{equation}
In both cases we write $(z^1,z^2)$ or $(y^1,y^2)$ whenever 
two coordinates are 
independent. When we glue these sequences together as step 3,
the sequences $\vec{\nn}$ and $(g_p^1,g_p^2)$ will contain
subsequences
\begin{equation}
\begin{array}{ccccccc}
\ldots & m_{r-2}/2 & m_{r-2}/2 & m_{r-2} & m_{r-1}/2 & 
m_{r-1} & \ldots
\\
\ldots & (Z_{r-2},Z_{r-2}) & (Y_{r-2},\eta Y_{r-2}) & 
(Y_{r-1}^1,Y_{r-1}^2) & (Z_{r-1},Z_{r-1}) &
(Y_{r}^1,Y_r^2) & \ldots
\end{array}
\label{seq1}
\end{equation}
in the case (\ref{a'}) and 
\begin{equation}
\begin{array}{ccccccc}
\ldots & m_{r-2}/2 & m_{r-2}/2 & m_{r-1}/2 & m_{r-2} &  
m_{r-1} & \ldots
\\
\ldots & (Z_{r-2},Z_{r-2}) & (Y_{r-2},\eta Y_{r-2}) & 
(Z_{r-1},Z_{r-1}) & (Y_{r-1}^1,Y_{r-1}^2) & 
(Y_{r}^1,Y_r^2) & \ldots
\end{array}
\label{seq2}
\end{equation}
in the case (\ref{a''}). 

Suppose that (\ref{seq1}) is the case. Then
the strict inequality will appear when we apply
equation (\ref{FarHolder}) at the step when 
$\nn_p$ is equal to $m_{r-1}/2.$ Indeed, at this step
$$
F_{p+1}^1 = F_{p+1}^1 (\ldots+ Z_{r-2} + Y_{r-2} + Y_{r-1}^1 +
Z_{r-1}),
$$
$$
F_{p+1}^2 = F_{p+1}^2 (\ldots+ Z_{r-2} + Y_{r-2} + Y_{r-1}^2 +
Z_{r-1})
$$
and $Y_{r-1}^1, Y_{r-1}^2$ are independent and nondegenerate
since $\e (Y_{r-1}^j)^2 = t(\xi'(q_{r-1}) - \xi'(|v|))>0$
by (\ref{a'}). Also, both functions $x\to F_{p+1}^j(x)$
are strictly convex because (\ref{lastconv}) is strictly convex
and iteration in the Parisi functional (\ref{Piter})
will preserve strict convexity. Therefore, $F_{p+1}^1$ and
$F_{p+1}^2$ are not collinear as functions of $Z_{r-1}$
with probability one over $(Y_{r-1}^1,Y_{r-1}^2)$ and,
therefore, H\"older's inequality in (\ref{FarHolder})
will be strict with probability one.

Now, suppose that (\ref{seq2}) holds. Then after using Lemma \ref{FarNT},
$\P(\vec{\m}) F^1$ will be defined in terms of the sequences that contain
subsequences
\begin{equation}
\begin{array}{ccccccc}
\ldots & m_{r-2} & m_{r-2} & m_{r-1} & m_{r-2} &  
m_{r-1} & \ldots
\\
\ldots & Z_{r-2} & Y_{r-2} & 
Z_{r-1} & Y_{r-1}^1  & 
Y_{r}^1 & \ldots
\end{array}
\end{equation}
In this case, 
$\m$ is not arranged in an increasing order, since $m_{r-1}>m_{r-2},$
and $Z_{r-1}, Y_{r-1}^1$ are nondegenerate. 
Therefore, when we rearrange these sequences in an increasing order
by applying Lemma \ref{FarTlem}, we will get strict inequality
in (\ref{FarT}).

{\bf Case (b).} In this case the scenario of (\ref{seq2}) can not occur
and the fact that $\eta=-1$ plays an important role.
Suppose for certainty that $q_{r-1}\leq |v|\leq q_{r}.$
(\ref{step1seq}) does not change but instead of (\ref{step2seq}) 
we will now have:
\begin{equation}
\begin{array}{cccccc}
\ldots & m_{r-2}/2  & m_{r-1}/2 & m_{r-1} & m_r & \ldots
\\
\ldots & (y_{r-2}, -y_{r-2}) & (y_{r-1},-y_{r-1}) & 
(y_{r}^1,y_r^2) & (y_{r+1}^1,y_{r+1}^2) & \ldots
%\\
%\ldots & \eta y_{r-2} & y_{r-1}^2 & y_{r}^2 & y_{r+1}^2 & \ldots
\end{array}
\label{step2seqb}
\end{equation}
When we glue this sequence with (\ref{step1seq}) we will get
\begin{equation}
\begin{array}{ccccccc}
\ldots & m_{r-2}/2 & m_{r-2}/2 & m_{r-1}/2 & m_{r-1}/2 & 
m_{r-1} & \ldots
\\
\ldots & (Z_{r-2},Z_{r-2}) & (Y_{r-2}, -Y_{r-2}) & 
(Z_{r-1},Z_{r-1}) & (Y_{r-1},-Y_{r-1}) &
(Y_{r}^1,Y_r^2) & \ldots
\end{array}
\label{seq1b}
\end{equation}
The strict inequality will appear when we apply
equation (\ref{FarHolder}) at the step when 
$\nn_p$ is equal to $m_{r-1}/2$ and 
$$
F_{p+1}^1 = F_{p+1}^1 (\ldots+ Z_{r-2} + Y_{r-2} + Z_{r-1} +
Y_{r-1}),
$$
$$
F_{p+1}^2 = F_{p+1}^2 (\ldots+ Z_{r-2} - Y_{r-2} + Z_{r-1} -
Y_{r-1}).
$$
Random variables $Y_{r-2},Z_{r-2}$ are independent and nondegenerate
and we can argue as in the case (\ref{seq1}) above.
All other cases in the proof of Proposition 5.6 in \cite{T-P} are
very similar and (\ref{casefar}) holds.

\subsection{Control of the close points.}\label{SecClose}

In Section \ref{SecFar} we obtained the control of the points 
$v$ far from $q_r$ and in this section we will consider 
the remaining cases when
$q_{r-1}\leq v\leq q_{r+1}$ or $-q_1\leq v< 0$ when $r=1.$
All arguments repeat the arguments of Section 5 in \cite{T-P},
so we will only consider the case when 
$$
q_{r-1}\leq v\leq q_r.
$$
As in the previous section, 
let $L_1,L_2,\ldots$ denote  constants that depend only on $\nu,\xi$ and $u.$
Consider a function
\begin{equation}
\Gamma(c)=\inf\Bigl
\{|\xi(y)-\xi(x)+(x-y)\xi'(y)| : 0\leq x,y \leq D, |x-y|\geq c
\Bigr\}.
\label{Gamma}
\end{equation}
Since $\xi''(x)>0$ we have $\Gamma(c)>0$ for $c>0.$
In the notations of Theorem \ref{Apriori} the following holds.
\begin{proposition}\label{Cprop1}
Suppose that $q_{r-1}\leq v\leq q_r.$
If $L_1\eps^{1/6}\leq 1-t_0$ then
\begin{equation}
L_1(q_r - v)\leq 1-t_0 \Rightarrow
\frac{1}{N}\e \P(\vec{n})F(A(v))\leq 2\psi(t)-\frac{(1-t_0)^2}{L_1}(v-q_r)^2
+\R,
\label{Close1}
\end{equation}
and if $L_2 \eps^{1/2}\leq (1-t_0)\Gamma((1-t_0)/L_1)$ then
\begin{equation}
L_1(q_r - v)\geq 1-t_0 \Rightarrow
\frac{1}{N}\e \P(\vec{n})F(A(v)) < 2\psi(t)
+\R.
\label{Close2}
\end{equation}
\end{proposition}
Together with a similar result for $q_r\leq v\leq q_{r+1}$
and the results of Section \ref{SecFar}, this proves Theorem
\ref{Apriori}.
We will again use Talagrand's interpolation for two copies.
Given
$$
\frac{m_{r-1}}{2}\leq m\leq m_r,
$$
let us define sequences $\vec{\nn}$ and $\vrho$ in (\ref{Tens})
and (\ref{Trho}) by
$$
0=\nn_0 = \frac{m_0}{2}, \nn_1=\frac{m_1}{2},\ldots, \nn_{r-1}=\frac{m_{r-1}}{2},
\nn_r = m, \nn_{r+1} = m_r, \ldots,\nn_{k+1}=m_k
$$
and
$$
\rho_0=q_0, \ldots, \rho_{r-1} = q_{r-1}, \rho_r=v,
\rho_{r+1} = q_r,\ldots,\rho_{k+2}=q_{k+1}.
$$
Since $\rho_r=v,$ we have $\tau=r.$
Consider a sequence  $(Y_p^1,Y_p^2)$ as in (\ref{yrand}), i.e.
$$
\e (Y_p^j)^2  = t(\xi'(\rho_{p+1}) - \xi'(\rho_{p})),
$$
$$
Y_p^1 = Y_p^2 \mbox{ for } p<r
\,\,\mbox{ and } Y_p^1, Y_p^2 \mbox{ are independent for }
p\geq r.
$$
Let $(Z_p^1,Z_p^2)$ be such that 
$$
\e (Z_p^j)^2 = (1-t)(\xi'(\rho_{p+1}) - \xi'(\rho_{p}))
$$
for $p< r-1 \mbox{ and } p> r,$
$$
\e (Z_{r-1}^j)^2 = (1-t)(\xi'(q_r) - \xi'(q_{r-1}))
$$
and $Z_r^j = 0.$ Let
$$
Z_p^1 = Z_p^2 \mbox{ for } p<r
\,\,\mbox{ and } Z_p^1, Z_p^2 \mbox{ are independent for }
p> r.
$$
If we denote
\begin{equation}
g_p^j = Y_p^j +Z_p^j 
\,\,\mbox{ for }\,\, 
0\leq p\leq k+1, j=1,2
\label{gYZ}
\end{equation}
then if follows from the construction that
$$
\e (g_p^j)^2 = \xi'(\rho_{p+1}) -\xi'(\rho_p)
\mbox{ for } p<r-1, p>r,
$$
\begin{equation}
\e (g_{r-1}^j)^2 = (\xi'(q_{r}) -\xi'(q_{r-1})) - 
t(\xi'(q_r) - \xi'(v)), 
\,\,\,\,
\e (g_{r}^j)^2 = t(\xi'(q_r) - \xi'(v))
\label{qrr}
\end{equation}
and
$$
g_p^1 = g_p^2 \mbox{ for } p<r
\,\,\mbox{ and } g_p^1, g_p^2 \mbox{ - independent for }
p> r.
$$
If we define a point $a\in[v,q_r]$ by
\begin{equation}
\xi'(a) = t\xi'(v) + (1-t)\xi'(q_r) 
\label{av}
\end{equation}
then (\ref{qrr}) can be rewritten as
$$
\e (g_{r-1}^j)^2 = \xi'(a) -\xi'(q_{r-1}),\,\,\,
\e (g_{r-1}^j)^2 = \xi'(q_{r}) -\xi'(a).
$$
If we define a new sequence $\vrho'$ by
$$
\rho_0'=q_0, \ldots, \rho_{r-1}' = q_{r-1}, \rho_r'=a,
\rho_{r+1}' = q_r,\ldots,\rho_{k+2}'=q_{k+1}
$$
obtained by inserting $a$ into the sequence $\vec{q},$
we finally get
\begin{equation}
\e (g_p^j)^2 = \xi'(\rho_{p+1}') -\xi'(\rho_p')
\mbox{ for all } p\leq k+1
\label{grho1}
\end{equation}
and
\begin{equation}
g_p^1 = g_p^2 \mbox{ for } p<r
\,\,\mbox{ and } g_p^1, g_p^2 \mbox{ are independent for }
p> r.
\label{grho2}
\end{equation}
Plugging the definition (\ref{gYZ}) into (\ref{Flg}) we get
\begin{equation}
F(\lambda,\gamma) = 
\log \int_{\Sigma\times\Sigma} 
\exp 
\Bigl(
\sum_{j\leq 2}\sigma_j \sum_{p\leq k+1}
g_{p}^j 
+\sum_{j\leq 2}\lambda (\sigma_j)^2
+\gamma\sigma_1\sigma_2
\Bigr)
d\nu(\sigma_1)d\nu(\sigma_2).
\label{Flgg}
\end{equation}
Let us define
\begin{equation}
V(\gamma,m,v)=\P(\vec{\nn})F(\lambda,\gamma)
\label{Vgmv}
\end{equation}
where we made the dependence of the right hand side
on the parameters $(\gamma,m,v)$ explicit.
In order to apply Corollary \ref{corT}, let us first note that
from the construction of sequences $\vec{\nn}$ and $\vrho$ 
we have
\begin{eqnarray*}
&&
-2t\sum_{l<\tau}\nn_{l}(\theta(\rho_{l+1}) - \theta(\rho_l)) 
-t\sum_{l\geq \tau} \nn_l(\theta(\rho_{l+1}) - \theta(\rho_l)) 
\\
&=&
-t\sum_{l\leq r-2} m_{l}(\theta(q_{l+1}) - \theta(q_l))
-t m_{r-1}(\theta(v)-\theta(q_{r-1})) 
\\
&&
-t m(\theta(q_r) - \theta(v))
-t\sum_{r\leq l\leq k} m_{l}(\theta(q_{l+1}) - \theta(q_l))
\\
&=&
-t\sum_{l\leq k} m_{l}(\theta(q_{l+1}) - \theta(q_l))
-t (m-m_{r-1})(\theta(q_r) - \theta(v)).
\end{eqnarray*}
Corollary \ref{corT} and (\ref{chione}) now imply that
\begin{eqnarray}
\frac{1}{N}\e\P(\vec{n})F(A(v)) 
&\leq&
-2\lambda u - \gamma v + V(\gamma,m,v)
\label{Ve}
\\
&&
-t\sum_{l\leq k} m_{l}(\theta(q_{l+1}) - \theta(q_l))
-t (m-m_{r-1})(\theta(q_r) - \theta(v))
+\R.
\nonumber
\end{eqnarray}
One can easily check using the argument of Lemma
\ref{Plem1} that 
\begin{equation}
V(0,m,v)=2T(m,\xi'(a))
\label{VT}
\end{equation}
where $T$ was defined in (\ref{Tmb}). 
(\ref{Tstable}) implies that
$$
V(0,m_{r-1},v) = 2 X_0(\vec{m},\vec{q},\lambda)
$$
and, therefore, (\ref{Ve}) with $\gamma=0,m=m_{r-1}$
implies that
$$
\frac{1}{N}\e\P(\vec{n})F(A(v))  \leq 2\psi(t) +\R.
$$
In order to prove Proposition \ref{Cprop1}, we will perturb
parameters $m$ and $\gamma$ around these values $0,m_{r-1}$
and use the properties of $\eps$-minimizer from the previous section.
The fundamental connection of the bound (\ref{Ve}) to the properties
of $\eps$-minimizer lies in the following fact:
\begin{equation}
\frac{\partial V}{\partial \gamma}(\gamma,m_{r-1},v)\bigr|_{\gamma=0}
= -U'(\xi'(a)).
\label{connect}
\end{equation}
The proof follows from  straightforward computation and
is given in Lemma 5.8 in \cite{T-P}. Also, similar to (\ref{fact}),
we have
\begin{equation}
\Bigl|
\frac{\partial^2 V}{\partial \gamma^2}
\Bigr|\leq L.
\label{fact2}
\end{equation}
We are now ready to prove Proposition \ref{Cprop1}.

{\bf Proof of Proposition \ref{Cprop1}.}
Let us consider a function
$$
\alpha(\gamma)=V(\gamma,m_{r-1},v) -\gamma v,
$$
which is the part of the bound (\ref{Ve}) for $m=m_{r-1}$
that depends on $\gamma.$
By (\ref{connect}) we have
$$
h(v) := \alpha'(0) = -U'(\xi'(a)) - v
=-U'\bigl(t\xi'(v) + (1-t)\xi'(q_r)\bigr) - v,
$$
since $a$ was defined in (\ref{av}). 
By (\ref{UAqr}) we have $h(q_r)=0$ and
\begin{equation}
h'(q_r)= -t\xi''(q_r)U''(\xi'(q_r))-1.
\label{hprime1}
\end{equation}
Using (\ref{f2}),
\begin{eqnarray}
h'(q_r)
&=& 
-\xi''(q_r)U''(\xi'(q_r))-1 + (1-t)\xi''(q_r)U''(\xi(q_r))
\nonumber
\\
&\leq&
\frac{L}{\xi''(q_r)}\eps^{1/6} + (1-t)\xi''(q_r)U''(\xi(q_r)).
\label{hprime2}
\end{eqnarray}
We will now show that
\begin{equation}
L\eps^{1/6}\leq 1-t_0 \Longrightarrow
h'(q_r)\leq -\frac{1-t_0}{4}.
\label{hco}
\end{equation}
If $-\xi''(q_r)U''(\xi'(q_r))\leq 1/2$ then (\ref{hprime1}) gives
$$
h'(q_r)\leq \frac{t}{2}-1\leq -\frac{1}{2}.
$$
If $-\xi''(q_r)U''(\xi'(q_r))\geq 1/2$ then
$$
\frac{1}{\xi''(q_r)} \leq -2U''(\xi'(q_r)) \leq L
$$
by (\ref{fact}), and (\ref{hprime2}) gives
$$
h'(q_r)\leq L\eps^{1/6} 
-\frac{1-t}{2}\leq -\frac{1-t_0}{4},
$$
where the last inequality holds if $4L\eps^{1/6}\leq 1-t_0.$ 
This proves (\ref{hco}).
Since (\ref{fact}) implies that $|h''(v)|\leq L$ and since
$h(q_r) = 0,$ we can write
$$
h(v)\geq (v-q_r)h'(q_r) - L(v-q_r)^2 \geq \frac{1}{8}(q_r - v)(1-t_0)
$$
if $q_r-v\leq (1-t_0)/8L$ and if (\ref{hco}) holds.
(\ref{fact2}) implies that $|\alpha''(\gamma)|\leq L$ and
we finally get
\begin{eqnarray}
\inf_{\gamma}\alpha(\gamma) 
&\leq&
\inf_{\gamma}\bigl(\alpha(0) + \alpha'(0)\gamma + L\gamma^2\bigr)
\leq \alpha(0) - \frac{\alpha'(0)^2}{L}
= \alpha(0) -\frac{h(v)^2}{L}
\nonumber
\\
&\leq& 
2X_0(\vec{m},\vec{q},\lambda) - \frac{1}{L}(1-t_0)^2(v-q_r)^2.
\label{infa}
\end{eqnarray}
Applying this to the bound (\ref{Ve}) proves (\ref{Close1}).
Note that the infimum was achieved on $\gamma=-\alpha'(0)/2L$ and
that (\ref{fact}) implies that $|\gamma|\leq L.$ 
As we explained in the remark following Corollary \ref{corT},
the bound (\ref{Tcor}) is used only for $|\gamma|\leq L.$

Next, we will prove (\ref{Close2}). If $-h(v) = U'(\xi'(a))+v \not= 0$
then we can simply use the first inequality in (\ref{infa}).
Let us assume now that $U'(\xi'(a))=-v.$ Let us set $\gamma=0$ in
the bound (\ref{Ve}) and consider the derivative of this bound
in $m$ at $m=m_{r-1},$ i.e.
\begin{eqnarray*}
D(t)
&=&
\frac{\partial V}{\partial m}(0,m,v)\bigr|_{m=m_{r-1}}
-t(\theta(q_r) - \theta(v))
%\\
%&=&
%2\frac{\partial T}{\partial m}(m,\xi'(a))\bigr|_{m=m_{r-1}}
%-t(\theta(q_r) - \theta(v))
%\mbox{\hspace{1cm} by (\ref{VT})}
\\
&=&
U(\xi'(a)) - t(\theta(q_r) - \theta(v))
\mbox{ \hspace{3.4cm} by (\ref{VT}) and (\ref{Ub}) }
\\
&=&
U(t\xi'(v)+(1-t)\xi'(q_r)) - t(\theta(q_r) - \theta(v)).
\mbox{ \hspace{0.6cm} by (\ref{av})}
\end{eqnarray*}
Since we assumed that $U'(\xi'(a))=-v,$
\begin{eqnarray*}
D'(t)
&=&
(\xi'(v)-\xi'(q_r))U'(\xi'(a)) - (\theta(q_r) - \theta(v))
\\
&=&
-(\xi'(v)-\xi'(q_r))v - (\theta(q_r) - \theta(v))
\\
&=&
\xi(q_r)-\xi(v) +(v-q_r)\xi'(q_r)\leq -\Gamma(q_r-v),
\end{eqnarray*}
where the last inequality follows from the definition (\ref{Gamma}).
By (\ref{Ubconc}), $D(t)$ is concave in $t$ and, therefore,
$$
D(1)\leq D(t) + (1-t)D'(t)\leq D(t) - (1-t)\gamma(q_r -v).
$$
By (\ref{f1}),
$$
D(1)=U(\xi'(v)) - (\theta(q_r)-\theta(v)) = 2f(v)\geq -L\eps^{1/2}.
$$
We get
$$
D(t)\geq (1-t)\gamma(q_r -v) - L\eps^{1/2}>0,
$$
if $(1-t_0)\gamma(q_r -v)> L\eps^{1/2},$ which is true
under the conditions in (\ref{Close2}) and one can finish
the proof as in (\ref{infa}).
\qed

\section{Cases reducible to the classical SK model.}\label{SecRed}

We will now show that only the case of $d<u<D$ in Theorem \ref{Th1}
is different from the classical SK model.
First of all, $d=D$ means that $\Sigma=\{-\sqrt{d},+\sqrt{d}\}$
which is precisely the case of the SK model.
If measure $\nu$ has nonzero mass at both points $\pm \sqrt{d}$ 
then $\nu(\sigma)$ is proportional to $\exp h\sigma$ for some 
external field parameter $h\in\Reals.$ 
Otherwise, if $\nu$ is concentrated at one point,
the statement of Theorem \ref{Th1} becomes trivial.
 
It remains to consider the cases of $d<D$ and $u=d$ or $u=D.$
We will only consider the case $u=d,$ since the case $u=D$ is similar.
Let us consider a set
$$
U_N(\eps)=\{\vsi : R_{1,1}\in[d,d+\eps]\}
$$
and a function
$$
F_N(\eps)=\frac{1}{N}\e\log\int_{U_N(\eps)}\exp H_N(\vsi) d\nu(\vsi).
$$
Since, by (\ref{correlation}),
$$
\e \int_{\Sigma^N}\exp H_N(\vsi)d\nu(\vsi)
<\infty,
$$
the function $\exp H_N(\vsi)$ is
$\nu$-integrable with respect to $\vsi$
almost surely and, therefore, by the monotone convergence theorem,
$$
\lim_{\eps\to 0} \int_{U_N(\eps)} 
\exp H_N(\vsi) d\nu(\vsi) = 
\int_{\{R_{1,1}=d\}} 
\exp H_N(\vsi) d\nu(\vsi) \,\,\mbox{ a.s.}
$$
Using the monotone convergence theorem once again implies
\begin{equation}
\lim_{\eps\to 0}F_N(\eps) = P_N : =
\frac{1}{N} \e\log
\int_{\{R_{1,1}=d\}} 
\exp H_N(\vsi) d\nu(\vsi).
\label{PN}
\end{equation}
If we choose a sequence $(\eps_N)$ so that
$|F_N(\eps_N) - P_N|\leq N^{-1},$ in order to prove
Theorem \ref{Th1} for $u=d$, it is enough to show that
\begin{equation}
\lim_{N\to\infty} P_N = \P(\xi,d).
\label{PNlim}
\end{equation}
We will prove this by considering two separate cases.

{\bf Case 1.} $\nu(\{\sigma^2=d\})=0.$ 
This means that the measure $\nu$ has no atoms at the points 
$\pm\sqrt{d}$ and, therefore, $\nu(\{R_{1,1}=d\}) = 0$
and $P_N=-\infty.$ To prove (\ref{PNlim}), we need to show
that $\P(\xi,d)=-\infty.$
Let us, for example, take $k=1$ and $q_1=0.$
For this choice of $\vec{m},\vec{q}$ we have 
\begin{eqnarray*}
X_0(\vec{m},\vec{q},\lambda) 
&=&
\log \e \int_{\Sigma}
\exp \Bigl(\sigma z\sqrt{\xi'(u)} 
+ \lambda\sigma^2 \Bigr)d\nu(\sigma).
\\
&=& 
\log\int_{\Sigma}
\exp \Bigl(\frac{1}{2} \sigma^2 \xi'(u) 
+ \lambda\sigma^2 \Bigr)d\nu(\sigma).
\end{eqnarray*}
Therefore,
$$
-\lambda d + X_0(\vec{m},\vec{q},\lambda) = 
\log\int_{\Sigma}
\exp \Bigl(\frac{1}{2} \sigma^2 \xi'(u) 
+ \lambda(\sigma^2 -d) \Bigr)d\nu(\sigma)
$$
and, by the monotone convergence theorem, 
$$
\lim_{\lambda\to -\infty} (-\lambda d + X_0(\vec{m},\vec{q},\lambda)) 
= \log\int_{\{\sigma^2 = d\}}\exp\Bigl(
\frac{1}{2}\sigma^2\xi'(u)
\Bigr)d\nu(\sigma) = -\infty,
$$
since we assumed that $\nu(\{\sigma^2 = d\}) = 0.$
Clearly, (\ref{Pu}) yields that
$\P(\xi,d) = -\infty.$

{\bf Case 2.} $\nu(\{\sigma^2 = d\})>0.$
This means that measure $\nu$ has at least one atom at 
the points $\pm \sqrt{d}.$ 
Consider a probability measure $\bar{\nu}$ defined by
\begin{equation}
\bar{\nu}(\{-\sqrt{d}\})=\frac{\nu(\{-\sqrt{d}\})}
{\nu(\{\sigma^2=d\})},\,\,\,
\bar{\nu}(\{\sqrt{d}\})=\frac{\nu(\{\sqrt{d}\})}
{\nu(\{\sigma^2=d\})}.
\label{barnu}
\end{equation}
Condition (\ref{diam1}) implies that $R_{1,1}=d$ if and only if
$\sigma_i^2 = d$ for all $i\leq N.$ In other words,
$$
\{R_{1,1}=d\}=\Sigma_d^N,\,\,\, \mbox{ where }\,\,\,
\Sigma_d=\{-\sqrt{d},\sqrt{d}\}.
$$
With these notations, $P_N$ in (\ref{PN}) can be written as
\begin{equation}
P_N = \log \nu(\{\sigma^2=d\})
+ \bar{P}_N
\label{barPN}
\end{equation}
where
$$
\bar{P}_N = 
\frac{1}{N} \e\log
\int_{\Sigma_d^N} 
\exp H_N(\vsi)d\bar{\nu}(\vsi).
$$
But $\bar{P}_N$ exactly falls into the case $d=D$ 
considered above because 
$\bar{\nu}(\{\sigma^2 = d\})=1,$ and, therefore, its limit
can be written as follows. If we consider
$$
\bar{Y}_{k+1}=\log 
\int_{\Sigma_d}
\exp \Bigl(\sigma\sum_{0\leq p\leq k} z_p\Bigr) 
d\bar{\nu}(\sigma),
$$
define $\bar{Y}_l$ recursively as in (\ref{Xl}) and let
$$
\bar{\P}_k(\vec{m},\vec{q})
=
\bar{Y}_0
-\frac{1}{2}
\sum_{1\leq l\leq k} m_l \bigl(\theta(q_{l+1}) - \theta(q_l)\bigr)
$$
then
$$
\lim_{N\to\infty} \bar{P}_N = \inf \bar{\P}_k(\vec{m},\vec{q}),
$$
where the infimum is over all choices of $k, \vec{m}$ and $\vec{q},$
which by (\ref{barPN}) implies
$$
\lim_{N\to\infty} P_N = 
\log \nu(\{\sigma^2 =d\}) + 
\inf \bar{\P}_k(\vec{m},\vec{q}).
$$
Equivalently, this can be written as follows. If we consider
\begin{equation}
Y_{k+1}=
\log \nu(\{\sigma^2=d\})
+ \log \int_{\Sigma_d}
\exp \Bigl(\sigma\sum_{0\leq p\leq k} z_p\Bigr) 
d\bar{\nu}(\sigma),
\label{eX}
\end{equation}
define $Y_l$ recursively as in (\ref{Xl}) and let
\begin{equation}
\P_k(\vec{m},\vec{q})
=
Y_0 -\frac{1}{2}
\sum_{1\leq l\leq k} m_l \bigl(\theta(q_{l+1}) - \theta(q_l)\bigr)
\label{paa1}
\end{equation}
then
\begin{equation}
\lim_{N\to\infty} P_N = \inf \P_k(\vec{m},\vec{q}).
\label{eP}
\end{equation}
By definition (\ref{barnu}) of measure $\bar{\nu},$
$Y_{k+1}$ in (\ref{eX}) can be also written as
\begin{equation}
Y_{k+1}=\log \int_{\{\sigma^2 = d\}}
\exp \Bigl(\sigma\sum_{0\leq p\leq k} z_p \Bigr)
d\nu(\sigma).
\label{eXX}
\end{equation}
In order to prove (\ref{PNlim}), we will show that
$\P(\xi,d)$ is equal to the right hand side of (\ref{eP}).
The definition of $\P(\xi,d)$ given by 
(\ref{Xlast}) - (\ref{Pu}) can be written equivalently
as follows. If we consider
\begin{equation}
X_{k+1}=\log 
\int_{\Sigma}
\exp\Bigl(\sigma\sum_{0\leq p\leq k} z_p 
+ \lambda (\sigma^2 - d) \Bigr)d\nu(\sigma),
\label{eXk}
\end{equation}
define $X_{l}$ recursively as in (\ref{Xl}),
and define
\begin{equation}
\P_k(\vec{m},\vec{q},\lambda,d)
=
X_0 -\frac{1}{2}
\sum_{1\leq l\leq k} m_l \bigl(\theta(q_{l+1}) - \theta(q_l)\bigr),
\label{paa2}
\end{equation}
then
\begin{equation}
\P(\xi,d)=\inf \P_k(\vec{m},\vec{q},\lambda,d),
\label{ePke}
\end{equation}
where the infimum is taken over all $\lambda, k, \vec{m}$
and $\vec{q}.$ Since $d\leq \sigma^2$ for $\sigma\in\Sigma,$
$X_{k+1}$ in (\ref{eXk}) is increasing in $\lambda$ which implies that
$X_0$ is also increasing in $\lambda.$ Therefore, for any fixed
$\vec{m}$ and $\vec{q}$ to minimize the right hand side of 
(\ref{ePke}) over $\lambda$ one should let $\lambda\to-\infty.$
By the monotone convergence theorem, almost surely,
$$
\lim_{\lambda\to-\infty} X_{k+1} = 
\int_{\{\sigma^2=d\}}
\exp \Bigl( \sigma\sum_{0\leq p\leq k} z_p \Bigr) 
d\nu(\sigma)=Y_{k+1}.
$$
Using the monotone convergence theorem repeatedly in the 
recursive construction (\ref{Xl}) gives
$\lim_{\lambda\to-\infty}X_0 = Y_0$
and comparing (\ref{paa1}) and (\ref{paa2}) we get
$$
\inf_{\lambda} \P_k(\vec{m},\vec{q},\lambda,d) 
= \P_k(\vec{m},\vec{q}).
$$
Combining this with (\ref{eP}) and (\ref{ePke}) gives
(\ref{PNlim}).
\qed


\begin{thebibliography}{99}

\bibitem{GS} Ghatak, S.K., Sherrington, D. (1977) 
Crystal field effects in a general $S$ Ising spin glass,
{\it J. Phys. C: Solid State Phys.} {\bf 10}, 3149.

\bibitem{Guerra} Guerra, F. (2003) Broken replica 
symmetry bounds in the mean field spin glass model. 
{\it Comm. Math. Phys.} {\bf 233}, no. 1, 1-12.

\bibitem{GT} Guerra, F., Toninelli, F.L. (2002)
The thermodynamic limit in mean field spin glass models,  
{\it Comm. Math. Phys.} {\bf 230}, 71-79.

\bibitem{MS} Mottishaw, P.J., Sherrington, D. (1985)
Stability of a crystal-field split spin glass. 
{\it J. Phys. C: Solid State Phys.} {\bf 18}, 5201-5213.

\bibitem{LT}
Ledoux, M., Talagrand, M. (1991)
Probability in Banach spaces. Isoperimetry and Processes.
Springer-Verlag.

\bibitem{P-PM} Panchenko, D. (2005) A question about Parisi functional.
Preprint.

\bibitem{SherK} Sherrington, D., Kirkpatrick, S. (1972)
Solvable model of a spin glass. 
{\it Phys. Rev. Lett.} {\bf 35}, 1792-1796.

\bibitem{SG} Talagrand, M. (2003)   
Spin Glasses: a Challenge for Mathematicians. 
Springer-Verlag.

\bibitem{T2} Talagrand, M. (2003) 
 On Guerra's broken replica-symmetry bound. 
{\it C. R. Math. Acad. Sci. Paris} {\bf 337}, no. 7, 477-480.

\bibitem{T-gP} Talagrand, M. (2003)
The generalized Parisi formula. 
{\it C. R. Math. Acad. Sci. Paris} {\bf 337}, no. 2, 111-114.

\bibitem{T-P} Talagrand, M. (2003)
Parisi formula. To appear in {\it Ann. Math.}

\bibitem{T-M} Talagrand, M. (2003)
On the meaning of Parisi's functional order parameter.
{\it C. R. Math. Acad. Sci. Paris} {\bf 337}, no. 9, 625-628.

\bibitem{T-Sph} Talagrand, M. (2004) Free energy of the spherical
mean field model. To appear in {\it Probab. Theory Related Fields}.

\bibitem{T-PM} Talagrand, M. (2004) Parisi measures. To appear in
{\it J. Funct. Analysis.}



\end{thebibliography}
\end{document}